\documentclass[12pt,leqno]{article}
\usepackage{amsmath,amssymb,graphicx,theorem,latexsym}
\usepackage[all]{xy}
\xyoption{2cell}
\UseAllTwocells
\usepackage[bbgreekl]{mathbbol}

\newcommand{\Appendix}[1]{%
  \refstepcounter{section}%
  \addcontentsline{toc}{section}%
    {\bfseries\appendixname~\thesection.\ #1}%
    {\medskip\noindent \Large\bfseries\appendixname\ \thesection.\ #1}%
\sectionmark{#1}\smallskip\noindent
\renewcommand{\theequation}{{\bf
{{\thesection}}.{\arabic{equation}}}}
}

\DeclareFontFamily{U}{rsf}{}
\DeclareFontShape{U}{rsf}{m}{n}{
  <5> <6> rsfs5 <7> <8> <9> rsfs7 <10->  rsfs10}{}
\DeclareMathAlphabet{\mathscr}{U}{rsf}{m}{n}
\newcommand{\mycal}[1]{\mathscr{#1}}
\newtheorem{thm}{Theorem}[section]

\newtheorem{lem}[thm]{Lemma}
\newtheorem{prop}[thm]{Proposition}
\newtheorem{defn}[thm]{Definition}
\newtheorem{claim}[thm]{Claim}
{\theorembodyfont{\rmfamily} \newtheorem{rem}{Remark}[section]}               
{\theorembodyfont{\rmfamily} \newtheorem{ex}{Example}[section]}
{\theorembodyfont{\rmfamily} }
\numberwithin{equation}{section}
\newcommand{\op}[1]{\operatorname{#1}}
\newcommand{\delq}[1]{\frac{\partial}{\partial q_{#1}}}
\newcommand{\delp}[1]{\frac{\partial}{\partial p_{#1}}}

\newcommand{\A}{\mycal{A}}

\newcommand{\Ch}{{\mathbb C}[[\hbar]]}
\newcommand{\bpi}{\boldsymbol{\Pi}}
\newcommand{\bB}{\boldsymbol{B}}
\newcommand{\be}{\boldsymbol{e}}
\newcommand{\bbf}{\boldsymbol{f}}

\newcommand{\Ah}[1]{\mycal{A}_{#1}}
\newcommand{\bbD}{{\mathbb D}}
\newcommand{\bbS}{{\mathbb{S}}}
\newcommand{\rar}[1]{\stackrel{#1}{\longrightarrow}}

\newcommand{\bLambda}{\mathbb{\Lambda}}
\newcommand{\bGamma}{\mathbb{\Gamma}}
\newcommand{\bbV}{\mathbb{V}}
\newcommand{\bbX}{\mathbb{X}}
\newcommand{\bbY}{\mathbb{Y}}

\newcommand{\zenter}{{\mathcal{Z}}}
\newcommand{\sff}[1]{{\sf{#1}}}
\newcommand{\lrmod}[2]{\mbox{${}_{#1}{\sff{Mod}}_{#2}$}}
\newcommand{\bbd}{\mathbb{d}}
\newcommand{\bbr}{\mathbb{r}}
\newcommand{\bbp}{\mathbb{p}}

\newcommand{\bGa}{\boldsymbol{\Gamma}}

\newcommand{\cEnd}{{\mycal{E}}nd}
\newcommand{\cAut}{{\mycal{A}}ut}

\newcommand{\bt}{\boldsymbol{t}}
\newcommand{\inv}{\op{inv}}
\newcommand{\bbu}{\mathbb{u}}
\newcommand{\oo}[1]{\otimes_{{}_{#1}}}
\newcommand{\AH}{{\sf{ah}}}
\newcommand{\qAH}{{\sf{qah}}}

\newcommand{\cntrct}                
{\hspace{2pt}\raisebox{1pt}{\text{$\lrcorner$}}\hspace{2pt}}
\newcommand{\cntrctother}               
{\hspace{2pt}\raisebox{1pt}{\text{$\llcorner$}}\hspace{2pt}}

 \newcommand{\les}[5]{\xymatrix{ 1 \ar[r] &  {#1}
 \ar[r] & {#2} \ar[r] & {#3} 
 \ar `r/10pt[d] `[l]  `^dl[dlll] `^dr/10pt[dll]   [dll] \\
 & {#4} \ar[r]  & {#5} \ar[r] & 1}}

\begin{document}
\title{Non-commutative tori and Fourier-Mukai duality} %
\author{O. Ben-Bassat \and  J. Block \and T. Pantev}%

\date{}

\maketitle

\begin{abstract}
The classical Fourier-Mukai duality establishes an equivalence of
categories between the derived categories of sheaves on dual complex
tori. In this article we show that this equivalence extends to an
equivalence between two dual objects. Both of these are generalized
deformations of the complex tori. In one case, a complex torus is
deformed formally in a non-commutative direction specified by a
holomorphic Poisson structure. In the other, the dual complex torus is
deformed in a B-field direction to a formal gerbe. We show these two
deformations are Fourier-Mukai equivalent.
\end{abstract}


\tableofcontents

\section{Introduction}

It is commonly believed that the most general deformations of a
complex algebraic space $X$ are captured by the deformations of some
version of the category of coherent sheaves on $X$.  Among the popular
choices are the abelian category of coherent sheaves $\op{Coh}(X)$, its
derived category $D^{b}(X)$ \cite{bondal-poisson,kontsevich-course,
kontsevich-poisson}, or a dg-enhancement of the latter. In this
paper we will look at the deformations of the derived categories and
the equivalences between them.

A particular family of infinitesimal deformations of $D^{b}(X)$ comes
from deforming the identity functor on $D^{b}(X)$. This family is
naturally parameterized by the second Hochschild cohomology
$HH^{2}(X)$ of $X$ \cite{bondal-poisson,keller-hochschild}. By
definition $HH^{i}(X)$ is the cohomology of 
$R\op{Hom}_{X\times X}(\Delta_{*}{\mathcal O}_{X},\Delta_{*}{\mathcal
O}_{X})$. If $X$ is a manifold, the geometric version of the
Hochschild-Kostant-Rosenberg theorem \cite{swan,gs,kontsevich-poisson}
identifies $HH^{i}(X)$ with the coherent cohomology of the holomorphic
polyvector fields on $X$. In particular
\begin{equation} \label{eq:decomposition}
HH^{2}(X) \cong H^{0}(X,\wedge^{2}T_{X})\oplus H^{1}(X,T_{X})\oplus
H^{2}(X,{\mathcal O}_{X}).
\end{equation}
Viewing $HH^{2}(X)$ as infinitesimal deformations of $D^{b}(X)$ we
can interpret the pieces in \eqref{eq:decomposition} as follows.
Elements in $H^{1}(X,T_{X})$ correspond to deformations of $X$ as
a complex manifold.  Elements in $H^{0}(X,\wedge^{2}T_{X})$
correspond to deforming the multiplication on ${\mathcal O}_{X}$
to a $\star$-product. Finally, elements in
$H^{2}(X,{\mathcal O}_{X})$ correspond to deforming the trivial
${\mathcal O}^{\times}$-gerbe on $X$.

Given two complex manifolds $X$ and $Y$ and an equivalence $\varphi :
D^{b}(X) \to D^{b}(Y)$ one obtains a natural  isomorphism $\widetilde{\varphi}
: HH^{2}(X) \widetilde{\to} HH^{2}(Y)$. In particular to every deformation
direction $\xi \in HH^{2}(X)$ for $D^{b}(X)$ we can associate a
deformation direction $\widetilde{\varphi}(\xi) \in HH^{2}(Y)$ for
$D^{b}(Y)$. The question we would like to investigate in general is
whether the equivalence $\widetilde{\varphi}$ deforms along with
$D^{b}(X)$ and $D^{b}(Y)$ in the directions $\xi$ and
$\widetilde{\varphi}(\xi)$.

In this paper we concentrate on the special case when $X$ is a
complex torus, $Y = X^{\vee}$ is the dual torus, and $\varphi$ is
the classical Fourier-Mukai equivalence. An interesting feature of
this case is that $\widetilde{\varphi}$ exchanges the non-commutative
deformations of $X$ with the gerby deformations of $Y$ and vise
versa. Thus the corresponding deformation of $\varphi$, if it
exists, will have to exchange sheaves of different geometric
origin.  We carry out this program to show that $\varphi$ deforms
to an equivalence of the derived category of a formal
non-commutative deformation of $X$  and the derived category of a
formal gerby deformation of $X^{\vee}$. Along the way we have to extend
some of the classical theory of complex tori to non-commutative
tori. In an attempt to make the exposition less cumbersome we have
collected the required technical results (and some generalizations) in
the appendices. 

The fact that a gerby complex torus should be Fourier-Mukai equivalent
to a non-commutative torus was originally conjectured by D.Orlov based
on the behavior of the map $\widetilde{\varphi}$ and on physical
considerations. Independently Kapustin, \cite{Kapustin}, analyzed this
setup as a duality transformation between $D$-branes in type {\sf II}
string theory. He studied branes on a complex manifold $X$ in the
presence of a $B$-field $\bB$ which is a $\overline{\partial}$-closed
$(0,2)$ form. In this case the branes admit two different
interpretations. On one hand they can be viewed as complexes of
$C^{\infty}$ complex vector bundles equipped with Hermitian
connections $D$ satisfying $F^{0,2}_{D} = \bB\cdot \op{id}$ and on the
other hand they can be viewed as complexes of coherent sheaves on the
topologically trivial holomorphic gerbe classified by $\exp(\bB) \in
H^{2}(X,\mathcal{O}^{\times})$.  On a complex torus Kapustin
investigated how such branes will transform under $T$-duality. It is
natural in this case to look for a Fourier transform that uses the
Poincar\'{e} line bundle on the product of the torus and its dual. In
order to set up such a transform for the first interpretation of
$\bB$-twisted branes, one needs to exhibit a Hermitian connection on
the Poincar\'{e} line bundle whose $(0,2)$ curvature is $-\bB$.
Kapustin searched for conditions that would ensure the existence of
such a connection. His calculations showed that if a connection
exists, then the variables on the dual torus can no longer commute. He
then gave a physical identification of the bundles on the $\bB$-twisted
torus and the bundles on the dual non-commutative torus. This led him
to conjecture that the Fourier-Mukai transform deforms to a full
equivalence of derived categories.

We are carrying out one interpretation of this conjecture. While
Kapustin works in a differential geometric context, involving vector
bundles and connections, we use the second interpretation of
$\bB$-twisted branes and hence work in an algebraic/complex geometric
context, where our tools are sheaves of algebras and modules. In order
to do this we are forced to work formally, that is, our non-commutative
torus is a formal deformation quantization of a classical complex
torus, and our gerbe is a formal deformation of the trivial
$\mathcal{O}^\times$-gerbe. The reason we are forced to work formally is
that by a theorem of Kontsevich, \cite{Kont-semi-algebraic}, the
complex torus has no algebraic, or even semi-algebraic deformations.
In order to obtain a duality statement which is not only formal but
analytic, one needs a different point of view on the category of
coherent sheaves on a torus. An appropriate  formalism
was recently developed by J.Block \cite{block-dg}. He gave an
interpretation of the categories of coherent sheaves and their
deformations as categories of twisted complexes over differential
graded algebras. For $C^{\infty}$ tori he proved a duality statement
compatible with both Kapustin's analysis and our formal duality.

As far as we know, ours is the first work on the derived
categories of modules over a deformation quantization. This is one
reason for the length of this article. We have to make sure that
much of the standard yoga of $\mathcal{O}$-modules extends to
deformations quantizations.

\subsection*{Relation to other works}

C\u{a}ld\u{a}raru \cite{caldararu} and Donagi and Pantev \cite{DP}
extended the classical Fourier Mukai equivalence to families of tori
(including some singular fibers).  Here also gerbes arise naturally. In
the case when the family has a section, the Fourier-Mukai duality
extends without sweat. When there are no sections, the natural dual
family must be interpreted as a gerbe. Our situation is much the
same. Our non-commutative torus can be viewed as a family over the
formal disc. The well-known fact that the non-commutative torus has no
points (that is no quasi-coherent sheaves supported at points)
manifests itself in the appearance of a gerbe on the dual side. In
fact, by Gabber's theorem, the support of a coherent sheaf of modules on the
non-commutative torus must be coisotropic for the complex Poisson
structure (and thus has a lower bound on the dimension of its
support). On the dual side, the support of a coherent sheaf on the gerbe must
be isotropic and thus has an upper bound on the dimension of its
support.

Polischuk and Schwarz investigated the geometry of holomorphic
structures on non-commutative real two tori
\cite{pol-sch,pol1,pol2}. They also studied the categories of sheaves
on the resulting non-commutative complex spaces. Our setup differs
from theirs in that our Poisson structures are holomorphic whereas
theirs are of type $(1,1)$. In particular in their case the derived
category does not deform in the non-commutative direction (the abelian
category of sheaves does deform however). A unifying approach to the
most general holomorphic structures on non-commutative deformations is
provided by \cite{block-dg}. 

Ours is not the first paper where non-commutative tori appear as duals
of gerbes.  In \cite{MR1,MR2} Mathai and Rosenberg, find that in some
cases families of non-commutative tori appear as the duals of families
of tori with a gerbe on the total space. Their context is completely
topological and their main result is an isomorphism of topological
$K$-theories of the two dual objects. In contrast, we work
holomorphically and our result is an equivalence of the full derived
categories of the two dual objects.

Recently, in a beautiful paper \cite{Toda}, Y.Toda proves a very
general result, related to our work. He constructs for any smooth
projective variety $X$, and a Hochschild class as above, a first order
deformation of the derived category of coherent sheaves on $X$. He
then shows that if there is an equivalence of derived categories
between $X$ and $Y$ then it deforms to the corresponding first order
deformations. It is not at all clear how to extend his results to
infinite order deformations in general. The main result of this paper
can be viewed as such an extension in the case of complex tori.

\subsection*{Notation and terminology}

\begin{description}
\item[$\Ah{X}$] a sheaf of associative flat $\Ch$ algebras on $X$
  satisfying $\A/\hbar \cong \mathcal{O}_{X}$.
\item[$\Ch$] the complete local algebra of formal power series in
  $\hbar$. 
\item[$\sff{D}^{*}$] the derived category of
  $*$-complexes of $\mathcal{O}$-modules.
  The decoration $*$ is in the set 
  $\{\varnothing,-,b\} = \{\text{unbounded},\text{bounded
  above},\text{bounded}\}$. 
\item[$\sff{D}^{*}_{c}$,$\sff{D}^{*}_{qc}$] the derived
  categories of complexes having coherent, respectively quasi-coherent
  cohomologies. 
\item[$\bbD$] the one dimensional formal disk.
\item[$\bGamma$] the Heisenberg group scheme $1 \to \mathbb{G}_{m} \to
  \bGamma \to \bLambda^{\vee} \to 0$ given by $\bB$.
\item[$\Lambda \subset V$] a free abelian subgroup of rank $2g$.
\item[$\Lambda^{\vee} \subset \overline{V}$]  the lattice of all $\xi
  \in \overline{V}^{\vee}$ satisfying $\op{Im}(\langle \xi, \lambda
  \rangle) \in \mathbb{Z}$ for all $\lambda \in \Lambda$.
\item[$\bLambda$, $\bLambda^{\vee}$] the constant group schemes
  $\Lambda\times \bbD$ and $\Lambda^{\vee}\times \bbD$ over $\bbD$.
\item[$\mathcal{P}$] the normalized Poincar\'{e} line bundle on
  $X\times X^{\vee}$.
\item[$\mathcal{P}_{\big\vert X\times \{\alpha\}}$] the degree zero
  line bundle $\mathcal{P}_{|X\times \alpha}$ on $X$.
\item[$\bpi$] a holomorphic Poisson structure on a complex manifold.
\item[$\phi_{K}^{[X\to Y]}$] the integral transform
  $\sff{D}^{*}(X) \to \sff{D}^{*}(Y)$ corresponding to a
  kernel object \linebreak $K \in \sff{D}^{*}(X\times Y)$.
\item[$\mathcal{S}_{X}$, $\mathcal{S}_{X^{\vee}}$] the Fourier-Mukai
  transforms $\sff{D}^{b}(X) \to \sff{D}^{b}(X^{\vee})$ and
  $\sff{D}^{b}(X^{\vee}) \to \sff{D}^{b}(X)$ respectively.
\item[$V$] a complex vector space of dimension $g$.
\item[$\overline{V}^{\vee}$] the complex space of conjugate linear
  homomorphisms from $V$ to $\mathbb{C}$. 
\item[$\bbV$,$\overline{\bbV}^{\vee}$] the formal spaces $V\times
  \bbD$ and $\overline{V}^{\vee}\times \bbD$.
\item[$X$] a complex torus of dimension $g$, typically $X = V/\Lambda$
\item[$X^{\vee}$] the dual complex torus of $X$, i.e. $X^{\vee} =
  \op{Pic}^{0}(X) = \overline{V}^{\vee}/\Lambda^{\vee}$.
\item[$\bbX$, $\bbX^{\op{op}}$] a formal non-commutative space $\bbX =
  (X,\Ah{X})$ and its opposite $\bbX^{\op{op}} = (X,\Ah{X}^{\op{op}})$
\item[$\bbX_{\bpi}$] the Moyal quantization of the Poisson torus
  $(X,\bpi)$.
\item[$\bbX^{\vee}$] the formal space $X^{\vee}\times \bbD$.
\end{description}

\section{The classical situation} \label{ssec:classicalfm}
First let us recall the basic properties of complex tori that we
will need. For a more detailed discussion of the properties of
complex tori the reader may consult
\cite{mumford,lange,polishchuk}.

A complex torus, is a compact complex manifold $X$ which is isomorphic
to a quotient $V/\Lambda$, where $V$ is a $g$-dimensional complex
vector space and $\Lambda \subset V$ is a free abelian subgroup of
rank $2g$. Note that by construction $X$ has a natural structure of an
analytic group induced from the addition law on the vector space $V$.

To any complex torus $X$ one can associate a dual complex torus
$X^{\vee}$. If $X$ is realized as $V/\Lambda$, the dual torus is
defined to be $X^{\vee} = \overline{V}^{\vee}/\Lambda^{\vee}$. Here
$\overline{V}^{\vee}$ denotes the space of conjugate linear
homomorphisms from $V$ to ${\mathbb C}$ and $\Lambda^{\vee} \subset
\overline{V}^{\vee}$ is the lattice defined by
\[
\Lambda^{\vee} = \left\{ \left. \xi \in \overline{V}^{\vee} \right|
\op{Im}(\xi(\lambda)) \in {\mathbb Z}, \text{ for all }
\lambda \in \Lambda \right\}.
\]
In fact the dual torus $X^{\vee}$ is intrinsically attached to $X$
and does not depend on the realization of $X$ as a quotient.
Namely, one can define $X^{\vee}$ as the torus $\op{Pic}^{0}(X)$
parameterizing all holomorphic line bundles $L \to X$ which are
invariant for the translation action of $X$ on itself.
Equivalently, these are the holomoprhic line bundles with the
property $c_{1}(L) = 0 \in H^{2}(X,{\mathbb Z})$. It is known (see
e.g. \cite{mumford,polishchuk}) that $X^{\vee}$ is a fine moduli space in
the sense that we can find a line bundle
\[
\mathcal{P} \to X\times X^{\vee}
\]
with the property that for any $\alpha \in X^{\vee}$, the
restriction $\mathcal{P}_{\big\vert X\times \{ \alpha \}}$ is
isomorphic to the degree zero line bundle corresponding to
$\alpha$. We will write $\mathcal{P}_{\alpha}$ for the line bundle
$\mathcal{P}_{\big\vert X\times \{ \alpha \}}$. Such a $\mathcal{P}$
is called a Poincar\'{e} line bundle. If we further normalize
$\mathcal{P}$ so that $\mathcal{P}_{\big\vert \{o\} \times X^{\vee}}$
is isomorphic to the trivial line bundle ${\mathcal
O}_{X^{\vee}}$, then $\mathcal{P}$ is uniquely determined.
Furthermore the assignment $\alpha \mapsto \mathcal{P}_{\alpha}$ for
the normalized Poincar\'{e} line bundle is compatible with group structures,
e.g.
\[
\mathcal{P}_{\alpha + \beta} \cong \mathcal{P}_{\alpha}\otimes
\mathcal{P}_{\beta}.
\]
Also, note that the normalized Poincar\'{e} sheaf gives rise to a
canonical isomorphism $X^{\vee\vee} \widetilde{\to} X$ \cite{lange}.

The main interest of this paper is a generalization of a powerful
duality theorem of Mukai \cite{mukai,VandenBergh}. First we need to
set things up. For a  
complex manifold $M$, let $\sff{D}^b(M)$ be the bounded derived category of
sheaves of $\mathcal{O}_{M}$-modules \cite{verdier}. An object in
$\sff{D}^b(M)$ 
is a bounded complex 
\[
\cdots\to \mycal{F}_{i-1}\to \mycal{F}_i\to \mycal{F}_{i+1}\to\cdots
\]
of analytic sheaves of $\mathcal{O}_{M}$-modules. The morphisms are
more complicated to define, see e.g. \cite{verdier,gelfand-manin}.

Given two compact complex manifolds $M$ and $N$ and an element $K\in
\sff{D}^b(M\times N)$, define the integral transform
\[
\phi_{K}^{[M\to N]} : \sff{D}^b(M)\to \sff{D}^b(N)
\]
by
\[
\phi_{K}^{[M\to N]}(G) = Rp_{N*}(p_{M}^{*}G\otimes^\mathbb{L}K).
\]
The integral transform has the following convolution property (see
\cite{mukai} or \cite[Proposition~11.1]{polishchuk}). If $M$, $N$ and
$P$ are compact 
complex manifolds and $K \in \sff{D}^{b}(M\times N)$ and $L \in
\sff{D}^{b}(N\times P)$, then one has a natural isomorphism of functors
\[
\phi_{L}^{[N\to P]} \circ \phi_{K}^{[M\to N]} \cong \phi_{K*L}^{[M\to P]},
\]
where
\[
K*L = Rp_{M\times P*}(p^{*}_{M\times N}K\otimes^{{\mathbb L}}
p^{*}_{N\times P}L) \quad \in \sff{D}^{b}(M\times P),
\]
and $p_{M\times N}$, $p_{M\times P}$, $p_{N\times P}$ are the natural
projections $M\times N\times P \to M\times N$, etc.

For a complex torus $X$, the Poincar\'{e} sheaf on $X\times X^{\vee}$
provides natural integral transforms
\[
\begin{split}
\mathcal{S}_{X} & :=
\phi_{\mathcal{P}}^{[X\to X^{\vee}]} : \sff{D}^{b}(X) \to
\sff{D}^{b}(X^{\vee}) \\
\mathcal{S}_{X^{\vee}} & :=
\phi_{\mathcal{P}}^{[X^{\vee}\to X]} : \sff{D}^{b}(X^{\vee}) \to
\sff{D}^{b}(X). 
\end{split}
\]
Now we are ready to state Mukai's duality theorem. Since we were
unable to find a reference treating this duality in the analytic
context, we have sketched below the necessary modifications of Mukai's
original proof \cite{mukai,polishchuk}:

\begin{thm} \label{thm-Mukai} The integral transform $\mathcal{S} :
  \sff{D}^{b}(X^{\vee}) \to \sff{D}^{b}(X)$ is an equivalence of triangulated
  categories.
\end{thm}
{\bf Proof.}  The theorem follows from the existence of natural isomorphisms
of functors:
\begin{equation} \label{eq-FMsquare}
\begin{split}
\mathcal{S}_{X^{\vee}}\circ \mathcal{S}_{X} & \cong (-1)^{*}_{X}[-g] \\
\mathcal{S}_{X}\circ \mathcal{S}_{X^{\vee}} & \cong (-1)^{*}_{X^{\vee}}[-g],
\end{split}
\end{equation}
where $g = \dim_{\mathbb C} X$.

The existence of \eqref{eq-FMsquare} follows from the identification $X
\cong X^{\vee\vee}$ and the fact that the canonical adjunction morphism
\begin{equation} \label{eq-adj}
\op{Id} \to \phi_{\mathcal{P}^{-1}[g]}^{[X \to X^{\vee}]}\circ
\phi_{\mathcal{P}}^{[X^{\vee}\to X]}
\end{equation}
is an isomorphism. This is a formal consequence
\cite[Theorem~11.4]{polishchuk} of the cohomology and base change
theorem for proper morphisms and of the fact that for a non-trivial
holomorphic line bundle $L \in \op{Pic}^{0}(X)$ one has
$H^{\bullet}(X,L) = 0$. The cohomology and base change theorem for
pushing forward a flat coherent sheaf under a proper analytic morphism
can be found in \cite{schneider}. To show that $H^{\bullet}(X,L) = 0$
for a non-trivial $L \in \op{Pic}^{0}(X)$ we first note that a degree
zero holomorphic line bundle on $X$ corresponds to a complex rank one
local system ${\mathbb L}$ on $X$ \cite{atiyah}. Since every complex
torus is K\"{a}hler, it follows that the natural map 
\[
H^{k}(X,{\mathbb
L}) \to H^{k}(X,L)
\]
is surjective for each $k$. If $L$ is non-trivial,
then so is ${\mathbb L}$ and so it corresponds to a non-trivial
character ${\mathbb Z}^{2g} \to {\mathbb C}^{\times}$. Now by using
the isomorphism $X \cong (S^{1})^{2g}$, the K\"{u}nneth formula, and
the vanishing of the cohomology of ${\mathbb Z}$ with
coefficients in a non-trivial character, we conclude that
$H^{k}(X,{\mathbb L}) = 0$ for all $k$.  \ \hfill $\Box$

\section{Non-commutative complex tori and $B$-fields} \label{sec:nctori}

In this section we introduce the relevant non-commutative and gerby
deformations of the two sides of the Fourier-Mukai equivalence.
The deformation of the equivalence itself will be studied in the
next section.

\subsection{Non-commutative complex tori}
Before we can extend the formalism of integral transforms to
the realm of non-commutative geometry we need to make precise the
notion of non-commutative space that we will be using. In the next
section we introduce the first  main player in our correspondence - the
deformation quantization  of a complex torus.

\

\subsubsection{Deformation quantization in the holomorphic setting}
\label{sssec-moyal}
Recall that a {\em deformation quantization}, see
 e.g. \cite{flato,NestTsygan}, 
 of a complex analytic
space $X$ is a formal one-parameter deformation of the structure sheaf
${\mathcal O}_{X}$. Explicitly this means that we are given a sheaf
$\Ah{X}$ of associative unital ${\mathbb C}[[\hbar]]$-algebras, flat
over ${\mathbb C}[[\hbar]]$, together with an algebra isomorphism
\[
\Ah{X}\oo{\Ch}  \mathbb{C} \cong \mathcal{O}_{X}.
\]
Geometrically one thinks of the data $(X,\Ah{X})$ as a non-commutative
formal deformation of $X$ over the one dimensional formal disk $\bbD
:= (\{o\},\Ch) = \op{Spf}({\mathbb C}[[\hbar]])$. In other words the data
$(X,\Ah{X})$ should be viewed as defining a non-commutative formal
space $\bbX$ which is equipped with a morphism $\bbu : \bbX \to
\bbD$ and which specializes to $X$ over the closed point $o \in \bbD$:
\[
\xymatrix@C-30pt{
X \ar[d] & \subset & \bbX \ar[d]^-{\bbu} \\
o & \in & \bbD
}
\]
Observe also that every formal non-commutative space $\bbX =
(X,\Ah{X})$ has a natural companion $\bbX^{\op{op}} :=
(X,\Ah{X}^{\op{op}})$ which has the same underlying analytic space but
is equipped with the opposite sheaf of algebras.

A {\em morphism ${\mathbb f} : {\mathbb X} \to {\mathbb Y}$ of
deformation quantizations} is defined to be a morphism between
ringed spaces $(f,f^{\sharp}) : (X,\Ah{X}) \to (Y,\Ah{Y})$, so that
$f^{\sharp} : f^{-1}\Ah{Y} \to \Ah{X} $ is continuous in the
$\hbar$-adic topology and the induced morphism
$(f,f^{\sharp}/\hbar) : (X,{\mathcal O}_{X}) \to (Y,{\mathcal O}_{Y})$
is a morphism of complex analytic spaces.  Note that the category
$\mycal{FS}/\bbD$ of
formal analytic spaces over $\bbD$ is equivalent to the
category of all commutative
deformation quantizations.

A particularly tractable class of deformation quantizations
$\bbX$ are the so called
$\star$-quan\-ti\-za\-tions (see e.g. \cite{NestTsygan}, 
\cite{BezrukavnikovKaledin},  
\cite[Definition~8.6]{Yekutieli}). 

\begin{defn} A {\em $\star$-quantization} \label{defn:star}
$(X,\Ah{X})$ of a complex manifold $X$ is a sheaf $\Ah{X}$ of
${\mathbb C}[[\hbar]]$-algebras, which is flat over ${\mathbb
    C}[[\hbar]]$ and for which:
\begin{enumerate}
\item[(a)] There is an isomorphism
$\varphi : \Ah{X}/\hbar \widetilde{\to} {\mathcal O}_{X}$
of ${\mathbb C}[[\hbar]]$-algebras.
\item[(b)] Locally on $X$ we have isomorphisms
$\psi_{U} : \Ah{X}|_{U} \widetilde{\to} {\mathcal O}_{X}|_{U}[[\hbar]]$ of
  sheaves of ${\mathbb
  C}[[\hbar]]$-modules, so that:
\begin{itemize}
\item The $\psi_{U}$'s are all compatible with $\varphi$.
\item Under this isomorphism $\psi_{U}$ the $1_{\Ah{X}}$ maps to
  $1_{{\mathcal O}_{X}}$ and the  product on
$\Ah{X}|_{U}$ becomes a product $\star$ on ${\mathcal O}_{X}|_{U}[[\hbar]]$
so that for all $a, b \in {\mathcal O}_{X}|_{U}$ we have
\[
a\star b = ab + \sum_{i = 1}^{\infty} \beta_{i}(a,b)\hbar^{i}
\]
with $\beta_{i} : {\mathcal O}_{X}|_{U} \otimes_{{\mathbb C}}
{\mathcal O}_{X}|_{U} \to {\mathcal O}_{X}|_{U}$ being bidifferential
operators.
\item The composition maps $\psi_{U} \circ {\psi_{V}}^{-1}$ are given by
a series in $\hbar$ of differential operators
${\mathcal O}_{X}|_{U \cap V} \to {\mathcal O}_{X}|_{U \cap V}$.
\end{itemize}
\end{enumerate}
\end{defn}

\

\medskip

\begin{rem}  It is a consequence of the definition that for all $a,
  b \in {\mathcal O}_{X}|_{U}[[\hbar]]$ we have $a\star b \equiv ab
  \text{ mod } \hbar$ and $\beta_{i}(1,a) = \beta_{i}(a,1) = 0$.
  Also, the term of order zero in $\hbar$ of the series of
  differential operators which give the transition isomorphisms is the
  identity map on the sheaf ${\mathcal O}_{X}|_{U \cap V}$.
\end{rem}

\

\bigskip

\noindent
A {\em morphism  between $\star$-quantizations}
${\mathbb f} : {\mathbb X} \to {\mathbb Y}$ is
defined to be a morphism between ringed
spaces $(f,f^{\sharp}) : (X,\Ah{X}) \to (Y,\Ah{Y})$ for which
$f^{\sharp} : f^{-1}\Ah{Y}  \to  \Ah{X} $ is continuous in the
$\hbar$-adic topology and such that the induced morphism
$(f,f^{\sharp}/\hbar) : (X,{\mathcal O}_{X}) \to (Y,{\mathcal O}_{Y})$
is a morphism of complex analytic spaces.  Furthermore, we
assume that $f^{\sharp}$ is {\em given by differential operators
with respect to $(f,f^{\sharp}/\hbar)$} in the following sense.

Suppose $(g, g^{\sharp}):(X,\mathcal O_{X}) \to (Y, \mathcal O_{Y})$
is a complex analytic morphism.  Then a differential operator of order
$0$ with respect to $(g, g^{\sharp})$ is defined to be a map
$g^{-1}\mathcal O_{Y} \to \mathcal O_{X}$ given locally by $a \mapsto
k \cdot (a \circ g)$ for some function $k \in {\mathcal O}_{X}$.  A
differential operator of order $j$ with respect to $(g, g^{\sharp})$
is defined inductively to be a ${\mathbb C}$-linear map $D :
g^{-1}\mathcal O_{Y} \to \mathcal O_{X}$ for which the assignment $a
\mapsto D(g^{-1}(q) \cdot a) - q \cdot D(h)$ is a differential
operator of order $j-1$.

Consider, for any structure maps
$\psi_{U} : \Ah{X}|_{U} \widetilde{\to} {\mathcal O_{X, U}}[[\hbar]]$ and
$\psi_{W} : \Ah{Y}|_{W} \widetilde{\to} {\mathcal O_{Y, W}}[[\hbar]]$ the
components of the composition
\[
\psi_{U} \circ f^{\sharp} \circ f^{-1}({\psi_{W}}^{-1})|_{f^{-1}(W) \cap U}:
f^{-1} \mathcal O_{Y, W}[[\hbar]]|_{f^{-1}(W) \cap U} \to \mathcal
O_{X, U}[[\hbar]]|_{f^{-1}(W) \cap U}
\]
\noindent
as maps from $f^{-1} \mathcal O_{Y, W}|_{f^{-1}(W) \cap U}$ to
$\mathcal O_{X, U}|_{f^{-1}(W) \cap U}$.  The
degree $0$ part is just given by $f^{\sharp}/\hbar$, the pullback
of functions, so it is a differential operator of order zero with
respect to $(f,f^{\sharp}/\hbar)$.  We say that
$\mathbb{f}$ is given by differential operators when all these components
are differential operators with respect to $(f,f^{\sharp}/\hbar)$. 

\

\medskip

\begin{rem}  \label{rem:formal} Every $\star$-quantization
  $(X,\Ah{X})$ induces a
  holomorphic Poisson structure $\bpi$ on $X$ defined by the formula
\[
 (df\wedge dg) \cntrct \, \bpi  = \frac{1}{2\hbar}\left(s(f)\star
s(g) - s(g)\star s(f) \right) \text{ mod } \hbar,
\]
for all local sections $f$ and $g$ in ${\mathcal O}_{X}$. A morphism
of $\star$-quantizations automatically induces a  Poisson morphism.

\end{rem}

\

\bigskip

\begin{ex} The basic example of a $\star$-quantization is  the
  standard Moyal product on the holomorphic functions on a complex
  vector space $V$ equipped with a constant Poisson structure $\bpi$
  \cite{moyal,flato2}. 
  By the constancy assumption, there are complex coordinates
  \[(q_1,\ldots, q_n, p_1,\ldots, p_n, c_1,\ldots, c_l)
  \] on
  $V$ so that the Poisson structure is diagonal, that is
\[
\bpi=\sum_{i=1}^{n} \delq{i}\wedge \delp{i}.
\]
Recall that a bidifferential operator on complex manifold $X$ is a
${\mathbb C}$-linear map $\varphi : {\mathcal
O}_{X}\otimes_{{\mathbb C}} {\mathcal O}_{X} \to {\mathcal O}_{X}$
which is a differential operator in each factor, i.e. for all $g
\in {\mathcal O}_{X}$ we have $\varphi(\bullet\otimes g) \in
\mycal{D}_{X}$ and $\varphi(g\otimes \bullet) \in \mycal{D}_{X}$.
Given a differential operator $D \in \mycal{D}_{X}$ we can promote
it to a bidifferential operator in two ways with $D$ acting on the
first and the second factor respectively. As usual we write
$\overleftarrow{D}$ for the bidifferential operator $D\otimes
\op{id}$ and $\overrightarrow{D}$ for the bidifferential operator
$\op{id}\otimes D$. Note that the assignment $D \to
\overrightarrow{D}$ is an algebra homomorphism, whereas the
assignment $D \to \overleftarrow{D}$ is an algebra
antihomomorphism.

With this notation we can now use $\bpi$ to define the
bidifferential operator $P$ by
\begin{equation} \label{eq:bidifferential}
P=\sum_i \left( \overleftarrow{\delq{i}}\,
\overrightarrow{\delp{i}} - \overleftarrow{\delp{i}}\,
\overrightarrow{\delq{i}}\right)
\end{equation}

Consider the sheaf $\mathcal{O}_V[[\hbar]]$ on $V$.
For any open set $U\subset V$, and any $f, g\in
\mathcal{O}_U[[\hbar]]$ we define their Moyal product
\[
f\star g=\sum_k\frac{\hbar^k}{k!} f\cdot P^k \cdot g=f\cdot\exp(\hbar
P)\cdot g = fg + \hbar\{ f, g \} + \cdots
\]
Since the $\star$-product is defined by holomorphic bidifferential
operators it maps holomorphic functions to holomorphic functions.
Moreover, since bidifferential operators are local, the product
sheafifies. We denote the resulting sheaf of $\Ch$-algebras on $V$
by $\Ah{V,\bpi}$.
\end{ex}

\subsubsection{Functional-analytic considerations} \label{sss:topologies}

For future reference we note that for a quantization, the sheaf
$\Ah{X} \to X$ is naturally a sheaf of multiplicatively convex nuclear
Frechet algebras.  Indeed, for a small enough open $U \subset X$ we
can topologize $\Ah{X}(U)$ by identifying it as a sheaf of vector
spaces with ${\mathcal O}(U)[[\hbar]]$ and using the uniform topology
on ${\mathcal O}(U)$ over compact subsets $K\subset U$ and the
$\hbar$-adic topology on ${\mathbb C}[[\hbar]]$. It is well known
\cite{treves} that both the uniform topology on holomorphic functions
and the $\hbar$-adic topology on ${\mathbb C}[[\hbar]]$ are nuclear
Frechet and so their completed tensor product is nuclear Frechet as
well. To check that a Frechet algebra is multiplicatively convex we
need in addition  to show the existence of a countable family of semi-norms
$p_{n}$, satisfying $p_{n}(a\star b) \leq p_{n}(a)p_{n}(b)$.   Choose a
exhaustion $\{ K_{n} \}_{n=0}^{\infty}$ of $U$ by nested compact
subsets. Given a local
section $f = \sum_{n} f_{n}\hbar^{n} \in \Ah{U} \cong
\mathcal{O}(U)[[\hbar]]$ we define 
\[
p_{n}(f) :=  c_{n} \sum_{i = 0}^{n} \sup_{x \in K_{i}} |f_{i}(x)|,
\]
where $c_{n} \in \mathbb{R}_{> 0}$ is an appropriately chosen
normalization constant. With this definition one checks immediately
that the $p_{n}$'s are multiplicatively convex seminorms.

\

\medskip

\noindent
{\bf Caution:}  The standard example of a nuclear Frechet algebra
which is not multiplicatively convex is the  Weyl algebra of a symplectic
vector space. This is because the Weyl relation $[x,y] = 1$ can never
hold in a Banach algebra. In the formal setting however we are saved
by the fact that the relations are of the form $[x,y] = \hbar$ and
$\hbar$ is a quasi-nilpotent element in our algebra.

\

\medskip

The multiplicative convexity property of $\Ah{U}$ is a key ingredient
in the analogue of Grauert's direct image theorem in the context of
formal deformation quantization. This theorem is
instrumental in setting up integral transforms between the coherent
derived categories.

\subsubsection{Formal Moyal products on the non-commutative torus}
\label{sec:moyal}
Let $X=V/\Lambda$ be a complex torus, with a holomorphic Poisson
structure $\bpi$. Since the holomorphic tangent bundle of a
complex torus is trivial, the bitensor $\bpi \in
H^{0}(X,\wedge^{2}T_{X})$ will necessarily be translation
invariant and hence will be of constant rank on $X$. The formal
$\star$-quantizations of a complex manifold equipped with a
Poisson structure of constant rank are known to be parameterized
\cite{NestTsygan,BezrukavnikovKaledin,Yekutieli} by an affine
space. This affine space is modeled on $F^{2}[[\hbar]]$, where $F^{2}$
is the second step of the
Hodge filtration on the second de Rham cohomology of the
symplectic Lie algebroid given by the sheaf of holomorphic vector
fields tangent to the leaves of the Poisson foliation. In the case
of a Poisson complex torus $(X,\bpi)$ the picture simplifies since
one can use the Moyal product to construct a canonical point in
the moduli space of quantizations of $(X,\bpi)$. We will call this
point {\em the Moyal quantization} of $(X,\bpi)$. The Moyal
quantization is very concrete and easier to work with than the
general constructions found for example in
\cite{Kontsevich,NestTsygan,BezrukavnikovKaledin,Yekutieli}. Since
all the essential features of the Fourier-Mukai duality are
already present in the context of Moyal deformations, we chose to
work mainly in this context.

To define the Moyal quantization $(X,\Ah{X,\bpi})$ of a holomorphic
Poisson torus $(X,\bpi)$ we will use the realization of $X$ as a
quotient $X = V/\Lambda$.  Let $\pi:V\to X$ be the covering
projection. Define the sheaf $\Ah{X,\bpi}$ of $\Ch$-algebras on $X$ as
follows. As a sheaf of $\mathbb{C}_{X}[[\hbar]]$-modules  it will
be just $\mathcal{O}_{X}[[\hbar]]$. To put a $\star$-product on
this sheaf one only has to use the natural identification
$\mathcal{O}_{X}[[\hbar]] := (\pi_{*}\mathcal{O}_{V})^{\Lambda}$
and note that the $\bpi$-Moyal product on $V$ is translation
invariant by construction. Explicitly the sections of $\Ah{X,\bpi}$
over $U\subset X$ can be described as the invariant sections
\begin{equation}\label{inv}
\Ah{X,\bpi}(U)=\Ah{V,\bpi}(\pi^{-1}(U))^\Lambda
\end{equation}
on the universal cover $V$. This is
well-defined since the Poisson structure $\bpi$ is constant and thus the
operator $P$ is translation invariant.

\subsubsection{The group structure on non-commutative tori}
\label{sssec:groupstr} 
For understanding the convolution of sheaves on $\bbX_{\bpi}$ it will
be useful to have a lift of the group structure on $X$ to a group law
on $\bbX_{\bpi}$. In contrast with the commutative case, we can not
hope for the multiplication to live on a single non-commutative
torus. This is because the multiplication on the commutative torus is
not a Poisson map. However, this problem can be easily rectified if we
replace the torus by $X$ by $X\times \mathbb{Z}$ equipped with the
Poisson structure which on the component $X\times \{k\}$ is $k\bpi$.

In this approach we view $\bbX_{\bpi}$ as a connected component of
a non-commutative space
\[
\aleph_{\bpi}=\coprod_{k\in\mathbb{Z}}\bbX_{k\bpi}
\]
This is a deformation quantization of the complex Poisson manifold
\[
\left(X\times \mathbb{Z}, \coprod_{k\in \mathbb{Z}}k\bpi\right)
\]
The group structure on the space $\aleph_{\bpi}$ is given by a map
$\mathbb{m}: \aleph_{\bpi}\times_{\bbD} \aleph_{\bpi}\to
\aleph_{\bpi}$. Viewed as a map of ringed spaces $\mathbb{m}$ is a
pair $\mathbb{m}=(m, m^{\sharp})$ where $m$ is the product group law
on $X\times \mathbb{Z}$ and $m^{\sharp}=\{
m^{\sharp}_{(a,b)}\}_{(a,b)\in \mathbb{Z}\times\mathbb{Z}}$ where
\[
m^{\sharp}_{(a,b)}:m_{(a,b)}^{-1}\Ah{(a+b)\bpi}\to
p_{1}^{-1}\Ah{a\bpi}\widehat{\otimes}_{{}_{\Ch}} p_{2}^{-1}\Ah{b\bpi},
\]
where  $m_{(a,b)}$ is the natural addition map from $(X\times
\{a\})\times (X\times \{b\})\to (X\times \{a+b\})$. For future
reference we will write $\mathbb{m}_{(a,b)}$ for the component
$(m_{(a,b)}, m_{(a,b)}^{\sharp})$. To define the
group structure, we have to define the map $m^{\sharp}$ and verify
that it gives a Hopf algebra structure on the structure sheaf of
$\aleph_{\bpi}$. Since the structure sheaf of $\aleph_{\bpi}$
descends from the cover $V\times \mathbb{Z}$ it suffices to define
the map $m^{\sharp}$ there.

Recall that on $V$ the sheaf $\Ah{k\bpi}$ is simply the sheaf
$\mathcal{O}_V[[\hbar]]$ equipped with the Moyal product
$\star_{k\bpi}$. Let $f\in \Ah{(a+b)\bpi}$, that is, $f$ is a
locally defined holomorphic function on $V$ with values in $\Ch$.
Now define
\[
m^{\sharp}_{(a,b)}(m^{-1}f)(v_1,v_2)=f(v_1+v_2).
\]
Here $m^{\sharp}$ is thought of as a map
$m^{\sharp}_{(a,b)}:m_{(a,b)}^{-1}\mathcal{O}_V[[\hbar]]\to
p_{1}^{-1}\mathcal{O}_V[[\hbar]]\widehat{\otimes}_{{}_{\Ch}}
p_{2}^{-1}\mathcal{O}_V[[\hbar]]$, and we use the fact that our
completed tensor product
$p_{1}^{-1}\mathcal{O}_V[[\hbar]]\widehat{\otimes}_{{}_{\Ch}}
p_{2}^{-1}\mathcal{O}_V[[\hbar]]$ is naturally identified with
$\mathcal{O}_{V\times V}[[\hbar]]$. Note that this product
descends to the torus since the addition on $V$ intertwines with
the covering actions. Note that with this definition $m^{\sharp}$ is a
coassociative and cocommutative coproduct.  We now verify that 
$m^{\sharp}$ is a map of algebras.

\begin{prop} \label{prop:Delta} The coproduct
\[
m^{\sharp} : m^{-1}\Ah{\aleph_{\bpi}} \to p_{1}^{-1}\Ah{\aleph_{\bpi}}
\widehat{\otimes}_{{}_{\Ch}} p_{2}^{-1} \Ah{\aleph_{\bpi}}
\]
is a morphism of sheaves of rings.
\end{prop}

\noindent
{\bf Proof.}  To check this, we need that for any
$f=\{f_a\}$ and $g=\{g_b\}$ with $f_a,g_b\in
\mathcal{O}_V[[\hbar]]$ we have
\begin{equation}\label{Hopf}
m^{\sharp}(m^{-1}f\star
m^{-1}g)_{(a,b)}(v_1,v_2)=(f_{a+b}\star_{(a+b)\bpi}
g_{a+b})(v_1+v_2)
\end{equation}
To check the property \eqref{Hopf} we will use the fact that the
$\star$-products on the different components of $\aleph_{\bpi}$
are all Moyal products built out of Poisson structures that are
proportional to $\bpi$. Since $\bpi$ is a constant Poisson
structure we can choose a system $\{ p_{1}, \ldots, p_{n}, q_{1},
\ldots q_{n}, \linebreak c_{1}, \ldots, c_{l} \}$ of linear coordinates on
the vector space $V$ so that
\[
\bpi = \sum_{i = 1}^{n} \frac{\partial}{\partial p_{i}} \wedge
\frac{\partial}{\partial q_{i}}.
\]
Let
\[
P = \sum_{i = 1}^{n} \left(
\overleftarrow{\delq{i}}\, \overrightarrow{\delp{i}} -
\overleftarrow{\delp{i}}\, \overrightarrow{\delq{i}}\right)
\]
The Moyal product $\star_{k\bpi}$ is given by
\[
f\star_{k\bpi} g = f \exp(k\hbar P) g
\]
for any two local sections $f$ and $g$ in ${\mathcal
O}_{V}[[\hbar]]$. Now, viewing $m^{\sharp}(f)$ and $m^{\sharp}(g)$ as
local sections in $\mathcal{O}_{V}[[\hbar]]$ we have
$m^{\sharp}(f)(v_{1},v_{2}) = f(v_{1}+v_{2})$ and
$m^{\sharp}(g)(v_{1},v_{2}) = g(v_{1}+v_{2})$.  In these terms
\eqref{Hopf} becomes
\[
m^{\sharp}\left( f \exp( (a+b)\hbar P ) g\right) = m^{\sharp}(f)\exp (\hbar
(aP\otimes 1 + 1\otimes bP))m^{\sharp}(g).
\]
This relation is a simple consequence of the chain rule as follows.

Consider a bi-differential operator $L$ with constant coefficients 
on a vector space $V$.  Denote the addition map on $V$ by 
$a: V \times V \to V$.  Applying the chain rule to  differentiation 
of the functions 
$m^{\sharp} (f) = f \circ a$, and $m^{\sharp} (g) = g \circ a$ with respect to 
coordinates on the first or second copy of $V$ we get
\[
m^{\sharp} (f) \cdot (L \otimes 1) \cdot m^{\sharp}(g) = 
m^{\sharp} (f \cdot L \cdot g)  
= m^{\sharp} (f) \cdot (1\otimes L) \cdot m^{\sharp}(g).
\]
Thus, we have established a group structure on $\aleph_{\bpi}$,  which
completes the proof of the proposition. \ \hfill $\Box$

\

\bigskip

\noindent
The usual definition of the antipode map of a (non necessarily
commutative) Hopf algebra allows us to define the inversion on the
group space $\aleph_{\bpi}$.  It is given by a morphism of ringed
spaces $\aleph_{\bpi} \to \aleph_{\bpi}^{\op{op}}$ represented by a
pair $\mathbb{inv} = (\inv, \inv^{\sharp})$.  Here $\inv : X\times
\mathbb{Z} \to X\times \mathbb{Z}$ is the group inversion and the
antipode $\inv^{\sharp} = \{\inv^{\sharp}_{a} \}_{a \in \mathbb{Z}}$
where $\inv_{a}$ is the algebra isomorphism
\[
\inv_{a} : \inv^{-1}\Ah{a\bpi} \to \Ah{-a\bpi}^{\op{op}}. 
\]
Similarly to the definition of $m^{\sharp}$, it suffices to define
$\inv^{\sharp}_{a}$ on $V$. Again we identify the sheaves $\Ah{a\bpi}$ and
$\Ah{-a\bpi}^{\op{op}}$ with $\mathcal{O}[[\hbar]]$ and define 
\[
\inv^{\sharp}_{a}(\inv^{-1}f)(v) := f(-v)
\]
for any local section $f \in \Ah{a\bpi}$ viewed as a $\Ch$-valued
locally defined holomorphic function on $V$. The same reasoning as in
the proof of Proposition~\ref{prop:Delta} now implies that
$\inv^{\sharp}_{a}$ is a ring homomorphism.

This concludes our general discussion on the properties of the
formal non-commutative torus $\bbX_{\bpi}$.
The Moyal quantization $\bbX_{\bpi} := (X,\Ah{X,\bpi})$ of the torus $X$ is
the first geometric input in the non-commutative Fourier-Mukai duality.
As we mentioned above, it is instructive to view the ringed space
$\bbX_{\bpi}$ as a formal deformation of the complex manifold $X$:
\[
\xymatrix@C-30pt{
X \ar[d] & \subset & \bbX_{\bpi} \ar[d]^-{\bbu_{\bpi}} \\
o & \in & {\mathbb D}
}
\]
which is parameterized
by the one dimensional formal disc $\mathbb{D} = \op{Spf}({\mathbb
C}[[\hbar]])$. Now this suggests that the non-commutative duality we
seek should be thought of as a formal deformation of the usual Mukai
equivalence (see Theorem~\ref{thm-Mukai}) between the derived
categories of the dual complex tori $X$ and $X^{\vee}$. Thus we need
to identify a dual object for $\bbX_{\bpi}$ which is again defined over
$\mathbb{D}$ and which specializes to $X^{\vee}$ at the closed point
$o \in {\mathbb D}$.

\subsection{Gerby complex tori}

The first clue for what the dual object should be, comes from the
fact that to first order in the formal parameter $\hbar$ this dual
object should again be determined by the Poisson structure $\bpi$.
It turns out that the correct dual object is an ${\mathcal
O}^{\times}$-gerbe on the formal space $X^{\vee} \times {\mathbb
D}$ which restricts to the trivial ${\mathcal O}^{\times}$-gerbe
on the reduced space $X^{\vee}\times \{ o \}$.

On the infinitesimal level this can be motivated as follows. As
explained in section~\ref{sssec-moyal}, the tangent space to
the moduli of $\star$-deformations of $X$ is
$H^{0}(X,\wedge^{2}T_{X})$. However, since $X$ is a complex torus, its
holomorphic tangent bundle is trivial and so we have an identification
\begin{equation} \label{eq:ncdefos}
H^{0}(X,\wedge^{2}T_{X}) = \wedge^{2}T_{X,0} = \wedge^{2}V.
\end{equation}
On the other hand, recall \cite{Griffiths-Harris} that if $Y$ is a
complex torus with universal cover $W$, then the Dolbeault
cohomology group $H^{p,q}_{\bar{\partial}}(Y) =
H^{q}(Y,\Omega^{p}_{Y})$ can be naturally identified with the
vector space $\wedge^{p}W^{\vee}\otimes
\wedge^{q}\overline{W}^{\vee}$. In particular the cohomology space
$H^{0,2}_{\bar{\partial}}(Y) = H^{2}(Y,{\mathcal O}_{Y})$ is
naturally identified with $\wedge^{2}\overline{W}^{\vee}$.
Applying this comment to the torus $X^{\vee} =
\overline{V}^{\vee}/\Lambda^{\vee}$ gives an identification
\begin{equation} \label{eq:gerbydefos}
H^{2}(X^{\vee},{\mathcal O}) = \wedge^{2} V.
\end{equation}
Thus we get an identification of the tangent space to the moduli of
$\star$-deformations of $X$ and the space $H^{2}(X^{\vee},{\mathcal
  O})$, which in turn can be viewed as the tangent space to
deformations of $X^{\vee}$ as an ${\mathcal O}^{\times}$-gerbe.
One can check that the identification
\begin{equation} \label{eq:match}
H^{0}(X,\wedge^{2}T_{X}) = H^{2}(X^{\vee},{\mathcal O})
\end{equation}
coming from \eqref{eq:ncdefos} and \eqref{eq:gerbydefos} is precisely
the identification between the pieces of the Hochschild cohomology of
$X$ and $X^{\vee}$ given by the cohomological Fourier-Mukai transform
$\widetilde{\mathcal{S}}_{X}$ that we discussed in the
introduction. To see this one simply has to note that on the level of
cohomology of polyvector fields the map $\widetilde{\mathcal{S}}_{X}$
is given by the cohomological Fourier-Mukai transform $\alpha \mapsto
p_{X^{\vee}*}(\exp(c_{1}(\mathcal{P}))\cup p_{X}^{*}\alpha)$.

Let $\bB \in H^{2}(X^{\vee},{\mathcal O})$ be the element
corresponding to $\bpi \in H^{0}(X,\wedge^{2}T_{X})$ via the
isomorphism \eqref{eq:match}. The class $\bB$ determines an ${\mathcal
O}^{\times}$-gerbe $_{\bB}\bbX^{\vee}$ over $\bbX^{\vee} :=
X^{\vee}\times \bbD$. In
fact this gerbe can be defined explicitly as a quotient gerbe. To
streamline the discussion we introduce special notation for the formal
analytic spaces we need. We will write
\[
\begin{split}
\overline{\mathbb V}^{\vee} &
:= \overline{V}^{\vee} \times \bbD \\
{\mathbb X}^{\vee} & :=
X^{\vee} \times \bbD
\end{split}
\]
for the formal analytic spaces which are
constant bundles over $\bbD$ with fibers $\overline{V}^{\vee}$ and
$X^{\vee}$ respectively. We will also write
\[
\bLambda^{\vee} :=
\Lambda^{\vee} \times \bbD
\]
for the constant group space over $\bbD$
with fiber $\Lambda^{\vee}$.

We can think of the formal
analytic space $\bbX^{\vee}$ as a quotient of the formal
analytic space $\overline{\bbV}^{\vee}$ by the free action
of the group space $\Lambda^{\vee}$, with $\Lambda^{\vee}$ acting
trivially on $\bbD$. Similarly, viewing $\bbX^{\vee}\to \bbD$ as a
relative space over $\bbD$ we can realize it  as the quotient of the
relative space
$\overline{\bbV}^{\vee} \to \bbD$ by the relative action of the
trivial bundle of commutative groups ${\bLambda}^{\vee}
 \to \bbD$. For the construction of
$_{\bB}\bbX^{\vee}$ we first define a bundle ${\bGamma} \to \bbD$ of
 non-commutative
groups on $\bbD$, which is a Heisenberg
extension:
\begin{equation}  \label{eq:heisenberg}
1\to {\mathbb G}_{m} \to {\bGamma}\to
{\bLambda}^\vee\to 1,
\end{equation}
where ${\mathbb G}_{m}$ is the multiplicative group scheme over
$\bbD$. As a formal space ${\bGamma} = {\mathbb
G}_{m}\times_{\bbD} {\bLambda}^{\vee}$.  Given a formal
space $S \to \bbD$ and sections $\xi, \xi' \in
{\bLambda}^{\vee}(S)$, $z, z' \in {\mathbb G}_{m}(S)$ the
multiplication on ${\bGamma}$ is given by the formula
\begin{equation} \label{eq:gammalaw}
(\xi,z)\cdot(\xi',z')=
(\xi+\xi',zz'c(\xi, \xi')).
\end{equation}
Where we define $c(\xi, \xi')$ by 
\begin{equation} \label{eq:defofc}
c(\xi, \xi') = \exp{\left(\hbar \pi^{2}
\bB(\xi', \xi)\right)}
\end{equation}
In this formula $\bB$ is interpreted as a group cocycle of
$\Lambda^{\vee}$ with values in ${\mathbb C}$. This involves two
steps. First  we use the
canonical splitting 
\[
\xymatrix@1{H^{2}(X^{\vee},\mathbb{C}) \ar[r]
  & H^{2}(X^{\vee},\mathcal{O})\ar@/^1pc/[l]  }
\] 
of the  Hodge filtration
on $H^{2}(X^{\vee},\mathbb{C})$ to interpret $\bB \in
H^{2}(X^{\vee},\mathcal{O})$ as an element in
$H^{2}(X^{\vee},\mathbb{C})$, and then we us the fact that $X^{\vee}$
is a $K(\Lambda^{\vee},1)$ to identify $H^{2}(X^{\vee},\mathbb{C})$
with the group cohomology $H^{2}(\Lambda^{\vee},\mathbb{C})$. 
Explicitly, in these terms,  $\bB$ is viewed as a skew-symmetric
biadditive map
\begin{equation} \label{eq:cocycleB}
\bB : \Lambda^{\vee} \times \Lambda^{\vee} \to \mathbb{C}
\end{equation}
which can be defined explicitly as follows. Recall that 
\[
\overline{V}^{\vee} = \left\{ l : V \to \mathbb{C} \left|
\begin{minipage}[b]{1.8in} 
\[
\begin{split}
l(v_{1} + v_{2})  & = l(v_{1}) + l(v_{2}) \\
l(c\cdot v) & = \bar{c}\cdot l(v)
\end{split}
\]
\end{minipage}
\right.
\right\},
\]
and
\[
\Lambda^{\vee} = \left\{ \left. \xi \in \overline{V}^{\vee} \right|
\op{Im}(\xi(\lambda)) \in \mathbb{Z}, \text{ for all } \lambda \in
\Lambda \right\}.
\]
Now, to every $l \in \overline{V}^{\vee}$ we can associate a
natural complex linear map $\bar{l} : V^{\vee} \to \mathbb{C}$, given
by $\bar{l}(v) := \overline{l(v)}$, and the alternating map
\eqref{eq:cocycleB} is given explicitly by
\[
\bB(\xi_{1},\xi_{2}) = \bpi \cntrct (\bar{\xi}_{1}\wedge \bar{\xi}_{2}).
\]

The group space ${\bGamma}$ still acts on $\overline{\bbV}^{\vee}$
by its image in ${\bLambda}^\vee$. Every section $x : S \to
\overline{\bbV}^{\vee}$ has a stabilizer equal to ${\mathbb
G}_{m}(S)$. The quotient $[\overline{\bbV}^{\vee}/{\bGamma}]$ is
therefore a ${\mathcal O}^{\times}$-gerbe. We will denote this
gerbe by ${}_{\bB}\bbu^{\vee} : {}_{\bB}\bbX^{\vee}\to \bbD$. Since
$_{\bB}\bbX^{\vee}$ is 
constructed as a quotient, we can compute the classifying element
of $_{\bB}\bbX^{\vee}$ in $H^{2}(\bbX^{\vee}, {\mathcal
O}^{\times})$ as the image of the ${\bLambda}^{\vee}$-torsor
$\overline{\bbV}^{\vee} \in
H^{1}(\bbX^{\vee},{\bLambda}^{\vee})$ under the boundary map
$H^{1}(\bbX^{\vee}, {\bLambda}^{\vee}) \to
H^{2}(\bbX^{\vee},{\mathcal O}^{\times})$ associated with
\eqref{eq:heisenberg}. From the definition of
\eqref{eq:heisenberg} it follows that $_{\bB}\bbX^{\vee}$ is
classified by $c \in H^{2}(\bbX^{\vee},{\mathcal
O}^{\times})$. More simply, if we ignore the map to $\bbD$, we can
think of the gerbe $_{\bB}\bbX^{\vee}$ as the quotient of the
formal space $\overline{\bbV}^{\vee}$ by the group $\Gamma =
H^{0}(\bbD,\bGamma)$. Explicitly $\Gamma$ is given as the central
extension
\[
1 \to {\mathbb C}[[\hbar]]^{\times} \to \Gamma \to  \Lambda^{\vee}
\to 0
\]
classified again by $c$ which is
now viewed as an 
element in the group cohomology  \linebreak
$H^{2}(\Lambda^{\vee},{\mathbb C}[[\hbar]]^{\times})$.

In section~\ref{s:Pic}  we will show that the stack $_{\bB}\bbX^{\vee}$
can be identified with the relative Picard stack
$\mycal{P}ic^{0}(\bbX_{\bpi}/\bbD)$ of the formal non-commutative
space $\bbX_{\bpi} \to \bbD$.

\section{Non-commutative line bundles and their moduli} \label{sec:Poincare}

In this section we investigate the family of translation
invariant line bundles on the non-commutative
torus $\bbX_{\bpi}$. In particular we exhibit a complete
(Poincar\'{e}) family of
such line bundles parameterized by the stack  ${}_{\bB}\bbX^{\vee}$.

\subsection{Line bundles and factors of automorphy} 

Recall the classical picture for line bundles on complex tori in terms
of factors of automorphy. We will describe this in a sufficiently
general context that it applies to all the situations we need. Suppose
$W$ is a locally compact space on which a discrete (not necessarily
commutative) group $\Upsilon$ acts freely and properly discontinuously
by homeomorphisms, (on the left). Denote by $Y$ the quotient and let
$\tau : W\to Y$ be the covering projection. Let $\Ah{W}$ be a sheaf of
unital not necessarily commutative algebras on $W$ which is
equivariant with respect to the action of $\Upsilon$. This means that
for every $\upsilon \in \Upsilon$ we are given an isomorphism
$\mathfrak{a}_{\upsilon} : \Ah{W} \to \upsilon^{*}\Ah{W}$ of sheaf of
algebras on $W$ and that these isomorphisms satisfy
$(\upsilon')^{*}(\mathfrak{a}_{\upsilon})\circ
\mathfrak{a}_{\upsilon'} = \mathfrak{a}_{\upsilon\upsilon'}$. Our
convention is that $\Upsilon$ acts on sections on the right. For any
open set $U \subset W$ and any 
$\upsilon \in \Upsilon$ we will write
\[
\xymatrix@R-2pc{
\Ah{W}(\upsilon(U)) \ar[r]^-{\cong} & \Ah{W}(U) \\
s \ar@{|->}[r] & s\cdot \upsilon
}
\]
for the action of $\upsilon$  on $\Ah{W}(\upsilon(U))$, i.e. $s\cdot
\upsilon := \upsilon^{*}(\mathfrak{a}_{\upsilon^{-1}})(s)$.

We will denote the product in $\Ah{W}$
by $a\star b$. The sheaf $\Ah{W}$ descends to a sheaf $\Ah{Y}$ of
algebras on $Y$ defined
by
\[
\Ah{Y} := (\tau_{*} \Ah{W})^{\Upsilon}.
\]
Explicitly, given an open set $U \subset Y$ we have
\[
\Gamma(U,\Ah{Y}) = \left\{ \left. s \in \Gamma(\tau^{-1}(U), \Ah{W})
\right| \,\, s\cdot \upsilon = s  \right\}.
\]
We are interested in describing sheaves of left (respectively right)
$\Ah{Y}$ modules that are locally free of rank one. In the usual way,
the isomorphism classes of such modules are described by elements in
the non-abelian cohomology set $H^1(Y,\Ah{Y}^\times)$. Note that a
\v{C}ech cocycle $[\gamma] \in Z^{1}({\mathfrak U},\Ah{Y}^\times))$ for
some open covering ${\mathfrak U}$ of $Y$ gives a left (respectively
right) rank one $\Ah{Y}$-module if we let $\gamma$ multiply elements
in $\Ah{Y}$ on the right (respectively left).

Alternatively we can describe sheaves of locally free rank one left
(respectively right) $\Ah{Y}$ modules in terms of factors of
automorphy.  The non-commutative factors of automorphy associated with
$\Upsilon$ and $W$ are just degree one group cocycles of $\Upsilon$
with values in the global sections of the sheaf $\Ah{W}^{\times}$.
Given such a cocycle $\be \in Z^{1}(\Upsilon,
H^{0}(W,\Ah{W}^{\times}))$ we can define a sheaf $L(\be)$ of left
$\Ah{W}$ modules of rank one (respectively a sheaf $R(\be)$ of right
$\Ah{W}$ modules of rank one) as follows. By definition a
non-commutative factor of automorphy $\be$ is a map $\be : \Upsilon \to
H^{0}(W,\Ah{W}^{\times}))$ satisfying the (left) cocycle condition
\begin{equation}\label{cocyclecondition}
\be(\upsilon_1\cdot\upsilon_2)=\be(\upsilon_2)\star
(\be(\upsilon_1)\cdot\upsilon_2)
\end{equation}
Now define the sheaves $L(\be)$ and $R(\be)$ by
\begin{equation} \label{eq:leftright}
\begin{split}
\Gamma(U,L(\be)) & = \left\{\left.  s\in
\Gamma(\tau^{-1}(U),\Ah{W})\,\, \right| \, s\cdot
\upsilon = s \star \be(\upsilon)\right\} \\
\Gamma(U,R(\be)) & = \left\{ \left. s\in\Gamma(\tau^{-1}(U),\Ah{W})\,\,
\right| \, s\cdot
\upsilon = (\be(\upsilon^{-1})\cdot \upsilon)\star s \right\}
\end{split}
\end{equation}
for an open $U\subset Y$. The fact that these formulas define sheaves
easily follows from the cocycle condition
\eqref{cocyclecondition}. Clearly $L(\be)$ and $R(\be)$ have
respectively
left and right $\Ah{Y}$-module structure since in the definition
\eqref{eq:leftright} the
non-commutative factor of automorphy $\be$ multiplies on the right and
left respectively. Modulo the obvious equivalences the assignment $\be
\mapsto L(\be)$ (respectively $\be \mapsto R(\be)$) gives the well
known map
\[
H^1(\Upsilon, H^{0}(W,\Ah{W}^\times))\to H^1(Y,\Ah{Y}^\times)
\]
from cohomology classes of factors of automorphy to isomorphisms
classes of rank one locally free left (respectively right) $\A_Y$-modules.

\begin{rem} \label{rem:two_actions} {\bf (i)} It is instructive to point out
  that the sheaves $L(\be)$ and $R(\be)$ can be written as the
  invariants of appropriately defined actions of $\Upsilon$ on
  $\tau_{*}\Ah{W}$. More precisely, given a non-commutative factor of
  automorphy $\be \in Z^{1}(\Upsilon,H^{0}(W,\Ah{W}))$ we can define
  two new $\Upsilon$-equivariant structures on the sheaf $\Ah{W}$ by
  the formulas
\[
\begin{split}
s{\scriptscriptstyle \lozenge} \upsilon& := (s\star
\be(\upsilon^{-1}))\cdot \upsilon, \\ s{\scriptscriptstyle \triangle}
\upsilon & := \be(\upsilon)\star (s\cdot \upsilon),
\end{split}
\]
for all $\upsilon \in \Upsilon$ and all sections $s$ of $\Ah{W}$ over
$\Upsilon$-invariant open sets on $W$. Now the automorphicity
conditions in \eqref{eq:leftright} become simply the condition of
invariance w.r.t. these actions and so can describe $L(\be)$ and
$R(\be)$ as the sheaves $(\tau_{*}\Ah{W})^{\Upsilon,
{\scriptscriptstyle \lozenge}}$ and $(\tau_{*}\Ah{W})^{\Upsilon,
{\scriptscriptstyle \triangle}}$ respectively.

\smallskip

\noindent
{\bf (ii)} The somewhat mysterious formulae defining the
automorphicity condition for $R(\be)$ and the ${\scriptscriptstyle
\triangle}$-action are forced on us by the non-commutative nature of
$\Ah{W}$. Indeed, the fact that the sheaf of groups $\Ah{W}^{\times}$
is in general non-commutative  implies that there are two natural
notions of a factor of automorphy. First we have the {\em left}
factors of automorphy $\be : \Upsilon \to H^{0}(W,\Ah{W}^{\times})$
satisfying the left cocycle condition \eqref{cocyclecondition}. By
the same token we have the {\em right} factors of automorphy $\bbf :
\Upsilon \to H^{0}(W,\Ah{W}^{\times})$ satisfying the right
cocycle condition
\[
\bbf(\upsilon_{1}\upsilon_{2}) = (\bbf(\upsilon_{1})\cdot
\upsilon_{2})\star \bbf(\upsilon_{2}).
\]
Clearly, given a right factor of automorphy we can write a right
automorphicity condition on $s$, namely $s\cdot \upsilon =
\bf(\upsilon)\star s$, and a new action of $\Upsilon$ on
$\Ah{W}$. However, the assignment $\be(\upsilon) \mapsto \bbf(\upsilon)
:= \be(\upsilon^{-1})\cdot \upsilon$ transforms bijectively left
factors of automorphy to right ones. Plugging this into the
formulas defining $R(\be)$ we obtain exactly the formulas in {\bf (i)}.
\end{rem}

\

\bigskip

It is well known \cite{mumford} that in the case of a complex torus $X =
V/\Lambda$ the map
\begin{equation} \label{eq:factors}
H^{1}(\Lambda, H^{0}(V,{\mathcal O}^{\times}_{V})) \to
H^{1}(X,{\mathcal O}^{\times})
\end{equation}
is an isomorphism and thus gives a group cohomology description of the
Picard group of $X$.  Our goal is to obtain an analogous description
for the non-commutative tori $\bbX_{\bpi}$. To set this up, note that
the relative spaces $\bbX_{\bpi} \to \bbD$ fit with the discussion of
modules in the previous paragraph. Indeed, by definition the
non-commutative torus $\bbX_{\bpi}$ is the ringed space $(X,\Ah{X})$,
which is the Moyal quantization of the Poisson torus $(X,\bpi)$. In
section \ref{sec:moyal} this ringed space was constructed as the
quotient of the Moyal ringed space $(V,\Ah{V})$ by the translation
action of the lattice $\Lambda \subset V$. In particular
$\Ah{V}^{\times}$-valued factors of automorphy for $\Lambda$ will
describe certain left
(or right) locally free rank one modules on $\bbX_{\bpi}$ and we will
have a map of cohomology sets:
\begin{equation} \label{eq:factorsXpi}
H^{1}(\Lambda,H^{0}(V,\Ah{V}^{\times})) \to H^{1}(X,\Ah{X}^{\times}).
\end{equation}
In fact, the map \eqref{eq:factorsXpi} is an isomorphism of pointed
sets and so every left (or right) locally free rank one module admits
a description via a factor of automorphy. The fact that
\eqref{eq:factorsXpi} is an isomorphism will follow by the standard
reasoning of \cite{mumford} from the fact that every $\Ah{V}$ locally
free left module of rank one is trivial. The latter statement can be
proven by an order by order analysis of the non-commutative cocycles
in $\check{Z}^{1}(V,\Ah{V}^{\times})$ or more generally in
$\check{Z}^{1}(V,\op{Aut}_{\Ah{V}\text{-mod}}(\Ah{V}^{\oplus n}))$. Since this
argument is somewhat technical we have included it in
Appendix~\ref{app:trivial}.

In contrast with the commutative case, in the non-commutative context,
the properties of being rank one and invertible no longer
coincide. Therefore it is important to differentiate between locally
free rank one left modules and invertible bimodules, both of which can
lay claim to be non-commutative line bundles.  An invertible bimodule
is often taken as the definition of a line bundle on a non-commutative
space. However, we have found that on our non-commutative tori
bimodules do not behave flexibly enough when one looks at
families. More precisely it turns out that for a non-degenerate
Poisson structure, non-constant holomorphic families of degree zero
line bundles on $X$ do not in general admit a consistent quantization:

\begin{prop} \label{prop:invertiblefamilies}  Let $X$ be a complex torus
  and suppose $\bpi \in H^{0}(X,\wedge^{2}T_{X})$ is a non-degenerate
  holomorphic Poisson structure. Let $S$ be a compact complex space and suppose
  $L \to S\times X$ is a holomorphic line bundle whose restriction to
  each slice $\{ s \} \times X$ has first Chern class zero. Suppose
  that we can find a holomorphic family $\mycal{L} \to S\times
  \bbX_{\bpi}$ of invertible bimodules with the property
  $\mycal{L}/\hbar \cong L$. Then the classifying map $\kappa_{L} : S
  \to X^{\vee}$ corresponding to $L$ is constant.
\end{prop}
{\bf Proof.} An invertible bimodule on $\bbX_{\bpi}$ is a sheaf of
$\Ah{X,\bpi}$ bimodules on $X$ which is locally free and of rank one
when considered both as a left and a right module.

Note that for any bimodule $\mycal{V}$ on $\bbX_{\bpi}$ the sheaf
$\mycal{V}/\hbar\mycal{V}$ is a Poisson module in the sense of
\cite[Appendix~A.5]{GinzburgKaledin}. Furthermore since $\bpi$ is
non-degenerate the category of Poisson modules on $(X,\bpi)$ is
equivalent to the category of $D$-modules on $X$
\cite[Appendix~A.6]{GinzburgKaledin}, \cite{Kaledin} and so we get a
well defined functor from the category of finitely generated bimodules
on the formal non-commutative space $\bbX_{\bpi}$ to the category of
finite rank complex local systems on $X$. Conversely, given a complex
local system $\boldsymbol{V}$ on $X$ we have an obvious
$\Ah{X,\bpi}$ bimodule structure on the sheaf
$\boldsymbol{V}\otimes_{\mathbb{C}} \Ah{X,\bpi}$. These functors are
easily seen to be inverse equivalences of each other. 

In particular given a holomorphic family $\mycal{L} \to S\times
\bbX_{\bpi}$ of invertible bimodules we get a holomorphic family
$\boldsymbol{L} \to S\times X$ of rank one complex local systems on $X$
which is parameterized by $S$. By assumption the family of holomorphic
line bundles underlying these local systems is precisely $L \to S\times
X$. In other words, the existence of $\mycal{L}$ implies that the
map $\kappa_{L} : S \to X^{\vee}$ lifts to a holomorphic
map $\kappa_{\boldsymbol{L}} : S \to \op{Loc}(X,1)$ from $S$ to the
moduli of rank one local systems on $X$. In other words we have a
commutative diagram of complex spaces
\[
\xymatrix{
& \op{Loc}(X,1) \ar[d]^-{p} \\
S \ar[ru]^-{\kappa_{\boldsymbol{L}}} \ar[r]_-{\kappa_{L}} & \, X^{\vee}
}
\]
with $p$ being the natural projection and $\kappa_{L}$,
$\kappa_{\boldsymbol{L}}$ being the classifying maps for $L$ and
$\boldsymbol{L}$. On the other hand the moduli space $\op{Loc}(X,1)$
is Stein (in fact isomorphic to $(\mathbb{C}^{\times})^{2g}$) and so
the map $\kappa_{\boldsymbol{L}}$ must be constant. Hence $\kappa_{L}$
is constant and the lemma is proven. \ \hfill $\Box$

\

\medskip

\noindent
Our aim is to deform the Picard variety $X^{\vee} = \op{Pic}^{0}(X)$
along with the non-commutative deformation $\bbX_{\bpi}$ of $X$ so that
the Fourier-Mukai transform deforms as well. A natural choice will be
to try and deform $X^{\vee}$ to the moduli of quantum line bundles on
$\bbX_{\bpi} \to \bbD$. If we attempt to do this with the
interpretation of a quantum line bundle as an invertible bimodule, then
we will end up with an obstructed moduli problem as explained in
Lemma~\ref{prop:invertiblefamilies}. As we will see below, this problem
does not occur if we work with rank one locally free left
$\Ah{X,\bpi}$-modules. This motivates the following:

\begin{defn} \label{defn:linebundle} A line bundle on a formal non-commutative
  space $(X,\Ah{X})$ is a sheaf \linebreak $\mycal{L} \to X$ of
  left-$\Ah{X}$-modules which is locally isomorphic as  a
  left module to $\Ah{X}$.
\end{defn}

\

\noindent
The bimodule properties of a line bundle are not
completely lost however.

\begin{prop}
Let $(X,\Ah{X})$ be a deformation quantization of
$(X,\mathcal{O})$ and let $L$ be a line bundle on $(X,\Ah{X})$.
Then $(X,{\cEnd}_{\Ah{X}} (L))$ is again a deformation quantization of
$(X,\mathcal{O})$ and thus $L$ is a left-$(X,\Ah{X})$ and a
right-$\cEnd_{\Ah{X}}(L)$ module. Furthermore these two actions
commute with each other and $L$ is a Morita equivalence bimodule.
\end{prop} 
{\bf Proof.}  The algebra ${\cEnd}_{\Ah{X}}(L)$ is naturally a
$\Ch$-algebra. Since flatness is a local condition and
${\cEnd}_{\Ah{}|_U}(L|_U)\cong \Ah{X}^{\op{op}}|_U$ we see that
${\cEnd}_{\Ah{}}{L}$ is a flat $\Ch$-module. Also
\[
{\cEnd}_{\Ah{X}}(L)/\hbar\cong {\cEnd}_{\Ah{X}/\hbar}(L/\hbar)\cong
\mathcal{O}
\]
and hence ${\cEnd}_{\Ah{X}}(L)$ is a deformation quantization 
of $X$. The Morita equivalence statement is straightforward.
\ \hfill $\Box$

\

\medskip

\begin{rem}
The previous proposition shows that with our definition, a
non-commutative line bundle $L$ implements an equivalence between the
category of $\Ah{X}$ modules and the category of
${\cEnd}_{\Ah{X}}(L)^{\op{op}}$ modules. In the commutative case this
reduces to the standard interpretation of a line bundle as an
autoequivalence of the category of sheaves. This idea of a line bundle
is very natural physically and was exploited before in
\cite{Wess}. Mathematically, it can be motivated by the natural
expectation that a non-commutative space should not just be taken to be
a ringed space $(X,\Ah{X})$ but should be a Morita equivalence class
of such spaces (perhaps only locally defined). There is other evidence
for this as well and we will hopefully pursue this in a future paper.
\end{rem}

\

\medskip

\noindent
For future reference we record some simple properties of the
deformation quantizations arising from non-commutative line bundles

\begin{prop} \label{prop:incenter}
Let $(X,\Ah{X})$ be a deformation quantization of
$(X,\mathcal{O})$. If $H^0(X,\mathcal{O})=\mathbb{C}$, then for
any line bundle $L$, the natural sheaf map $\zenter(\Ah{X})\to
{\cEnd}_{\Ah{X}}(L)$ 
given by the center of the algebra acting by left multiplication 
induces an isomorphism
\[
H^0(X;\zenter(\Ah{X}))\to H^0(X,{\cEnd}_{\Ah{X}}(L))
\]
\end{prop}
{\bf Proof.} We will prove that the composition 
\begin{equation} \label{eq:throughiso}
\Ch \to  H^0(X, \zenter(\Ah{X})) \hookrightarrow
H^0(X,{\cEnd}_{\Ah{X}}(L)) 
\end{equation}
 is an isomorphism. 

Suppose $(X,\mathcal{O})$ is a complex manifold with
$H^0(X;\mathcal{O})=\mathbb{C}$ and let $(X,\mycal{B})$ be a
deformation quantization of $(X,\mathcal{O})$. Then the map
$\Ch\to \mycal{B}$ induces an isomorphism $H^0(X;\Ch)\to
H^0(X,\mycal{B})$.

Indeed, both $\Ch$ and $H^0(X;\mycal{B})$ are complete filtered
algebras where the filtrations are given by
\[
H^0(X;\mycal{B})_k=\hbar^k\cdot H^0(X;\mycal{B})
\]
and similarly for $\Ch$. The map $\Ch\to H^0(X;\mycal{B})$ is a
filtered map and induces a map
\[
\op{gr}(f):\op{gr}\Ch\cong \mathbb{C}[\hbar]\to
\op{gr}H^0(X;\mycal{B})
\]
Now,
\begin{equation}
\begin{split}
\op{gr}_{k}H^0(X;\mycal{B}) & =
H^0(X,\mycal{B})_{k}/H^0(X,\mycal{B})_{k+1 }\\
 & =
H^0(X,\hbar^k\mycal{B})/H^0(X,\hbar^{k+1}\mycal{B})\\
 & \cong
H^0(X,\mycal{B}/\hbar\mycal{B}) \cong H^0(X,\mathcal{O}) \cong \mathbb{C}
\end{split}
\end{equation}
because $H^0(\hbar^k\mycal{B})$ surjects onto
$H^0(\mycal{B}/\hbar\mycal{B})$, and because
$\lambda\in\mathbb{C}$ is covered by $\hbar^k\lambda$. However if $f :
A\to B$ is a 
filtered map of complete filtered algebras and $\op{gr}(f)$ is
an isomorphism, then $f$ is also an isomorphism.  Indeed, recall that a
complete filtered algebra $A$ is a $\mathbb{C}$-algebra equipped 
with a decreasing filtration by $2-$sided ideals $A_{i}$, $A = A_{0}
\supseteq A_{1} \supseteq \cdots$ which satisfy $A_{i} A_{j} \subseteq
A_{i+j}$, and completeness: $A = \lim_{k\to
  \infty}\displaylimits A/A_{k}$ 

Now by completeness, it is enough to show that
the maps $f_{k} : A/A_{k} \to B/B_{k}$ induced by $f$ are isomorphisms
for all $k$.  Since by assumption
$\op{gr}_{k}(f) : \op{gr}_{k}(A) \to \op{gr}_{k}(B)$ is an
isomorphism for all $k$ 
and since 
\[ A/A_{1} = \op{gr}_{1}(A) \to \op{gr}_{1}(B) = B/B_{1}\] 
is an isomorphism, the claim follows by induction using
the commutative diagram of short exact sequences 
\[
\xymatrix{
0 \ar[r] &\ar[d] A/A_{k-1} \ar[r] \ar[d] & A/A_{k} \ar[r] \ar[d] &
\op{gr}_{k} A \ar[r] \ar[d] &  0 \\ 
0 \ar[r] & B/B_{k-1}\ar[r]  & B/B_{k} \ar[r]  & \op{gr}_{k} B \ar[r]  &  0 
}
\]
This shows that \eqref{eq:throughiso} is an isomorphism and completes
the proof of the lemma. \ \hfill $\Box$

\

\bigskip

\noindent
We now begin to develop the necessary tools to analyze families of
line bundles. We have the following statement, which parallels the
classical see-saw lemma \cite{mumford}: 

\begin{prop} \label{prop:see-saw}
Let $\bbX = (X,\Ah{X})$ be a deformation quantization of
$(X,\mathcal{O}_{X})$ and let \linebreak 
$\bbY = (Y,\mathcal{O}_{\bbY})$ be a complex
manifold over 
$\bbD$. Also, assume that $H^0(X,\mathcal{O})\cong\mathbb{C}$.
Consider the deformation quantization $\bbX\times_{\bbD} \bbY = 
(X\times Y,\Ah{X\times
Y})=(X\times Y,p_1^{-1}\Ah{X}\widehat{\otimes}_{{}_{\Ch}}
p_2^{-1}\mathcal{O}_{\bbY})$. 
Let $L_1$ and $L_2$ be two line bundles on $\bbX\times_{\bbD} \bbY$
such that for all $y\in Y$, there is an isomorphism 
\[
\phi_y\in \mathcal{I}som(L_1|_{X\times \{y\}}, L_2|_{X\times \{y\}}).
\]
Then:
\begin{itemize}
\item[{\bf (a)}] There is a global isomorphism of sheaves of algebras
\[
{\cEnd}_{\Ah{X\times Y}}(L_1)\cong {\cEnd}_{\Ah{X\times Y}}(L_2)
\]
which on each fiber $X \times \{y\}$ restricts to the isomorphism
$\phi_{y}\circ (-)\circ\phi_{y}^{-1}$.
\item[{\bf (b)}] There exists a line bundle $M$ on $\bbY$ and an
  isomorphism  $\bbp_{\bbY}^{*} M\otimes_{\Ah{X\times Y}} L_{2}
  \widetilde{\to} 
  L_{1}$. 
\end{itemize}
\end{prop}
{\bf Proof.} Part {\bf (a)} of this proposition is a new element of
the see-saw principle, specific to the deformation quantization
situation.  To prove {\bf (a)}, chose a cover $\mathfrak{U}$
of $Y$ and elements $\phi_{U} \in \mathcal{I}som_{\Ah{X \times Y}}
(L_{1}, L_{2}) (X \times U)$ for all $U \in \mathfrak{U}$.  Denote by
$\psi_{U}: {\cEnd}_{\Ah{X \times Y}}(L_{1})|_{X \times U} \to
{\cEnd}_{\Ah{X \times Y}}(L_{2})|_{X \times U}$ the induced local
isomorphisms of algebras $\psi_{U}(f) = \phi_{U} \circ f \circ
\phi_{U}^{-1}$.  Using Proposition \ref{prop:incenter} we see that the
elements
\[
\phi_{V}^{-1} \circ \phi_{U} \in 
H^{0}(X \times (U \cap V), \mathcal{I}som_{\Ah{X \times Y}}(L_{1}))
\] 
are in fact in 
$H^{0}(X \times (U \cap V), p_{1}^{-1} 
\zenter(\Ah{X}) \widehat{\otimes}_{{}_{\Ch}} p_{2}^{-1} 
\mathcal{O}_{U \cap V})$.  Therefore 
the $\psi_{U}$ patch to a global isomorphism 
${\cEnd}_{\Ah{X\times Y}}(L_1)\cong {\cEnd}_{\Ah{X\times Y}}(L_2)$ of 
sheaves of algebras.

\

\smallskip

\noindent
The proof of  {\bf (b)} is essentially the same as the see-saw proof
    found in \cite[Section~II.5]{mumford}. We will only discuss the
    modifications of the argument needed in the non-commutative
    setting.  Set $\mycal{E} :=
    {\cEnd}_{\Ah{X\times 
    Y}}(L_1)\cong {\cEnd}_{\Ah{X\times Y}}(L_2)$. Then $L_{1}$ and
    $L_{2}$ are sheaves on $X\times Y$ which are locally-free rank one
    right modules over the sheaf of algebras $\mycal{E}$. Hence the
    sheaf $L_{2}^{\vee} := \mycal{H}om_{\Ah{X\times
    Y}}(L_{2},\Ah{X\times Y})$ has a natural structure of a left
    $\mycal{E}$-module and a right $\Ah{X\times Y}$-module. Consider
    now the tensor product $L_{1}\otimes_{\mycal{E}}
    L_{2}^{\vee}$. This is a sheaf on $X\times Y$ which is a
    $\Ah{X\times  Y}$-bimodule and for every $y \in Y$ satisfies 
\[
(L_{1}\otimes_{\mycal{E}}
    L_{2}^{\vee})_{|X\times \{y \}} \cong \Ah{X}
\]
as an $\Ah{X}$-bimodule. By Proposition \ref{prop:incenter} applied to
    the trivial line bundle on $X$ we have $H^{0}(X,\Ah{X}) = \Ch$ and
    therefore $M :=
    p_{Y*}(L_{1}\otimes_{\mycal{E}} L_{2}^{\vee})$ is a rank one
    locally free
    left module over $\mathcal{O}_{\bbY}$. By adjunction the
    identity endomorphism 
\[
\op{id}_{M} \in \op{Hom}_{\bbY}(M,M) =
    \op{Hom}_{\bbY}(p_{Y*}(L_{1}\otimes_{\mycal{E}} L_{2}^{\vee}),
    p_{Y*}(L_{1}\otimes_{\mycal{E}} L_{2}^{\vee}))
\] 
corresponds to a map 
\[
\bbp_{\bbY}^{*}M\otimes_{\Ah{X\times Y}} L_{2} \to L_{1}.
\]
The check that this map is an isomorphism is exactly the same as in
the commutative situation.
\ \hfill $\Box$

\subsection{The Poincar\'{e} sheaf} \label{ss:poincare}

We now describe the factor of automorphy that defines the Poincar\'{e}
sheaf in our context. Let $\bpi$ be a holomorphic Poisson structure on
$X$ and let $\bB$ denote the corresponding $B$-field on $X^\vee$. The
Poisson structure $\bpi$ lifts to a Poisson structure on $V$ which
will be denoted by the same letter. Consider the Poisson structure on
$V\times \overline{V}^{\vee}$ which is $\bpi$ on $V$ and $0$ on
$\overline{V}^{\vee}$.  We can then form the corresponding Moyal
quantizations $(V\times \overline{V}^{\vee},
\Ah{V\times\overline{V}^{\vee},(\bpi,0)})$ and $(X\times
X^{\vee},\Ah{X\times X^{\vee},(\bpi,0)})$ (see
Section~\ref{sssec-moyal}). In what follows we will drop the labels
$(\bpi,0)$ from our algebras since they will be clear from the
context.  Note that $(X\times X^{\vee},\Ah{X\times
X^{\vee},(\bpi,0)})$ is just the non-commutative space
$\bbX_{\bpi}\times_{\bbD} \bbX^{\vee}$ which is the moduli
space of the stack $\bbX_{\bpi}\times_{\bbD}
\left({}_{\bB}\bbX^{\vee}\right)$.

Our goal is to construct a deformation of the Poincar\'{e} line bundle
$\mathcal{P} \to X\times X^{\vee}$ to a line bundle on the
non-commutative stack $\bbX_{\bpi} \times_{\bbD}
\left({}_{\bB}\bbX^{\vee}\right)$. 
Why a stack? Even classically the moduli problem of topologically
trivial line bundles on a space $Z$ leads most naturally to an
analytic stack $\mycal{P}ic^{0}(Z)$ which is an
$\mathcal{O}^{\times}$-gerbe over the usual Picard variety
$\op{Pic}^{0}(Z)$. This is rarely discussed since the gerbe
$\mycal{P}ic^{0}(Z) \to \op{Pic}^{0}(Z)$ is trivial. Indeed, we can
construct a trivialization of this gerbe by looking at the moduli
problem of framed line bundles on $Z$ with the framing specified at a
fixed point $z \in Z$. This feature of the moduli of line bundles 
does not persist in families. If we look at the relative Picard
problem for a smooth family $Z \to B$, then the moduli stack is not
necessarily a trivial gerbe since now a trivialization depends on a
framing along a section $B \to Z$, which may not exist.
In our case, $\bbX^{\vee} \to \bbD$ should be thought of as the
relative Picard variety $\op{Pic}^{0}(\bbX_{\bpi}/\bbD)$ and
${}_{\bB}\bbX^{\vee} \to \bbD$ should be thought of as the stack of
relative line bundles $\mycal{P}ic^{0}(\bbX_{\bpi}/\bbD)$ (this will
be justified in section~\ref{s:Pic}). By
construction the gerbe ${}_{\bB}\bbX^{\vee} \to \bbX^{\vee}$ is
non-trivial, and indeed we don't expect
$\mycal{P}ic^{0}(\bbX_{\bpi}/\bbD) \to \op{Pic}^{0}(\bbX_{\bpi}/\bbD)$
to be trivial since $\bbX_{\bpi} \to \bbD$ has no sections.

To define the deformation $\mycal{P} \to \bbX_{\bpi}\times_{\bbD}
\left({}_{\bB}\bbX^{\vee}\right)$ of $\mathcal{P} \to X\times
X^{\vee}$ we will use 
the description of line bundles via factors of automorphy. Since the
$\bbX_{\bpi}$ is the $\Lambda$-quotient of the non-commutative space
$(V,\Ah{V,\bpi})$ and ${}_{\bB}\bbX^{\vee}$ is the $\Gamma$ quotient
of $(\overline{V}^{\vee},\mathcal{O}[[\hbar]])$ we will need a
$\Lambda\times \Gamma$-factor of automorphy:
\[
\phi:\Lambda\times \Gamma\to
H^{0}(V\times V^{\vee}, \cAut(\Ah{V\times
  \overline{V}^{\vee}})) = 
\Ah{V\times \overline{V}^{\vee}}^\times(V\times 
\overline{V}^{\vee}).
\]
Here $ \op{Aut}(\Ah{V\times
  \overline{V}^{\vee}})$ denotes the automorphisms of $\Ah{V\times
  \overline{V}^{\vee}}$ considered as a left module over itself.
We define $\phi$
by the formula
\begin{equation}\label{PSFA}
\begin{split}
\phi(\lambda,(\xi, & z))(v,l) =
\\
& =
z\exp\left({\pi\sqrt{-1}\op{Im}(\langle\xi,\lambda\rangle)}\right)
\exp\left({\pi\left(\langle  
   l,\lambda\rangle 
  +\overline{\langle\xi,v\rangle}\right)}\right)
\exp\left({\frac{\pi}{2}\left(\langle\xi,  
  \lambda\rangle
  +\overline{\langle\xi,\lambda\rangle}\right)}\right)
\\
& = z \exp\left(\pi \left(\langle l + \xi, \lambda \rangle +\overline{
  \langle \xi,v \rangle}\right)\right) 
\end{split}
\end{equation}
where $\lambda\in \Lambda$, $(\xi,z)\in\Gamma$ and $(v,l)\in
V\times \overline{V}^{\vee}$, and $\langle l, v \rangle = l(v)$. 
In particular $\langle l, v \rangle$ is complex anti-linear as a function of
$v$ and complex linear as a function of $l$. We note that the only term in this
product that involves $\hbar$ is $z$. The exponentials and the
products of terms in this formula can be viewed either as
$\star$-exponentials and products, or as ordinary exponentials and
products of functions on $V\times \overline{V}^{\vee}$. This is
unambiguous as we will see in the proof of the next proposition.

\begin{prop}
$\phi$ is a factor of automorphy.
\end{prop}
{\bf Proof.} 
To see that $\phi$ is a factor of automorphy, we have to check the
condition \eqref{cocyclecondition}, which says
\[
\phi(\lambda_{1} + \lambda_{2},\xi_{1}
+ \xi_{2} ,z_{1} z_{2} c(\xi_{1}, \xi_{2}) )
= \phi(\lambda_{2}, \xi_{2}, z_{2})
\star ([\phi(\lambda_{1}, \xi_{1}, z_{1})] \cdot (\lambda_{2}, \xi_{2}, z_{2}))
\]
The computation which establishes the factor of automorphy condition for the 
classical Poincar\'e cochain 
\cite{polishchuk} immediately reduces us 
in this case to showing the equality 

\begin{equation}\label{needtoshow} 
c(\xi_{1}, \xi_{2})
\exp\left({\pi\overline{\langle\xi_1 +
\xi_2,v\rangle}}\right) =
\exp\left({\pi\overline{\langle\xi_2,v\rangle}}\right)\star 
\exp\left({\pi\overline{\langle\xi_1,v\rangle}}\right).
\end{equation}
To that end consider two points $f$ and $g$ in $V^{\vee}$ viewed
as linear functions on $V$. If $P$ is the bidifferential operator
associated with $\bpi$, then using the formula
\eqref{eq:bidifferential} we compute (directly, or by the
Campbell-Baker-Hausdorff formula) 
\[
(\exp{f})\cdot P \cdot (\exp{g}) = \{ f,g \} \exp{(f+g)}.
\]
Therefore
\begin{equation} \label{eq-star-linear}
(\exp{f}) \star (\exp{g}) = \exp{\left(\hbar  \{
    f,g \}\right)}\exp{(f+g)}. 
\end{equation}
Given $\xi \in \overline{V}^{\vee}$ let  $f_\xi  :=
\pi\overline{\langle\xi, \bullet\rangle}$ denote the corresponding
linear function on $V$. Then \eqref{eq-star-linear} implies that
the $\star$-inverse of $\exp(f_\xi)$ is $\exp(f_{-\xi})$.

\noindent
Using  \eqref{eq-star-linear} in the right hand side of equation
\eqref{needtoshow}  and canceling $\exp(f_{\xi_{1}} + f_{\xi_{2}})$
from both sides, it remains only to show that
\begin{equation}\label{showme}
c(\xi_{1}, \xi_{2})
= \exp \left({\hbar \{
  f_{\xi_2},f_{\xi_1} \}}\right).
\end{equation}
\noindent
However, 

\begin{equation}\label{blah}
\{ f_{\xi_2},f_{\xi_1} \} = \bpi \cntrct ({df_{\xi_{2}}}\wedge
   {df_{\xi_{1}}}) =  
\pi^{2} \bpi \cntrct (\bar{\xi}_{2}\wedge \bar{\xi}_{1}) = \pi^{2}
   \bB(\xi_{2}, \xi_{1}).
\end{equation}

Now we simply recall that we defined $c$ in equation \eqref{eq:defofc} by 
$c(\xi, \xi') = \exp\left({\hbar \pi^{2}
\bB(\xi', \xi)}\right)$ to see that equation \eqref{showme} is satisfied so
we are done.

\ \hfill $\Box$

\section{The equivalence of categories}

We are now ready to compare the sheaf theories of the non-commutative
torus $\bbX_{\bpi}$ and the gerbe $_{\bB}\bbX^{\vee}$. Just as in the
commutative case the equivalence question can be posed for different
derived categories of sheaves.

\subsection{Categories of sheaves} \label{ss:categories}  The basic
categories we need for 
$\bbX_{\bpi}$ will be:
\begin{description}
\item[$\sff{Mod}(\bbX_{\bpi})$] = the category of all sheaves of left
$\Ah{\bpi}$-modules on $X$.
\item[$\sff{Coh}(\bbX_{\bpi})$] = the category of all sheaves of
  coherent left $\Ah{\bpi}$-modules on $X$. By definition these are
  sheaves in $\sff{Mod}(\bbX_{\bpi})$ which locally in $X$ admit a
  finite presentation by free modules. 
\end{description}
We will write $\sff{D}^{*}(\bbX_{\bpi})$ for the derived
category of $\sff{Mod}(\bbX_{\bpi})$ and
$\sff{D}_{c}^{*}(\bbX_{\bpi})$ for the derived version of
$\sff{Coh}(\bbX_{\bpi})$. Here the decoration $*$ can be
anything in the set $\{ \varnothing, -, b\}$.

It is important to note that $\Ah{X,\bpi}$ is a coherent and
Noetherian sheaf of rings, see \cite[AII.6.27, AIII.3.24]{Bjork} for
the definitions. This implies that $\sff{Coh}(\bbX_{\bpi})$ is a full
abelian subcategory in $\sff{Mod}(\bbX_{\bpi})$ and also that
$\sff{D}_{c}^{*}(\bbX_{\bpi})$ is equivalent to the full subcategory
of $\sff{D}^{*}(\bbX_{\bpi})$ consisting of complexes with coherent
cohomology.  To check that $\Ah{X,\bpi}$ is coherent and Noetherian we
use a very general criterion due to Bj\"{ork}
\cite[Lemma~8.2,Theorem~9.6]{Bjork-rings}. Even though the
$\hbar$-filtration on $\Ah{X,\bpi}$ is decreasing, we can still apply
the technology of \cite{Bjork-rings} since his setup allows for
filtrations infinite in both directions and we can simply relabel the
filtration to make it increasing. To check the hypotheses of
\cite[Lemma~8.2]{Bjork-rings} we need to make sure that the stalks of
$\op{gr} \Ah{X,\bpi}$ are left and right Noetherian. This is clear
since (stalks of $\op{gr} \Ah{X,\bpi}$) = (stalks of
$\mathcal{O}_{X}[[\hbar]]$). To check the hypotheses of
\cite[Theorem~9.6]{Bjork-rings} we need to show that given a left
ideal $L$ in a stalk of $\Ah{X,\bpi}$ and elements $a_{1}, \ldots,
a_{s} \in L$, such that $\sigma(a_{1}), \ldots, \sigma(a_{s})$
generate $\sigma(L)$, then
\begin{equation} \label{eq:bjorkcondition}
\Sigma_{\nu} \cap L = \Sigma_{\nu-\nu_{1}}a_{1} + \cdots +
\Sigma_{\nu-\nu_{s}}a_{s} 
\end{equation}
holds for all $\nu$. Here $\Sigma_{\nu}$ is the $\nu$-th step of the
increasing filtration on the stalk of $\Ah{X,\bpi}$, i.e.
\[
\Sigma_{\nu} = \begin{cases} \hbar^{-\nu}\mycal{A} & \nu \leq 0 \\
  \mycal{A} & \nu > 0, \end{cases}
\]
$\sigma : \Ah{X,\bpi} \to \op{gr}\Ah{X,\bpi}$ is the symbol map, and
$\nu_{i}$ is the order of $a_{i}$, i.e. $a_{i} \in
\Sigma_{\nu_{i}}\setminus \Sigma_{\nu_{i-1}}$.

In fact the condition \eqref{eq:bjorkcondition} holds for any
deformation quantization. Indeed given $a_{1}, \ldots, a_{s}$ as above
and $\ell \in \Sigma_{\nu}\cap L$, 
we can choose $\alpha_{1}, \ldots, \alpha_{s}$ in the stalk of
$\Ah{X,\bpi}$ so that $\nu(\alpha_{i}) = \nu - \nu_{i}$ satisfying
\[
\ell - (\alpha_{1}\star a_{1} + \cdots \alpha_{s}\star a_{s}) \in
\Sigma_{\nu-1}\cap L.
\]
Now induction and the $\hbar$-adic completeness of
$\Ah{X,\bpi}$ finish the verification of \eqref{eq:bjorkcondition}.

\

\noindent
For the gerbe $_{\bB}\bbX^{\vee}$ the relevant categories are
\begin{description}
\item[$\sff{Mod}( _{\bB}\bbX^{\vee})$] = the category of sheaves of ${\mathcal
  O}_{ _{\bB}\bbX^{\vee}}$-modules.
\item[$\sff{Coh}( _{\bB}\bbX^{\vee})$] = the category of coherent
  ${\mathcal
  O}_{ _{\bB}\bbX^{\vee}}$-modules.
\end{description}

\

\noindent
Since the gerbe $_{\bB}\bbX^{\vee}$ admits a presentation
$_{\bB}\bbX^{\vee} = [\overline{\bbV}^{\vee}/\bGamma]$ as a global
quotient we can explicitly describe $\sff{Mod}( _{\bB}\bbX^{\vee})$
and $\sff{Coh}( _{\bB}\bbX^{\vee})$ the categories of
$\bGamma$-equivariant ${\mathcal O}_{\overline{\bbV}^{\vee}}$-modules
and coherent $\bGamma$-equivariant ${\mathcal
O}_{\overline{\bbV}^{\vee}}$-modules respectively. Since
$_{\bB}\bbX^{\vee}$ is an ${\mathcal O}^{\times}$-gerbe, these
categories decompose into orthogonal direct sums $\sff{Mod}(
_{\bB}\bbX^{\vee}) = \coprod_{k \in {\mathbb Z}} \sff{Mod}(
_{\bB}\bbX^{\vee},k)$ and $\sff{Coh}( _{\bB}\bbX^{\vee}) = \coprod_{k
\in {\mathbb Z}} \sff{Coh}( _{\bB}\bbX^{\vee},k)$, labeled by the
character $k$ of the stabilizer ${\mathbb C}^{\times}$. Explicitly
$\sff{Mod}_{k}( _{\bB}\bbX^{\vee})$ and $\sff{Coh}_{k}(
_{\bB}\bbX^{\vee})$ are respectively the categories of
$\bGamma$-equivariant sheaves and coherent, $\bGamma$-equivariant
sheaves for which the action of the center ${\mathcal
O}^{\times}_{\bbD}$ is the $k$-th power of the tautological
action. These also admit an alternative interpretation as categories
of $k\bB$-twisted sheaves and coherent, $k\bB$-twisted sheaves on
$\bbX^{\vee}$ \cite{giraud,caldararu}. Finally, we will write
$\sff{D}^{*}( _{\bB}\bbX^{\vee}, k)$ for the corresponding derived
categories.

\

\noindent
{\bf Digression on quasi-coherent sheaves} \ An unpleasant phenomenon
of the sheaf theory in analytic geometry is the fact that for generic
analytic spaces the categories of analytic coherent sheaves tend to be
fairly small and boring \cite{Verbitsky}. In particular the category
of analytic coherent sheaves is not a good invariant in the analytic
world. This contrasts with the algebraic category where a Noetherian
scheme can be reconstructed from its category of coherent sheaves
\cite{Gabriel}. Therefore the common wisdom is that for an analytic
space $X$ the geometry is captured better by more general categories
of $\mathcal{O}_{X}$-modules, e.g. quasi-coherent ones or all
$\mathcal{O}_{X}$-modules.

Unfortunately the notion of a quasi-coherent analytic sheaf is a bit
murky and there are several competing definitions \cite{ramis-ruget},
\cite{putinar},\cite[Chapter~11.9]{taylor},\cite{orlov-qc} in the
literature. The relationships among these definitions are not clear in
general. We will adopt the definition of \cite[Chapter~11.9]{taylor}
which is best suited to our setup and we will comment how this
definition compares to the others in our specific case.

\begin{defn} \label{defn:taylor-qc} Suppose $X$ is an analytic space
  and let $F$ be an analytic sheaf of $\mathcal{O}_{X}$-modules. The
  sheaf $F$ is called {\em quasi-coherent} if for every point $x \in
  X$ we can find a compact Stein neighborhood $x \in K \subset X$
  having the Noether property, and a module $M$ over
  $\Gamma(K,\mathcal{O}_{K})$, so that  
\[
F_{|K} \cong \widetilde{M} :=
\mathcal{O}_{K}\otimes_{\Gamma(K,\mathcal{O}_{K})} M.
\]
As usual, in this formula, $\Gamma(K,\mathcal{O}_{K})$ and $M$ are
viewed as constant sheaves on $K$. 

Similarly, if $\bbX_{\bpi} = (X,\Ah{X,\bpi})$ is a holomorphic
deformation quantization of $X$, we will say that a sheaf $\mycal{F}
\in \sff{Mod}(\bbX_{\bpi})$ of left $\Ah{X,\bpi}$-modules is {\em
quasi-coherent} if for every $x \in X$ we can find a compact Stein
neighborhood $x \in K \subset X$ having the Noether property, and a
left module $\mycal{M}$ over $\Gamma(K,\Ah{X,\bpi})$, so that
\[
\mycal{F}_{|K} \cong \widetilde{\mathcal{M}} :=
\left({\Ah{X,\bpi}}_{|K} \right)
\otimes_{\Gamma(K,\Ah{X,\bpi})} \mathcal{M}.
\]
\end{defn}
Recall that a compact Stein set is a compact analytic subspace in some
$\mathbb{C}^{n}$ which can be realized as the intersection of a nested
sequence of Stein neighborhoods. The Noether property of a compact
Stein space $K$ is that $\Gamma(K,\mathcal{O}_{K})$ a Noetherian Stein
algebra. It is known \cite{frisch,langmann} that compact analytic
polyhedra (i.e.  subsets of a Stein space defined by finitely many
inequalities of the form $|f|\leq 1$, for $f$ holomorphic) are compact
Stein sets that have the Noether property. In particular polydiscs are
compact Stein and Noether. More generally every point in an analytic
space has a basis of compact Stein Noether neighborhoods (see also
\cite{gr-stein}).

The above notion of quasi-coherence is somewhat unconventional from
the point of view of Grothendieck sheaf theory on complex spaces. We
chose to work with it, since it is compatible with the more standard
points of view on quasi-coherence and in addition turns out to have a
very good behavior with respect to pullbacks and pushforwards. 

It is instructive to compare the quasi-coherence in the sense of
Definition~\ref{defn:taylor-qc} to the quasi-coherence of
\cite{ramis-ruget}, \cite{putinar}, and \cite{orlov-qc}. 

The most
general notion of quasi-coherence is the one given by Orlov in
\cite[Definition~2.6]{orlov-qc}. Orlov's definition works on an
arbitrary ringed site and is conceptually the closest to the usual
notion of quasi-coherence in algebraic geometry. In fact
Definition~\ref{defn:taylor-qc} is a special case of
\cite[Definition~2.6]{orlov-qc}. For any complex space $X$ we can
consider the site $\sff{cSt}/X$ of compact Stein spaces, taken
with the analytic topology and ringed by the sheaf of analytic
functions. Any analytic sheaf $F$ of $\mathcal{O}_{X}$-modules is 
a sheaf on $\sff{St}/X$ gives rise to a sheaf $\sff{c}F$ on
$\sff{cSt}/X$ and the quasi-coherence of $F$ in the sense of
Definition~\ref{defn:taylor-qc} is simply the quasi-coherence of
$\sff{c}F \to  \sff{cSt}/X$ in the sense of
\cite[Definition~2.6]{orlov-qc}. 

The other definition is the Ramis-Ruget definition of quasi-coherence
\cite{ramis-ruget}, \cite{putinar} which is historically the first
one. In their definition, a sheaf $F$ of $\mathcal{O}_{X}$ is called
quasi-coherent if locally on $X$ we can can realize $F$ as a
transversal localization of a module. More precisely, $F$ is
{\em algebraically Ramis-Ruget quasi-coherent} if for every two Stein opens
$V \subset U$ we have:
\begin{itemize}
\item $\Gamma(U,F)$ and $\Gamma(V,\mathcal{O})$ are Verdier
  transversal over $\Gamma(U,\mathcal{O})$, that is
\[
\Gamma(V,\mathcal{O})\stackrel{L}{\otimes}_{\Gamma(U,\mathcal{O})}
 \Gamma(U,F) =
\Gamma(V,\mathcal{O})\otimes_{\Gamma(U,\mathcal{O})} \Gamma(U,F). 
\]
\item The natural map
  $\Gamma(V,\mathcal{O})\stackrel{L}{\otimes}_{\Gamma(U,\mathcal{O})}
  \Gamma(U,F) \to \Gamma(V,F)$ is an isomorphism.
\end{itemize}
In fact the Ramis-Ruget notion of quasi-coherence requires that the
sheaf $F$ is an analytic sheaf of nuclear Frechet
$\mathcal{O}$-modules, and all the tensor products are the completed
ones. This is necessary in their setup since they are concerned with
the Grothendieck duality theory and in particular need the double dual
of the space of sections of a sheaf to be isomorphic to itself. For
our purposes, the nuclear Frechet condition is irrelevant, and so we
talk only about the algebraic Ramis-Ruget quasi-coherence.

Note that the transversality condition in the Ramis-Ruget definition
implies that any sheaf $F \in \sff{Mod}(X)$ which is algebraically
Ramis-Ruget quasi-coherent is also quasi-coherent in the sense of 
Definition~\ref{defn:taylor-qc}. So our notion of quasi-coherence is
sandwiched between the Ramis-Ruget's function theoretic notion and
Orlov's general categorical notion.

The main advantage of Definition~\ref{defn:taylor-qc} is that the
localization functor on compact Stein spaces having the Noether
property is an exact functor \cite[Chapter~11.9]{taylor}. In
particular, given an analytic morphism $f : X \to Y$ the pushforward
$f_{*} : \sff{Mod}(X) \to \sff{Mod}(Y)$ preserves quasi-coherence. In
the terminology of \cite{orlov-qc} this means that the natural map
from $F$ to its coherator is an isomorphism. We write $\sff{Qcoh}(X)$
and $\sff{Qcoh}(Y)$ for the categories of quasi-coherent sheaves and
$f_{*}$ and $f^{*}$ for the corresponding pullback and pushforward
functors.

Finally, observe that the exactness of compact Stein localization of
\cite[Chapter~11.9]{taylor} also holds for the sheaves in the category
$\sff{Qcoh}(\bbX_{\bpi}))$. This is an immediate consequence of the
the exactness in \cite[Chapter~11.9]{taylor} and the fact that
$\Ah{X,\bpi}$ is a coherent and Noetherian sheaf of complete
$\Ch$-algebras.

For future reference we will write
$\sff{D}^{*}_{\op{qc}}(\bbX_{\bpi})$ and
$\sff{D}^{*}_{\op{qc}}({}_{\bB}\bbX^{\vee},k)$ for the derived
categories of analytic sheaves with quasi-coherent cohomologies.

\subsection{The main theorem} \label{ss:maintheorem}
In section \ref{ss:poincare} we defined a Poincar\'{e} sheaf
$\mycal{P}$ on $\bbX_{\bpi}\times_{\bbD} {_{\bB}\bbX^{\vee}}$
deforming the classical Poincar\'{e} line bundle. By definition
$\mycal{P}$ is a sheaf on $X\times ( _{\bB}\bbX^{\vee})$ which is a
left $p_{1}^{-1}\Ah{\bpi}$-module and a right $p_{2}^{-1}{\mathcal
O}_{ _{\bB}\bbX^{\vee}}$-module, i.e. $\mycal{P} \in
\lrmod{p_{1}^{-1}\Ah{\bpi}}{p_{2}^{-1}{\mathcal O}_{
{}_{\bB}\bbX^{\vee}}}$. We will also need the (algebraic) dual sheaf
$\mycal{Q}$ of $\mycal{P}$. For our purposes it  will be convenient to
define $\mycal{Q}$
on the product of $X$ and ${}_{\bB}\bbX^{\vee}$ in which the order of
the factors is transposed. In other words we define $\mycal{Q}$ as the
sheaf on $({}_{\bB}\bbX^{\vee})\times X$ given by
\[
\mycal{Q} = \mycal{P}^{\vee} :=
  \op{Hom}_{{}_{p_{2}^{-1}\Ah{\bpi}}{\scriptstyle
  \sff{Mod}}_{p_{1}^{-1}{\mathcal O}_{  
  {}_{\bB}\bbX^{\vee}}}}(\bt^{*}\mycal{P},p_{2}^{-1}\Ah{\bpi}\otimes_{\Ch}
p_{1}^{-1}{\mathcal O}_{ _{\bB}\bbX^{\vee}}).
\]
Here $\bt : ({}_{\bB}\bbX^{\vee})\times X \to X\times
({}_{\bB}\bbX^{\vee})$ is the isomorphism transposing the factors. The
sheaf $\mycal{Q}$ is in $\lrmod{p_{1}^{-1}{\mathcal O}_{
  _{\bB}\bbX^{\vee}}}{p_{2}^{-1}\Ah{\bpi}}$ by definition.  The
left-right modules 
  $\mycal{P}$ and $\mycal{Q}$ can be used as kernels of integral
  transforms between derived categories. More precisely we have
  functors
\begin{equation} \label{eq:fmgerbetonc}
\xymatrix@R-2pc{
\phi_{\mycal{P}}^{[ _{\bB}\bbX^{\vee} \to \bbX_{\bpi}]} :
\hspace{-2pc} &
\sff{D}^{*}( _{\bB}\bbX^{\vee}, -1) \ar[r] &
\sff{D}^{*}(\bbX_{\bpi}) \\
 & \mycal{F} \ar[r] & Rp_{1*}( \mycal{P}
\otimes^{\mathbb L}_{p_{2}^{-1}
  {\mathcal O}_{ _{\bB}\bbX^{\vee}}} p_{2}^{-1} \mycal{F})
}
\end{equation}
and
\begin{equation} \label{eq:fmnctogerbe}
\xymatrix@R-2pc{
\phi_{\mycal{Q}}^{[\bbX_{\bpi}\to _{\bB}\bbX^{\vee}]} : \hspace{-2pc} &
 \sff{D}^{*}(\bbX_{\bpi})\ar[r] &
 \sff{D}^{*}( _{\bB}\bbX^{\vee}, -1)\\
 & \mycal{F} \ar[r] & Rp_{1*}(\mycal{Q} \otimes^{\mathbb
   L}_{p_{2}^{-1} \Ah{\bpi}} p_{2}^{-1} \mycal{F} )
}
\end{equation}
Here $p_{1}$ and $p_{2}$ denote the projections on the first and
second factors of the product $X \times (
_{\bB}\bbX^{\vee})$. Alternatively they can be thought of as the
projections onto the two factors in the fiber product $\bbX_{\bpi}
\times_{\bbD} { _{\bB}\bbX^{\vee}}$, but from that point of view we have
to use the fact that ${\mathcal O}_{ _{\bB}\bbX^{\vee}}$ is
commutative and regard $\mycal{P}$ as a left ${\mathcal
O}_{\bbX_{\bpi} \times_{\bbD} { _{\bB}\bbX^{\vee}}} =
p_{1}^{-1}\Ah{\bpi}\otimes_{{\mathbb C}[[\hbar]]} p_{2}^{-1}
{\mathcal
O}_{ { _{\bB}\bbX^{\vee}}}$-module. The two points of view are equivalent
but the second one introduces a certain asymmetry in the treatment of
$\mycal{P}$ and $\mycal{Q}$ and so we chose to consistently work with
left-right modules rather than left modules over tensor product
algebras.

The integral transform functors \eqref{eq:fmgerbetonc} and
\eqref{eq:fmnctogerbe} are well defined functors between all flavors
of derived categories that we consider. Indeed these integral
transforms are compositions of sheaf-theoretic pullbacks, tensor
products over sheaves of rings and derived direct images. The
sheaf-theoretic pullbacks are always exact so they do not cause any
trouble in the definition. The tensoring with the Poincar\'{e} sheaves
$\mycal{P}$ (or with $\mycal{Q}$) is also exact since by definition
$\mycal{P}$ is a flat $p_{1}^{-1}\Ah{\bpi}\otimes_{{\mathbb
C}[[\hbar]]} p_{2}^{-1} {\mathcal O}_{ { _{\bB}\bbX^{\vee}}}$ module
and both $\Ah{\bpi}$ and ${\mathcal O}_{ { _{\bB}\bbX^{\vee}}}$ are
flat over $\Ch$. Finally since the derived pushforward is defined by
means of injective resolutions, it always makes sense as a functor
between bounded derived categories. In fact, the push-forward makes
sense as a functor between unbounded or bounded above derived
categories since our maps are proper maps of finite homological
dimension and thus satisfy the sufficient condition
\cite[Section~II.2]{hartshorne-RD}. Alternatively one can use the
technique of Spaltenstein which allows us to define all functors on
the unbounded derived category by means of $K$-injective and $K$-flat
resolutions \cite{spaltenstein}.

Note also that by the discussion in the quasi-coherent digression
above all these functors preserve quasi-coherence. To see that the
integral transforms restrict to functors between the corresponding
coherent categories $\sff{D}^{*}_{c}$ we have to argue that the
analogue of Grauert's direct image theorem holds in our case. By now
there is plenty of technology in the literature to handle this, see
for example \cite{schneiders}. The details are routine so to keep down
the size of the paper we do not
include them. The essential hypotheses to check in
order to apply \cite{schneiders}, is that $\Ah{\bpi}$ is a sheaf of nuclear
Frechet, multiplicatively convex algebras. We have noted all these
conditions above. 

Now we are ready to state our main equivalence result

\begin{thm} \label{thm:main} Suppose $X$ is a $g$-dimensional
complex  torus equipped with 
a holomorphic Poisson structure $\bpi$. Let $\bbX_{\bpi} \to \bbD$ be the
  corresponding Moyal quantization, and let $ _{\bB}\bbX^{\vee} \to
  \bbD$ be the dual ${\mathcal O}^{\times}$-gerbe on
  $\bbX^{\vee}$. Then  we have isomorphisms of functors
\[
\begin{split}
\phi_{\mycal{Q}}^{[\bbX_{\bpi}\to {_{\bB}\bbX^{\vee}}]}
\circ \phi_{\mycal{P}}^{[ { _{\bB}\bbX^{\vee}} \to
    \bbX_{\bpi}]} & \cong \op{id}_{\sff{D}^{*}( {
    _{\bB}\bbX^{\vee}},-1)}[-g] \\ \smallskip
\phi_{\mycal{P}}^{[ _{\bB}\bbX^{\vee} \to
    \bbX_{\bpi}]} \circ \phi_{\mycal{Q}}^{[\bbX_{\bpi}\to
    _{\bB}\bbX^{\vee}]} & \cong \op{id}_{\sff{D}^{*}( {
    \bbX_{\bpi}})}[-g]
\end{split}
\]
In particular $\sff{D}^{*}(\bbX_{\bpi})$ and
$\sff{D}^{*}( _{\bB}\bbX^{\vee})$ are triangulated
equivalent. Similarly $\sff{D}^{*}_{c}(\bbX_{\bpi})$ (resp. 
$\sff{D}^{*}_{qc}(\bbX_{\bpi})$) and
$\sff{D}^{*}_{c}( _{\bB}\bbX^{\vee},-1)$ 
(resp. $\sff{D}^{*}_{qc}( _{\bB}\bbX^{\vee},-1)$) are triangulated
equivalent, and so $\bbX_{\bpi}$ and $ _{\bB}\bbX^{\vee}$ are
Fourier-Mukai partners.
\end{thm}
\noindent
{\bf Proof.}  Similarly to the classical case discussed in section
\ref{ssec:classicalfm}, the theorem will follow from the convolution
property of the integral transform functors. Specifically, we have
natural isomorphisms of functors:
\begin{equation} \label{eq:convolutions}
\begin{split}
\phi_{\mycal{Q}}^{[\bbX_{\bpi}\to {_{\bB}\bbX^{\vee}}]}
\circ \phi_{\mycal{P}}^{[ { _{\bB}\bbX^{\vee}} \to
    \bbX_{\bpi}]} & \cong \phi_{\mycal{Q}*\mycal{P}}^{[ {
      _{\bB}\bbX^{\vee}} \to { _{\bB}\bbX^{\vee}}]} \\ \smallskip
\phi_{\mycal{P}}^{[ { _{\bB}\bbX^{\vee}} \to
    \bbX_{\bpi}]} \circ \phi_{\mycal{Q}}^{[\bbX_{\bpi}\to
    {_{\bB}\bbX^{\vee}}]} & \cong
\phi_{\mycal{P}*\mycal{Q}}^{[\bbX_{\bpi} \to \bbX_{\bpi}]}.
\end{split}
\end{equation}
Our strategy will be to first compute the convolution
$\mycal{P}*\mycal{Q}$ of the kernel objects
$\mycal{P}$ and $\mycal{Q}$ and use the result to show that 
\[
\phi_{\mycal{P}}^{[ {}_{\bB}\bbX^{\vee} \to
    \bbX_{\bpi}]} \circ \phi_{\mycal{Q}}^{[\bbX_{\bpi}\to
    _{\bB}\bbX^{\vee}]}  \cong \op{id}_{\sff{D}^{*}( {
    \bbX_{\bpi}})}[-g].
\]
After that we finish the proof of the theorem by using the
Bondal-Orlov criterion to check that
$\phi_{\mycal{P}}^{[ {}_{\bB}\bbX^{\vee} \to \bbX_{\bpi}]}$ is fully
faithful.

To compute $\mycal{P}*\mycal{Q}$ consider the triple product
$\bbX_{\bpi}\times_{\bbD} {}_{\bB}\bbX^{\vee} \times_{\bbD}
\bbX_{\bpi}^{\op{op}}$, which for the purposes of handling 
left-right modules we
will view as the product $X \times {}_{\bB}\bbX^{\vee} \times X$
equipped with the structure sheaf
\[
{\mathcal O}_{\bbX_{\bpi}\times_{\bbD} {}_{\bB}\bbX^{\vee} \times_{\bbD}
\bbX_{\bpi}} = p_{1}^{-1} \Ah{\bpi} 
\widehat{\otimes}_{{}_{\Ch}} p_{2}^{-1} {\mathcal
  O}_{{}_{\bB}\bbX^{\vee}} \widehat{\otimes}_{{}_{\Ch}}  p_{3}^{-1} 
\Ah{\bpi}^{\op{op}}.
\]
Where we write $p_{1} : X \times {}_{\bB}\bbX^{\vee} \times X \to X$,
$p_{12} : X \times {}_{\bB}\bbX^{\vee} \times X \to X\times
{}_{\bB}\bbX^{\vee}$, etc. for the projections onto the corresponding
spaces or stacks. We will also write $\bbp_{1} :
\bbX_{\bpi}\times_{\bbD} {}_{\bB}\bbX^{\vee} \times_{\bbD} 
\bbX_{\bpi}^{\op{op}}
\to \bbX_{\bpi}$, $\bbp_{12} : \bbX_{\bpi}\times_{\bbD}
    {}_{\bB}\bbX^{\vee} \times_{\bbD} \bbX_{\bpi}^{\op{op}}
\to \bbX_{\bpi}\times_{\bbD} {}_{\bB}\bbX^{\vee}$ for the
corresponding maps of ringed objects.

Now, according to section~\ref{sssec:groupstr} the difference map $d
  : X\times X \to X$,
  $d(x,y) = x-y$ for the additive structure on the complex torus $X$,
  induces a natural map of ringed spaces
\[
\bbd_{\bpi} : \bbX_{\bpi}\times_{\bbD} \bbX_{\bpi}^{\op{op}} \to \bbX.
\]
After inserting ${}_{\bB}\bbX^{\vee}$ as a middle factor in the triple
product $\bbX_{\bpi}\times_{\bbD} {}_{\bB}\bbX^{\vee} \times_{\bbD}
\bbX_{\bpi}^{\op{op}}$ we will write $\bbp_{1-3,2}$ for the composition map
\[
\xymatrix@C+5pc{
{\bbX_{\bpi}\times_{\bbD} {}_{\bB}\bbX^{\vee} \times_{\bbD}
  \bbX_{\bpi}^{\op{op}}} 
\ar[r]^-{(\bbd_{\bpi}(\bbp_{1}(\bullet),\bbp_{3}(\bullet)),
  \bbp_{2}(\bullet))} 
\ar[dr]_-{\bbp_{1-3,2}} & \bbX\times_{\bbD} {}_{\bB}\bbX^{\vee}
\ar[d]^-{\bbp_{1}\times \bbr} \\
& \bbX\times_{\bbD} \bbX^{\vee}
}
\]
where $\bbr : {}_{\bB}\bbX^{\vee} \to \bbX$ is the natural structure
map. Being a morphism of ringed objects the map $\bbp_{1-3,2} :
{\bbX_{\bpi}\times_{\bbD} {}_{\bB}\bbX^{\vee} \times_{\bbD}
  \bbX_{\bpi}^{\op{op}}} \to \bbX\times_{\bbD} \bbX^{\vee}$ is given as a pair
$\bbp_{1-3,2} = (p_{1-3,2},p_{1-3,2}^{\sharp})$ where
$p_{1-3,2} : X\times {}_{\bB}\bbX^{\vee} \times X \to X\times
X^{\vee}$ and a morphism of sheaves of algebras
\[
\xymatrix@1@M+1pc{
p_{1-3,2}^{\sharp} : \hspace{-2pc} &
p_{1-3,2}^{-1}{\mathcal O}_{X\times X^{\vee}} \ar[r] &  p_{1}^{-1} \Ah{\bpi}
\widehat{\otimes}_{{}_{\Ch}} p_{2}^{-1} {\mathcal
  O}_{{}_{\bB}\bbX^{\vee}} 
\widehat{\otimes}_{{}_{\Ch}} p_{3}^{-1} \Ah{\bpi}^{\op{op}}.
}
\]
As usual, given a coherent sheaf $\mathcal{M}$ on $X\times X^{\vee}$
we define its pullback via $\bbp_{1-3,2}$ to be the sheaf on $X\times
{}_{\bB}\bbX^{\vee} \times X$ given by:
\[
\bbp_{1-3,2}^{*}\mathcal{M} := \quad \left(p_{1-3,2}^{-1} \mathcal{M} \right)
\otimes_{p_{1-3,2}^{-1}{\mathcal O}_{X\times X^{\vee}}} 
\left( p_{1}^{-1} \Ah{\bpi}
\widehat{\otimes}_{{}_{\Ch}} p_{2}^{-1} {\mathcal
  O}_{{}_{\bB}\bbX^{\vee}} \widehat{\otimes}_{{}_{\Ch}} 
p_{3}^{-1} \Ah{\bpi}^{\op{op}}\right).
\]
Since ${\mathcal O}_{X\times X^{\vee}}$ is commutative, the sheaf
$\bbp_{1-3,2}^{*}\mathcal{M}$ will have a natural structure of a left
$p_{1}^{-1}\Ah{\bpi}$-module and a right $p_{3}^{-1}\Ah{\bpi}$ module.

With this notation we have

\begin{lem} \label{lem:PQ} Let $\mathcal{P} \to X\times X^{\vee}$ be
  the classical Poincar\'{e} bundle for the pair of dual complex tori
  $(X,X^{\vee})$. On the triple product
  $X \times ({}_{\bB}\bbX^{\vee}) \times X$, there is a natural
  isomorphism of sheaves
\begin{equation} \label{eq:isomPQ}
p_{12}^{-1}\mycal{P}\otimes_{p_{2}^{-1} {\mycal
    O}_{{}_{\bB}\bbX^{\vee}}} p_{23}^{-1}\mycal{Q} \cong \bbp_{1-3,2}^{*}
    \mathcal{P}[[\hbar]].
\end{equation}
which is also an isomorphism in
$\lrmod{p_{1}^{-1}\Ah{\bpi}}{p_{3}^{-1}\Ah{\bpi}}$.
\end{lem}

\noindent
{\bf Proof.} First observe  that the sheaves
appearing in the two sides of the identity
\eqref{eq:isomPQ}  can be all specified via factors of
automorphy. Thus the question of proving the lemma reduces to
computing and comparing the factors of automorphy of the left hand
side and the right hand side of \eqref{eq:isomPQ}. 

As explained in section \ref{ss:poincare} the Poincar\'{e} sheaf
$\mycal{P}$ is defined by a $\Lambda\times \bGa$ factor of  automorphy
$\phi_{\mycal{P}}$  
on the vector space $V\times \overline{V}^{\vee}$ taking values in the
invertible global sections of the sheaf of algebras
$p_{1}^{-1}\mycal{A}_{\bpi}
\widehat{\otimes}_{{\mathbb C}} p_{2}^{-1}{\mathcal
  O}_{\overline{V}^{\vee}}$. Explicitly
\begin{equation} \label{eq:phiP}
\xymatrix@R-2pc{
\phi_{\mycal{P}} : \hspace{-2pc} &  \Lambda\times \bGa \ar[r] &
\Gamma\left(V\times 
\overline{V}^{\vee}, \; \left(p_{1}^{-1}
\mycal{A}_{\bpi}\widehat{\otimes}_{{\mathbb C}} 
p_{2}^{-1}{\mathcal 
  O}_{\overline{V}^{\vee}}\right)^{\times}\right) \\
& (\lambda; (\xi,z)) \ar[r] & \left( (v,l) \mapsto z e^{\pi\langle \xi,
\lambda \rangle}  
  e^{\pi\overline{\langle \xi, v \rangle }} e^{\pi\langle l, \lambda
    \rangle}\right). 
}
\end{equation}
Similarly 
\begin{equation} \label{eq:phiQ}
\xymatrix@R-2pc{ \phi_{\mycal{Q}} : \hspace{-2pc} & \bGa\times \Lambda
\ar[r] & \Gamma\left(\overline{V}^{\vee}\times V, \;
\left(p_{1}^{-1}{\mathcal
O}_{\overline{V}^{\vee}} \widehat{\otimes}_{{\mathbb C}}
p_{2}^{-1}\mycal{A}_{\bpi}\right)^{\times}\right)  \\ & (\lambda;
(\xi,z)) \ar[r] & \left(  
(l,v) \mapsto z^{-1} e^{- \pi\langle \xi, \lambda \rangle} e^{- \pi
\langle l, \lambda \rangle}e^{\pi \overline{\langle \xi, v \rangle
}}\right).  }
\end{equation}
Let 
\[
\phi : \Lambda \times \bGa \times \Lambda \to \Gamma\left( V\times
\overline{V}^{\vee} \times V, \; 
\op{Aut}\left(p_{1}^{-1}\Ah{\bpi}\widehat{\otimes}_{{}_{\Ch}}
p_{2}^{-1}\mathcal{O}_{\overline{V}^{\vee}}[[\hbar]] 
\widehat{\otimes}_{{}_{\Ch}}
p_{3}^{-1} \Ah{\bpi}^{\op{op}} \right)\right)
\] 
denote the factor of automorphy for the left-right module
$p_{12}^{-1}\mycal{P}\otimes_{p_{2}^{-1} {\mycal
    O}_{{}_{\bB}\bbX^{\vee}}} p_{23}^{-1}\mycal{Q}$. The formulas
\eqref{eq:phiP} and \eqref{eq:phiQ} now combine in an 
explicit expression for $\phi$.  If $f$ is a local section of 
$p_{1}^{-1}\Ah{\bpi}\otimes_{\Ch}
p_{2}^{-1}\mathcal{O}_{\overline{V}^{\vee}}[[\hbar]] \otimes_{\Ch}
p_{3}^{-1} \Ah{\bpi}$, then 
\[
[\phi(\lambda,(\xi,z),\mu)(f)](v,x,w) = [p_{2}^{-1}(e^{\pi\langle
\bullet +\xi, \lambda- \mu\rangle}) p_{3}^{-1}
(e^{-\pi\overline{\langle \xi, \bullet \rangle}}) \star  f
\star p_{1}^{-1}(e^{\pi\overline{\langle \xi, \bullet \rangle}})](v,x,w).
\] 
We now compute the factor of automorphy for 
$\bbp_{1-3,2}^{*}\mathcal{P}[[\hbar]]$.  
Let $\phi_{\mathcal{P}}$ denote the factor of automorphy for 
the classical 
Poincar\'{e} sheaf $\mathcal{P} \to X \times X^{\vee}$

\begin{equation} \label{eq:phiPclassical}
\xymatrix@R-2pc{
\phi_{\mathcal{P}} : \hspace{-2pc} &  \Lambda\times \Lambda^{\vee} \ar[r] &
\Gamma\left(V\times 
\overline{V}^{\vee}, \; 
\left(p_{1}^{-1}{\mathcal O}_{V}\widehat{\otimes}_{{\mathbb C}} 
p_{2}^{-1}{\mathcal 
  O}_{\overline{V}^{\vee}}\right)^{\times}\right) \\
& (\lambda, \xi) \ar[r] & \left( (v,l) \mapsto  e^{\pi\langle l+ \xi,
\lambda \rangle}  
  e^{\pi\overline{\langle \xi, v \rangle }} \right). 
}
\end{equation}
\noindent
Consider the map 
\begin{equation} \label{eq:diffoncover}
\xymatrix@R-2pc{
{\sf{diff}} : \hspace{-2pc} & V \times \overline{V}^{\vee} \times V \to
V \times \overline{V}^{\vee} \\ 
&(v, x, w) \mapsto (v-w,x) 
}
\end{equation}
\noindent
The factor of automorphy $\psi$ for $\bbp_{1-3,2}^{*}\mathcal{P}[[\hbar]]$ is 
a map 
\[
\psi : \Lambda \times \bGa \times \Lambda \to \Gamma\left( V\times
\overline{V}^{\vee} \times V, \; 
\cAut \left(p_{1}^{-1}\Ah{\bpi}\widehat{\otimes}_{{}_{\Ch}}
p_{2}^{-1}\mathcal{O}_{\overline{V}^{\vee}}[[\hbar]] 
\widehat{\otimes}_{{}_{\Ch}}
p_{3}^{-1} \Ah{\bpi}^{\op{op}} \right)\right), 
\] 
where the automorphism $\psi(\lambda; (\xi, z) ; \mu)$ is given by
right-left multiplication by the invertible section  
\[
p_{1-3,2}^{\sharp} (p_{1-3,2}^{-1} \phi_{\mathcal{P}} \circ
{\sf{diff}})(\lambda; (\xi, z) ; \mu)
\]
of the sheaf of algebras $p_{1}^{-1}\Ah{\bpi}\widehat{\otimes}_{{}_{\Ch}}
p_{2}^{-1}\mathcal{O}_{\overline{V}^{\vee}}[[\hbar]] 
\widehat{\otimes}_{{}_{\Ch}}
p_{3}^{-1} \Ah{\bpi}$.
Since the map  $p_{1-3,2}^{\sharp}$ is defined
in terms  of the coproduct structure of section~\ref{ss:poincare} we
compute explicitly:
\begin{equation} \label{eq:p132}
\begin{split}
[p_{1-3,2}^{\sharp} (p_{1-3,2}^{-1} \phi_{\mathcal{P}} & \circ
{\sf{diff}})(\lambda; (\xi, z) ; \mu)] (v,l,w)  = \\
& = \left( 1\otimes e^{\pi\langle l+ \xi,
\lambda - \mu \rangle} \otimes 1\right)\star \left(  
  e^{\pi\overline{\langle  \xi, v\rangle}}\otimes 1 \otimes
  1\right)\star \left(1 \otimes 1
      \otimes e^{-\pi\overline{\langle  \xi, w \rangle }}\right) \\
& =  
      e^{\pi\overline{\langle  \xi, v\rangle}} \otimes e^{\pi\langle l+ \xi,
\lambda - \mu \rangle} \otimes 
      e^{-\pi\overline{\langle  \xi, w \rangle }} 
\end{split}
\end{equation}
Thus $\psi(\lambda; (\xi, z) ; \mu)$ is the automorphism of
$p_{1}^{-1}\Ah{\bpi}\widehat{\otimes}_{{}_{\Ch}}
p_{2}^{-1}\mathcal{O}_{\overline{V}^{\vee}}[[\hbar]]
\widehat{\otimes}_{{}_{\Ch}} p_{3}^{-1} \Ah{\bpi}^{\op{op}}$ given by
$p_{1}^{-1}e^{\pi\overline{\langle \xi, \bullet \rangle}}$ acting by a
right multiplication, $p_{3}^{-1} e^{-\pi\overline{\langle \xi,
\bullet\rangle}}$ acting by a left multiplication, and the function
$p_{2}^{-1}e^{\pi\langle \bullet + \xi, \lambda - \mu \rangle}$ acting
in the unambiguous, obvious way.  Since this agrees precisely with the
description of $\phi$, we are done.  \ \hfill $\Box$

\medskip

To finish the computation of the convolution $\mycal{P} * \mycal{Q}$
we use the fiber-product diagram:
\[
\xymatrix{
\bbX_{\bpi} \times_{\bbD} {}_{\bB}\bbX^{\vee} \times_{\bbD}
\bbX_{\bpi} \ar[r]^-{\bbp_{1-3,2}} \ar[d]_-{p_{13}} & \bbX_{0}\times_{\bbD}
    {}_{0}\bbX^{\vee} \ar[d]^-{p_{1}} \\
\bbX_{\bpi} \times_{\bbD}
\bbX_{\bpi} \ar[r]_-{\bbp_{1-2}} & \bbX_{0}
}
\]
where $\bbX_{0} = X\times \bbD$ denotes the trivial formal deformation
of $X$ and ${}_{0}\bbX^{\vee}$ denotes the trivial
$\mathcal{O}^{\times}$-gerbe on the space $\bbX^{\vee} =
X^{\vee}\times \bbD$. 

Alternatively, viewing the maps $\bbp_{1-3,2} =
(p_{1-3,2},p_{1-3,2}^{\sharp})$, etc. as morphisms of ringed spaces we
can apply the base change property and Lemma~\ref{lem:PQ} to conclude that 
\[
\begin{split}
\mycal{P}* \mycal{Q} & := Rp_{13*} \left( p_{12}^{-1}\mycal{P}
 \otimes_{ p_{2}^{-1} \mathcal{O}_{{}_{\bB}\bbX^{\vee}}}
 p_{23}^{-1}\mycal{Q}  \right) \\
 & =  Rp_{13*}\bbp_{1-3,2}^{*} \mathcal{P}[[\hbar]] \\
& = \bbp_{1-2}^{*}\left(Rp_{1*}\mathcal{P}[[\hbar]]\right) \\
& = \bbp_{1-2}^{*}\left( \mathcal{O}_{0}[[\hbar]]\right)[-g],
\end{split}
\]
where $\mathcal{O}_{0}$ denotes the skyscraper sheaf on $X$ supported
at the origin $0 \in X$. In particular, the identity 
\[
\phi_{\mycal{P}}^{[ _{\bB}\bbX^{\vee} \to
    \bbX_{\bpi}]} \circ \phi_{\mycal{Q}}^{[\bbX_{\bpi}\to
    _{\bB}\bbX^{\vee}]}  \cong \op{id}_{\sff{D}^{*}( {
    \bbX_{\bpi}})}[-g]
\]
will follow immediately from the convolution property
\eqref{eq:convolutions} and the following

\begin{lem} \label{lem:diagonal}
\[
\phi_{\bbp_{1-2}^{*}\mathcal{O}_{0}[[\hbar]]}^{[\bbX_{\bpi} \to
    \bbX_{\bpi}]} = \op{id}_{\sff{D}^{*}( {
    \bbX_{\bpi}})}.
\]
\end{lem}
\noindent
{\bf Proof.} Consider the sheaf
$\bbp_{1-2}^{*}\mathcal{O}_{0}[[\hbar]]$ on the topological space
$X\times X$. By construction it is naturally a left
$p_{1}^{-1}\Ah{\bpi}$ module and a right
$p_{2}^{-1}\Ah{\bpi}$ module. The element
$\bbp_{1-2}^{*}\mathcal{O}_{0}[[\hbar]] \in 
\lrmod{p_{1}^{-1}\Ah{\bpi}}{p_{2}^{-1}\Ah{\bpi}}$ is easy to
compute. Recall that the sheaf
$\bbp_{1-2}^{*}\mathcal{O}_{0}[[\hbar]]$ is defined as the tensor
product 
\[
\bbp_{1-2}^{*}\mathcal{O}_{0}[[\hbar]] = \left( p_{1}^{-1}
\Ah{\bpi}\widehat{\otimes}_{\Ch} p_{2}^{-1}\Ah{\bpi}^{\op{op}}\right)
\otimes_{p_{1-2}^{-1}\mathcal{O}_{X}[[\hbar]]} \left(
p_{1-2}^{-1}\mathcal{O}_{0}[[\hbar]] \right)
\]
of $p_{1-2}^{-1}\mathcal{O}_{0}[[\hbar]]$ with the sheaf of
algebras 
$p_{1}^{-1}\Ah{\bpi}\otimes_{\Ch} p_{2}^{-1}\Ah{\bpi}^{\op{op}}$, where
the tensor product is taken over the algebra
$p_{1-2}^{-1}\mathcal{O}[[\hbar]]$ via the coproduct homomorphism
\[
p_{1-2}^{\sharp} : p_{1-2}^{-1}\mathcal{O}_{X}[[\hbar]] \to
p_{1}^{-1}\Ah{\bpi}\widehat{\otimes}_{\Ch} p_{2}^{-1}\Ah{\bpi}^{\op{op}} 
\]
defined in section \ref{sssec:groupstr}. In these terms the left-right
module structure on $\bbp_{1-2}^{*}\mathcal{O}_{0}[[\hbar]]$ arises
from the left multiplication of $p_{1}^{-1}\Ah{\bpi}
\widehat{\otimes}_{\Ch}
p_{2}^{-1}\Ah{\bpi}^{\op{op}}$ on itself.

\

\medskip

\noindent
Consider the commutative diagram of topological spaces 
\[
\xymatrix@C+5pc{
X
\ar[r]^-{\Delta}
\ar[d]_-{} &  X\times X
\ar[d]^-{p_{1-2}} \\
{0} \ar[r]^-{} &  X
}
\]
From this diagram, we deduce $p_{1-2}^{-1} \mathcal{O}_{0}[[\hbar]]
\cong \Delta_{*} \mathbb{C}[[\hbar]]$. Below we will use this
identification to argue that the element 
$\bbp_{1-2}^{*}\mathcal{O}_{0}[[\hbar]] \in
\lrmod{p_{1}^{-1}\Ah{\bpi}}{p_{2}^{-1}\Ah{\bpi}}$ is given as
\begin{equation} \label{eq:basechange}
\bbp_{1-2}^{*}\mathcal{O}_{0}[[\hbar]]  \cong  \Delta_{*}\Ah{\bpi}.
\end{equation}
Here we view $\Delta_{*}\Ah{\bpi}$ as a sheaf on $X\times X$ supported on the
diagonal and equipped with natural left-right module structure. Namely
the left $p_{1}^{-1}\Ah{\bpi}$-module structure corresponding to left
multiplication by elements in $\Delta^{-1}p_{1}^{-1} \Ah{\bpi} =
\Ah{\bpi}$ and the right $p_{2}^{-1}\Ah{\bpi}$-module structure
corresponding to right 
multiplication by elements in $\Delta^{-1}p_{2}^{-1} \Ah{\bpi} =
\Ah{\bpi}$.  

The isomorphism claimed in equation \ref{eq:basechange} is given by the
following mutually inverse maps in
$\lrmod{p_{1}^{-1}\Ah{\bpi}}{p_{2}^{-1}\Ah{\bpi}}$
\begin{equation} \label{eq:firstmap}
\xymatrix@R-2pc{
\left(
p_{1}^{-1}\Ah{\bpi}\widehat{\otimes}_{\Ch} 
p_{2}^{-1}\Ah{\bpi}^{\op{op}} \right)
\otimes_{p_{1-2}^{-1}\mathcal{O}_{X}[[\hbar]]} 
\left(\Delta_{*}
\mathbb{C}[[\hbar]]\right)
 \ar[r] &  \Delta_{*}\Ah{\bpi} \\
a \otimes b \otimes \left(\Delta_{*}c\right) \ar@{|->}[r] & 
\Delta_{*}\left((\Delta^{-1}(b)) \star (\Delta^{-1}(a))\star c \right) 
}
\end{equation}
and 
\begin{equation} \label{eq:secondmap}
\xymatrix@R-2pc{
\Delta_{*}\Ah{\bpi}   \ar[r] &  \left(
p_{1}^{-1}\Ah{\bpi}\widehat{\otimes}_{\Ch} 
p_{2}^{-1}\Ah{\bpi}^{\op{op}} \right)
\otimes_{p_{1-2}^{-1}\mathcal{O}_{X}[[\hbar]]}
\left(\Delta_{*}
\mathbb{C}[[\hbar]]\right) \\
\Delta_{*}f \ar@{|->}[r] & 1 \otimes p_{1}^{-1} (f) \otimes 1 
}
\end{equation}
The composition of \eqref{eq:secondmap} followed by
\eqref{eq:firstmap} tautologically gives the identity. Composing in
the reverse order also gives the identity. Indeed, the composed map is
\[
a \otimes b \otimes \left(\Delta_{*}c\right) \quad  \mapsto \quad 
p_{1}^{-1}\left( \Delta^{-1}(b)\star \Delta^{-1}(a)\star c \right)
\otimes 1 \otimes 1.
\]
Consider now the element
\[
a\otimes 1 \otimes \Delta_{*}c \in p_{1}^{-1}\Ah{\bpi}\widehat{\otimes}_{\Ch}
p_{2}^{-1}\Ah{\bpi}^{\op{op}}\otimes_{p_{1-2}^{-1}\mathcal{O}_{X}[[\hbar]]}
\left( \Delta_{*}\Ch \right).
\]
Taking into account that $p_{1}^{-1}\left(\Delta^{-1}(a)\right) = a$,
$p_{2}^{-1}\left(\Delta^{-1}(b)\right) = b$, and  using the definition of
the left-right module structure on $p_{1}^{-1}\Ah{\bpi}\widehat{\otimes}_{\Ch}
p_{2}^{-1}\Ah{\bpi}^{\op{op}}\otimes_{p_{1-2}^{-1}\mathcal{O}_{X}[[\hbar]]}
\left( \Delta_{*}\Ch \right)$, we get
\[
\begin{split}
a\otimes b \otimes \Delta_{*}c & =
p_{1}^{-1}\left(\left(\Delta^{-1}(a)\right)\star c \right) \otimes b
\otimes 1 \\
& = \left( p_{1}^{-1}\left(\left(\Delta^{-1}(a)\right)\star c \right)
\otimes 1 
\otimes 1\right)\cdot b.
\end{split}
\]
Moreover, for any $x \in p_{1}^{-1}\Ah{\bpi}$, we have the identity
\begin{equation} \label{eq:diagid}
(x\otimes 1\otimes 1)\cdot b = p_{1}^{-1}(\Delta^{-1}(b))\cdot
  (x\otimes 1\otimes 1),
\end{equation}
valid in the left-right module $p_{1}^{-1}\Ah{\bpi}\widehat{\otimes}_{\Ch}
p_{2}^{-1}\Ah{\bpi}^{\op{op}}\otimes_{p_{1-2}^{-1}\mathcal{O}_{X}[[\hbar]]}
\left( \Delta_{*}\Ch \right)$. This identity follows immediately from
the definition of the coproduct map $p_{1-2}^{\sharp}$ and from the
fact that this module is supported on the diagonal in $X\times X$. 

Combining the identity \eqref{eq:diagid} with the definition of the left
$p_{1}^{-1}\Ah{\bpi}$ action on \linebreak $p_{1}^{-1}\Ah{\bpi}
\widehat{\otimes}_{\Ch}
p_{2}^{-1}\Ah{\bpi}^{\op{op}}\otimes_{p_{1-2}^{-1}\mathcal{O}_{X}[[\hbar]]}
\left( \Delta_{*}\Ch \right)$, yields
\[
(x\otimes 1\otimes 1)\cdot b = p_{1}^{-1}(\Delta^{-1}(b))\cdot
  (x\otimes 1\otimes 1) = \left( p_{1}^{-1}(\Delta^{-1}(b))\star x
  \right) \otimes 1 \otimes 1,
\]
and so  $a\otimes b \otimes \Delta_{*}c = p_{1}^{-1}\left(
\Delta^{-1}(b)\star \Delta^{-1}(a) \star c\right)\otimes 1 \otimes
1$.

Now that we have the isomorphism 
$\bbp_{1-2}^{*}\mathcal{O}_{0}[[\hbar]]  \cong  \Delta_{*}\Ah{\bpi}$ we can easily 
compute that for any element  $\mycal{M} \in \sff{D}^{*}(\bbX_{\bpi})$

\[
Rp_{1*}(\bbp_{1-2}^{*} \mathcal{O}_{0}[[\hbar]] \otimes_{p_{2}^{-1}\Ah{\bpi}} p_{2}^{-1} \mycal{M} )
\cong Rp_{1*}(R\Delta_{*} \Ah{\bpi} \otimes_{p_{2}^{-1}\Ah{\bpi}} p_{2}^{-1} \mycal{M} )
\]
\[\cong Rp_{1*} R\Delta_{*}(\Ah{\bpi} \otimes_{\Delta^{-1} p_{2}^{-1} 
\Ah{\bpi}} \Delta^{-1} p_{2}^{-1} \mycal{M} )
\]
in other words we have 
\[
Rp_{1*}(\bbp_{1-2}^{*} \mathcal{O}_{0}[[\hbar]] \otimes_{p_{2}^{-1}\Ah{\bpi}} p_{2}^{-1} \mycal{M} )
\cong
Rp_{1*} R\Delta_{*}(\Ah{\bpi} \otimes_{\Ah{\bpi}} \mycal{M} ) 
\cong Rp_{1*} R\Delta_{*}( \mycal{M} ) \cong \mycal{M}
\]
This completes the proof of the lemma.
\ \hfill $\Box$  

\

\bigskip

In the opposite direction we must verify the identity
\[
\phi_{\mycal{Q}}^{[\bbX_{\bpi}\to {_{\bB}\bbX^{\vee}}]}
\circ \phi_{\mycal{P}}^{[ { _{\bB}\bbX^{\vee}} \to
    \bbX_{\bpi}]}  \cong \op{id}_{\sff{D}^{*}( {
    _{\bB}\bbX^{\vee}},-1)}[-g].
\]
By Lemma~\ref{lem:diagonal} we have
\[
\phi_{\mycal{P}}^{[ { _{\bB}\bbX^{\vee}} \to
    \bbX_{\bpi}]}  \circ \phi_{\mycal{Q}}^{[\bbX_{\bpi}\to
    {_{\bB}\bbX^{\vee}}]}  \cong \op{id}_{\sff{D}^{*}( {
    \bbX^{\bpi}})}[-g],
\]
and so  $\phi_{\mycal{P}} := \phi_{\mycal{P}}^{[ {
    _{\bB}\bbX^{\vee}} \to 
    \bbX_{\bpi}]}$ is essentially surjective. Therefore, if we can
    check that $\phi_{\mycal{P}}$ is fully faithful we can conclude
    that $\phi_{\mycal{P}}$ is an equivalence and that $\phi_{\mycal{Q}}[g]$
    is the inverse. 

\begin{lem} \label{lem:phiPfaithful}
The functor $\phi_{\mycal{P}}^{[ { _{\bB}\bbX^{\vee}} \to
    \bbX_{\bpi}]}$ is fully faithful. 
\end{lem}
To argue that $\phi_{\mycal{P}}$ is fully faithful we will first
 identify an orthogonal spanning class of objects for the category
 $\sff{D}^{*}( {{}_{\bB}\bbX^{\vee}},-1)$, and then check that
 $\phi_{\mycal{P}}$ satisfies the Bondal-Orlov faithfulness criterion
 \cite{BO,bridgeland} on the spanning class. To apply the criterion we
 will need to know that the functor $\phi_{\mycal{P}}$ has left and
 right adjoints. For this we only need to note that $\phi_{\mycal{P}}$
 is given as the composition of a pullback, a tensoring with
 $\mycal{P}$ and a pushforward. Since the pullback and the pushforward
 have left and right adjoints and the tensoring with $\mycal{P}$ has a
 left and right adjoint given by the tensoring with $\mycal{Q}$ we get
 that $\phi_{\mycal{P}}$ has both left and right adjoints.  This puts
 us in a position to apply the Bondal-Orlov criterion.  We proceed in
 several steps:

\

\bigskip 

\noindent
{\bfseries Step 1.} Consider the category $\sff{D}^{*}(
{{}_{\bB}\bbX^{\vee}},-1)$ of weight $(-1)$ sheaves on the gerbe
${}_{\bB}\bbX^{\vee}$. We have a natural collection of 
objects in this category labeled by the points of the torus $X^{\vee}$.

Indeed, suppose $s \in X^{\vee}$ is a point in the torus $X^{\vee}$
and let $\mathfrak{s} : \bbD \to \bbX^{\vee}$ be the constant section
of $\bbX^{\vee} = X^{\vee}\times \bbD \to \bbD$ passing through
$s$. The first observation is that the pullback of the gerbe
${}_{\bB}\bbX^{\vee}$ by $\mathfrak{s}$ is a trivial ${\mathcal
O}^{\times}$-gerbe on $\bbD$. Indeed, the gerbe ${}_{\bB}\bbX^{\vee}$
was defined as a global quotient $[\overline{\bbV}^{\vee}/\bGamma]$
with $\bGamma$ acting through the projection $\bGamma \to
\bLambda^{\vee}$ and the natural free action of $\bLambda^{\vee}$ on
$\overline{\bbV}^{\vee}$. In particular, for any space $S$ and any
map $f : S \to X^{\vee}$ for which $S\times_{X^{\vee}}
\overline{V}^{\vee}$ is a trivial $\Lambda^{\vee}$ cover of $S$, the
$f$-pullback of ${}_{\bB}\bbX^{\vee}$ will be a trivial ${\mathcal
  O}^{\times}$-gerbe. This implies that for any weight $k \in
\mathbb{Z}$ we can view the sheaf
$\mathfrak{s}_{*}{\mathcal O}_{\bbD} = \mathcal{O}_{s \times \bbD}$
as a sheaf on ${}_{\bB}\bbX^{\vee}$ of pure weight $k$.

In fact, the presentation of ${}_{\bB}\bbX^{\vee}$ as a global
quotient provides a canonical way of endowing
$\mathfrak{s}_{*}{\mathcal O}_{\bbD}$ with the structure of a weight
$k$ sheaf on the gerbe. To that end, consider the universal covering
map $\pi : \overline{V}^{\vee} \to X^{\vee}$. The preimage $F_{s} :=
\pi^{-1}(s) \subset \overline{V}^{\vee}$ is a $\Lambda^{\vee}$-orbit
in $\overline{V}^{\vee}$, and the pullback gerbe 
\[
\mathfrak{s}^{*}
({}_{\bB}\bbX^{\vee}) :={}_{\bB}\bbX^{\vee}
\times_{\bbX^{\vee},\mathfrak{s}} \bbD 
\]
 is naturally realized as the global quotient
$[(F_{s}\times \bbD)/\bGamma]$. Looking at $\bbD$-points it is clear
that in order equip $\mathfrak{s}_{*}\mathcal{O}_{\bbD}$ with the
structure of a weight $(-1)$-sheaf on ${}_{\bB}\bbX^{\vee}$, it
suffices to describe a $\Gamma$ equivariant structure on
$\mathcal{O}_{F_{s}}[[\hbar]]$ in which every central element $z \in
\Ch^{\times} \subset \Gamma$ acts as multiplication by $z^{-1}$ viewed
as a section in $\mathcal{O}_{F_{s}}[[\hbar]]$.

To achieve this we let as before $c :
\Lambda^{\vee}\times \Lambda^{\vee} \to \Ch^{\times}$ denote the cocycle
defining the group $\Gamma$. Since $c$ was given as the exponential
of a multiple of the ${\mathbb R}$-linear map $B : \wedge^{2}
\overline{V}^{\vee} \to {\mathbb C}$, we can use ${\mathbb
  R}$-linearity to extend $c$ to a multiplicative map $\tilde{c} :
\Lambda^{\vee}\otimes \overline{V}^{\vee} \to \Ch^{\times}$. Now,
suppose $(\xi,z) \in \Gamma$ and let $f \in
\mathcal{O}_{F_{s}}[[\hbar]]$. Since $F_{s}$ is a discrete set of 
points, we can write $f$ as a collection $\{ f_{w} \}_{w \in F_{s}}$
  with $f_{w} \in \mathcal{O}_{F_{s},w}[[\hbar]] = \Ch$. With this
  notation  $(\xi,z)$ gives rise to an automorphism
\[
\rho_{(\xi,z)} : H^{0}(F_{s},\mathcal{O}_{F_{s}}[[\hbar]]) \to 
H^{0}(F_{s},\mathcal{O}_{F_{s}}[[\hbar]]),
\]
where for each $w \in F_{s}$ we define
\[
\rho_{(\xi,z)} = \left\{ \left( \rho_{(\xi,z)}f \right)_{w} \right\}_{w
  \in F_{s}}, \quad \text{with} \quad \left( \rho_{(\xi,z)}f
  \right)_{w} := z^{-1}\tilde{c}(w,\xi) f_{w-\xi}. 
\]
To check that $\rho$ is a $\Gamma$-action we compute
\[
\begin{split}
\left(\rho_{(\xi',z')}\rho_{(\xi,z)}f\right)_{w} & =
z^{'-1}\tilde{c}(w,\xi')\left( \rho_{(\xi,z)}f \right)_{w - \xi'} \\
& = z^{'-1}z^{-1}\tilde{c}(w,\xi')\tilde{c}(w - \xi',\xi) f_{w - \xi'
  - \xi} \\
& = (z'z)^{-1}\tilde{c}(-\xi',\xi)\tilde{c}(w,\xi+\xi') f_{w -
  \xi' - \xi}. 
\end{split}
\]
On the other hand $(\xi',z')\cdot (\xi, z) = (\xi+\xi',zz'c(\xi',\xi))$
and so
\[
\left(\rho_{(\xi',z')\cdot (\xi, z)}f \right)_{w} =
  (z'z)^{-1}\tilde{c}(\xi',\xi)^{-1} \tilde{c}(w,\xi+\xi') f_{w -
  \xi' - \xi}.
\]
taking into account the fact that $c$ is bilinear on $\Lambda^{\vee}$ we
conclude that $\rho_{(\xi',z')\cdot (\xi, z)} =
\rho_{(\xi',z')}\rho_{(\xi, z)}$. 

For future reference we will denote the weight $(-1)$-sheaf given by
$\rho$ by $(\mathfrak{s},-1)_{*}\mathcal{O}_{\bbD}$. Clearly this
generalizes to all weights $k \in \mathbb{Z}$ yielding weight $k$
sheaves $(\mathfrak{s},k)_{*}\mathcal{O}_{\bbD}$ on
${}_{\bB}\bbX^{\vee}$. In fact we have a pair of adjoint functors
$(L(\mathfrak{s},k)^{*},(\mathfrak{s},k)_{*})$ between the categories
$D^{*}(\bbD)$ and $\sff{D}^{*}({}_{\bB}\bbX^{\vee}, k)$,
defined for all integers $k$. 

The existence of such functors is a
basic fact about any morphism between a space and an
$\mathcal{O}^{\times}$ gerbe. Indeed, suppose $S$ is a space and
$\mycal{T}$ is an $\mathcal{O}^{\times}$ gerbe on a space $T$. Suppose
we are given a morphism of spaces $f : S \to T$ with the property that 
$f^{*}\mycal{T}$ is trivializable. Choosing a trivialization
$f^{*}\mycal{T} \cong {}_{0}S$ we get a well defined pair of adjoint
functors $(Lf^{*},Rf_{*})$ between $D^{*}(S)$ and
$\sff{D}^{*}(\mycal{T})$. By construction these functors are
compatible with the weight decompositions so we get get induced
adjoint pairs between the corresponding weight pieces. These can in
turn be combined with the canonical equivalences $D^{*}(S) =
\sff{D}^{*}({}_{0}S,k)$ which are also defined for all $k$.
This results in adjoint pairs of functors
\begin{equation} \label{eq-adjoint}
\xymatrix@1@M+1pc@C+1pc{
D^{*}(S) \ar@<1ex>@{->}[r]^-{R(f,k)_{*}} &
\sff{D}^{*}(\mycal{T},k) \ar@<1ex>@{->}[l]^-{L(f,k)^{*}} 
}
\end{equation}
which we will frequently use below. Note that the construction of
\eqref{eq-adjoint} depends on the choice of trivialization of the
gerbe $f^{*}\mycal{T}$. In the particular case of the map
$\mathfrak{s} : \bbD \to \bbX^{\vee}$, we used a
special trivialization constructed out of the quotient presentation of
${}_{\bB}\bbX^{\vee}$. This trivialization is precisely encoded in the
map $\tilde{c}$ used above. 

Consider now the collection of objects $\left\{
(\mathfrak{s},-1)_{*}\mathcal{O}_{\bbD} \right\}_{s \in X^{\vee}}
\subset \op{ob}
\left(\sff{D}^{*}({{}_{\bB}\bbX^{\vee}},-1)\right)$. We will
argue that this is an orthogonal spanning class of
$\sff{D}^{*}( {{}_{\bB}\bbX^{\vee}},-1)$.

The orthogonality is obvious since for any two points $s \neq t \in
X^{\vee}$ the supports of the sheaves
$(\mathfrak{s},-1)_{*}\mathcal{O}_{\bbD}$ and
$(\mathfrak{t},-1)_{*}\mathcal{O}_{\bbD}$ are disjoint substacks in 
${}_{B}\bbX^{\vee}$. 

To show that these sheaves span the category, we need to check
that if $A$ is a complex of sheaves on ${}_{B}\bbX^{\vee}$ of pure
weight $(-1)$, with the property that $A$ is left (respectively right)
orthogonal to all $(\mathfrak{s},-1)_{*}\mathcal{O}_{\bbD}$, then $A =
0$ in $\sff{D}^{*}( {{}_{\bB}\bbX^{\vee}},-1)$.  Suppose first
$R\op{Hom}((\mathfrak{s},-1)_{*}\mathcal{O}_{\bbD}, A) = 0$ 
for all $s \in X^{\vee}$. Let $j : X^{\vee} \to {}_{B}\bbX^{\vee}$ be
the natural closed immersion. We have a distinguished triangle 
\begin{equation} \label{eq:triangleA}
\xymatrix@1{A \ar[r]^-{\hbar} & A \ar[r] & (j,-1)_{*}L(j,-1)^{*} A
  \ar[r] & A[1]} 
\end{equation}
and since the complex of vector spaces
$R\op{Hom}((\mathfrak{s},-1)_{*}\mathcal{O}_{\bbD}, A)$ is
exact, it follows that the complex
$R\op{Hom}((\mathfrak{s},-1)_{*}\mathcal{O}_{\bbD},
(j,-1)_{*}L(j,-1)^{*}A)$ is exact.

By adjunction we get that
\begin{equation} \label{eq-adjtoA}
R\op{Hom}_{X}(L(j,-1)^{*}(\mathfrak{s},-1)_{*}\mathcal{O}_{\bbD},
L(j,-1)^{*}A) = 0.
\end{equation}
 However $L(j,-1)^{*}(\mathfrak{s},-1)_{*}\mathcal{O}_{\bbD}$
can be computed explicitly:
\begin{equation} \label{eq:point}
L(j,-1)^{*}(\mathfrak{s},-1)_{*}\mathcal{O}_{\bbD} =
Lj^{*}\mathfrak{s}_{*} \mathcal{O}_{\bbD} = 
\mathcal{O}_{s}.
\end{equation}
The first  equality in \eqref{eq:point} follows tautologically from
the definition of the functors $L(j,-1)^{*}$ and
$(\mathfrak{s},-1)_{*}$ and the second follows immediately from the
base change identity $Lj^{*}\mathfrak{s}_{*}\mathcal{O}_{\bbD} =
Rs_{*}Li^{*}\mathcal{O}_{\bbD} = Rs_{*}\mathcal{O}_{0} =
\mathcal{O}_{s}$. Here $i : 0 \to \bbD$ denotes the inclusion of the
closed point and $s : 0 \to X^{\vee}$ is the map given by the point
$s$.

Now  \eqref{eq-adjtoA} and \eqref{eq:point} imply that 
$R\op{Hom}_{X}(\mathcal{O}_{s},
L(j,-1)^{*}A) = 0$ for all $s \in X^{\vee}$. Since the structure sheaves of
points form a  spanning class in the derived category of $X^{\vee}$ it
follows that $L(j,-1)^{*} A = 0$. Thus $(j,-1)_{*}L(j,-1)^{*} A = 0$ in 
$\sff{D}^{*}( {{}_{\bB}\bbX^{\vee}},-1)$. Now the exact
triangle \eqref{eq:triangleA} implies that multiplication by $\hbar$
is an isomorphism on all cohomology sheaves of the complex $A$. By
Nakayama's lemma this implies that the cohomology sheaves of $A$ are
all zero and so $A$ is quasi-isomorphic to the zero complex.

\

\bigskip

\noindent
{\bfseries Step 2.} Given a section $\mathfrak{s} : \bbD \to
\bbX^{\vee}$ of $\bbX^{\vee} \to \bbD$, we will
check that 
\begin{equation} \label{eq:Ps}
\phi_{\mycal{P}}((\mathfrak{s},-1)_{*}\mathcal{O}_{\bbD})
\cong \mycal{P}_{\mathfrak{s}},
\end{equation}
 where $\mycal{P}_{\mathfrak{s}}$ is
the line bundle  
\[
\mycal{P}_{\mathfrak{s}}:= L(\op{id} \times \mathfrak{s},1)^{*}\mycal{P} =
 (\op{id} \times \mathfrak{s},1)^{*}\mycal{P}
\]
on
$\bbX_{\bpi}$. 

Indeed, if we write $\bbu_{\bpi}  : \bbX_{\bpi} \to \bbD$ for the
structure morphism of $\bbX_{\bpi}$ we get a natural cartesian
diagram
\[
\xymatrix@C+2pc{
\bbX_{\bpi} \ar[d]_-{\bbu_{\bpi}} \ar[r]^-{\op{id}\times \mathbb{s}}
 & \bbX_{\bpi}\times_{\bbD}
    {}_{\bB}\bbX^{\vee} \ar[d]^-{\bbp_{2}} \\
\bbD \ar[r]_-{\mathbb{s}} &  {}_{\bB}\bbX^{\vee}
}
\] 
where $\mathbb{s} : \bbD \to {}_{\bB}\bbX^{\vee}$ denotes the map
corresponding to our preferred trivialization of the
$\mathcal{O}^{\times}$-gerbe
$\mathfrak{s}^{*}({}_{\bB}\bbX^{\vee})$ on $\bbD$. 

Base changing along this diagram now gives
$p_{2}^{-1}(\mathfrak{s},-1)_{*}\mathcal{O}_{\bbD} = R(\op{id}\times
\mathfrak{s},-1)_{*}\bbu_{\bpi}^{-1}\mathcal{O}_{\bbD}$. Thus 
\begin{equation} \label{eq:projectionP}
\begin{split}
\mycal{P}\otimes_{p_{2}^{-1}\mathcal{O}_{{}_{\bB}\bbX^{\vee}}}
p_{2}^{-1}(\mathfrak{s},-1)_{*}\mathcal{O}_{\bbD} & =
\mycal{P}\otimes_{p_{2}^{-1}\mathcal{O}_{{}_{\bB}\bbX^{\vee}}} R(\op{id}\times
\mathfrak{s},-1)_{*}\bbu_{\bpi}^{-1}\mathcal{O}_{\bbD} \\
& = R(\op{id}\times \mathfrak{s},0)_{*}(L(\op{id}\times
\mathfrak{s},1)^{*}\mycal{P}\otimes_{\bbu_{\bpi}^{-1}\mathcal{O}_{\bbD}}
\bbu_{\bpi}^{-1}\mathcal{O}_{\bbD}) \\
& = R(\op{id}\times \mathfrak{s},0)_{*}\mycal{P}_{\mathfrak{s}},
\end{split}
\end{equation}
where for the second equality we used the projection formula applied
to the map $\op{id}\times \mathbb{s}$.

We are now in a position to check \eqref{eq:Ps}. By definition, we have 
\[
\phi_{\mycal{P}}((\mathfrak{s},-1)_{*}\mathcal{O}_{\bbD}) = 
Rp_{1*}(\mycal{P} \otimes_{p_{2}^{-1} 
\mathcal{O}_{{}_{\bB}\bbX^{\vee}}} 
p_{2}^{-1}((\mathfrak{s},-1)_{*}\mathcal{O}_{\bbD}))
\]
and therefore by \eqref{eq:projectionP} 
\[
\begin{split}
\phi_{\mycal{P}}((\mathfrak{s},-1)_{*}\mathcal{O}_{\bbD}) & =
Rp_{1*}R(\op{id}\times \mathfrak{s},0)_{*}\mycal{P}_{\mathfrak{s}}  \\
& = R(p_{1}\circ (\op{id}\times \mathfrak{s},0))_{*}
\mycal{P}_{\mathfrak{s}} \\
& = R(\op{id})_{*} \mycal{P}_{\mathfrak{s}} \\
& = \mycal{P}_{\mathfrak{s}}.
\end{split}
\]
\

\begin{rem} This calculation only uses the fact that
  $\mathfrak{s} : \bbD \to \bbX^{\vee}$ is a section and is
  insensitive to whether this section is constant or not. In
  particular our proof shows that the identity \eqref{eq:Ps} is valid
  for all (not necessarily constant) sections $\mathfrak{s} : \bbD \to
  \bbX^{\vee}$. 
\end{rem}

\

\bigskip

\noindent
{\bfseries Step 3.} Finally we check that
\[
R\op{Hom}_{\sff{D}^{*}({}_{\bB}{\bbX^{\vee}})}(
(\mathfrak{s},-1)_{*}\mathcal{O}_{\bbD},  
(\mathfrak{t},-1)_{*}\mathcal{O}_{\bbD}) = 
R\op{Hom}_{\sff{D}^{*}(\bbX_{\bpi})}(
\mycal{P}_{\mathfrak{s}},\mycal{P}_{\mathfrak{t}})
\]
for all constant sections $\mathfrak{s}, \mathfrak{t} : \bbD \to
{}_{\bB}\bbX^{\vee}$. 

Since $\mycal{P}_{\mathfrak{s}}$ and $\mycal{P}_{\mathfrak{t}}$ are
translation invariant line bundles on $\bbX_{\bpi}$, they are
naturally $\Ah{\bpi}$ bimodules. In particular there is a well defined
inner hom  
\[
\mycal{H}om_{{\mycal{A}_{\bpi}}-\op{mod}}(\mycal{P}_{\mathfrak{s}}, 
\mycal{P}_{\mathfrak{t}}) =
\mycal{P}_{\mathfrak{s}}^{\vee}\otimes_{\Ah{\bpi}} \mycal{P}_{\mathfrak{t}}
\]
which is also a translation invariant line bundle on $\bbX_{\bpi}$. In
fact, by writing the factors of automorphy for
$\mycal{P}_{\mathfrak{s}}^{\vee}$ and $\mycal{P}_{\mathfrak{t}}$ it is
clear that $\mycal{P}_{\mathfrak{s}}^{\vee}\otimes_{\Ah{\bpi}}
\mycal{P}_{\mathfrak{t}} = \mycal{P}_{\mathfrak{t} - \mathfrak{s}}$.

Thus for the global homomorphisms in the derived category we get
\[
R\op{Hom}^{\bullet}(\mycal{P}_{\mathfrak{s}}, 
\mycal{P}_{\mathfrak{t}}) \cong  H^{\bullet}(X,\mycal{P}_{\mathfrak{t}
  - \mathfrak{s}}).
\]
Now by the computation in subsection~\ref{ss:cohomologyleft} we
conclude that $H^{\bullet}(X,\mycal{P}_{\mathfrak{t}
  - \mathfrak{s}}) = 0$ unless $\mycal{P}_{\mathfrak{t}
  - \mathfrak{s}}/\hbar \cong \mathcal{O}_{X}$, or equivalently
$\mathfrak{s}(0) = \mathfrak{t}(0)$. Since our sections $\mathfrak{s}$
and $\mathfrak{t}$ are constant, this can happen only when
$\mathfrak{s} = \mathfrak{t}$.

Finally, when $\mathfrak{s} = \mathfrak{t}$, we have
$\mycal{P}_{\mathfrak{t} - \mathfrak{s}} = \mycal{P}_{\mathfrak{0}} =
\Ah{\bpi}$, and so as a sheaf on $X$ we have $\mycal{P}_{\mathfrak{0}}
= \mathcal{O}_{X}[[\hbar]]$. In other words
$H^{\bullet}(X,\mycal{P}_{\mathfrak{0}}) =
H^{\bullet}(X,\mathcal{O})[[\hbar]]$.  This completes the proof of
  Step 3. \ \hfill $\Box$

\

\medskip 

Lemma~\ref{lem:diagonal} together with  Lemma~\ref{lem:phiPfaithful}
now yield a proof of the theorem 
\ \hfill $\Box$

\subsection{Remarks on classical supports} \label{ss:supports}

The reader may have noticed by now that there are some suggestive
similarities between the module theory on a deformation
quantization and the theory of D-modules. Here we point out a
particular aspect of this similarity that has to do with the supports
of quantum modules.

\begin{defn} \label{def:support} {\bf (i)} Suppose $\mathbb{M} = (M,\Ah{M})$
    is a deformation quantization of a complex manifold $M$. Let 
\[
\xymatrix@C-2pc{
M \ar[d] & \stackrel{i}{\subset} & \mathbb{M} \ar[d] \\
o & \in & \bbD
}
\] 
be the inclusion of the closed fiber and let 
$F$ be  a coherent sheaf on $\mathbb{M}$.  We define the {\em
    classical support of $F$} as the support of the complex
$Li^{*}F \in \sff{D}^{b}_{c}(M)$.

\

\smallskip

\noindent
{\bf (ii)} Suppose $N$ is a complex manifold and suppose $\mathfrak{N}
\to \mathbb{N} = N\times \bbD$ is an $\mathcal{O}^{\times}$-gerbe which
  is trivialized on the closed fiber $N\times \{ o \} \subset
  \mathbb{N}$. Let 
\[
\xymatrix@C-2pc{
N \ar[d] & \stackrel{i}{\subset} & \mathfrak{N} \ar[d] \\ 
o & \in & \bbD
}
\] 
be the inclusion corresponding to this trivialization and let $G$ be a
coherent sheaf of pure weight on $\mathfrak{N}$. We define the {\em
  classical support of $G$} as the support of the complex 
$Li^{*}G \in \sff{D}^{b}_{c}(N)$.
\end{defn}

Observe that if the sheaf $F$ on $\mathbb{M}$ is flat over $\bbD$, then the
classical supports of $F$ is just the support of the sheaf
$i^{*}F$. Similarly if $G$ is flat over $\mathfrak{N}$, then the
classical support of $G$ is the support of the sheaf $i^{*}G$.

\begin{claim} \label{claim:support} {\bf (a)} Suppose $\mathbb{M} =
  (M,\Ah{M})$ is a deformation quantization of a complex manifold $M$
and suppose $F$ is a coherent sheaf on $\mathbb{M}$ which is flat over
$\bbD$.  Let
$\bpi \in H^{0}(M,\wedge^{2}T_{M})$ be the holomorphic Poisson
structure associated with $\mathbb{M}$,  and let
$\left\{\bullet,\bullet \right\}_{\bpi} : \mathcal{O}_{M}\otimes_{\mathbb{C}}
  \mathcal{O}_{M} \to \mathcal{O}_{M}$ be the corresponding bracket on
  functions.   Then the classical support $S \subset M$  of $F$ is
  coisotropic with respect to $\bpi$, i.e. the ideal sheaf
  $\mathcal{I}_{S} \subset \mathcal{O}_{M}$ satisfies
  $\left\{\mathcal{I}_{S}, \mathcal{I}_{S} \right\}_{\bpi} \subset
  \mathcal{I}_{S}$. 

\

\smallskip
\noindent
{\bf (b)} Suppose $N$ is a complex manifold and suppose $\mathfrak{N}
\to \mathbb{N} = N\times \bbD$ is an $\mathcal{O}^{\times}$-gerbe which
  is trivialized on the closed fiber. Let $\bB \in
  H^{2}(N,\mathcal{O})$ be the infinitesimal class of the gerbe
  $\mathfrak{N}$ and let $G$ be a pure weight coherent sheaf on
  $\mathfrak{N}$ which is flat over $\bbD$. Then the classical support
  $T \subset N$ of $G$ is $\bB$-isotropic, i.e. $\bB_{|T} = 0 \in
  H^{2}(T,\mathcal{O}_{T})$. 
\end{claim}

\begin{rem} \label{rem:isotropic} There is a puzzling asymmetry in the
notions of coisotropic and isotropic defined above. Whereas the
property of being coisotropic is geometric, the property of being
isotropic appears to be only homological.  To justify why our notion
of isotropic is meaningful note that if $\bB \in
H^{2}(N,\mathcal{O}_{N})$, and if we choose a Hermitian metric on $N$,
then $\bB$ can be represented by a $\bar{\partial}$-harmonic
$(0,2)$-form $\boldsymbol{\beta}$. When the metric is K\"{a}hler and
$N$ is compact, the form $\boldsymbol{\beta}$ is $d$-harmonic and
hence closed. In other words, $\bB$ is represented by a presymplectic
form $\boldsymbol{\beta}$ on the $C^{\infty}$-manifold $N$. If now $T
\subset N$ is a complex submanifold for which $\bB_{|T} = 0 \in
H^{2}(N,\mathcal{O}_{N})$, it follows that the restriction
$\boldsymbol{\beta}_{|T} \in \Gamma(T,A^{0,2}_{T})$ is the zero
form, i.e. $T$ is isotropic in the sense of presymplectic geometry. 

Note also, that when $N = X^{\vee}$ was a complex torus and
$\mathfrak{N} = {}_{\bB}\bbX^{\vee}$, then the class $\bB$ had
a canonical harmonic representative, since $\bB$ was given by an element
in $\wedge^{2}V$. 
\end{rem}

\

\smallskip

\noindent
{\bf Proof of the Claim.} Part {\bf (a)} is a deformation quantization
analogue of the corresponding result for D-modules. It is in fact a
special case of a general version of Gabber's theorem
\cite[Appendix~III,Theorem~3.7]{Bjork}. The result
\cite[Appendix~III,Theorem~3.7]{Bjork} implies the usual Gabber's
theorem on D-modules when applied to $\mathcal{R}/t^{2}$, where
$\mathcal{R}$ is the Rees algebra of the sheaf of differential
operators. It also implies part {\bf (a)} of our claim when applied to
$\Ah{M}/\hbar^{2}$. The only thing we need to check is that hypothesis
$\ker\left[F/\hbar^{2} \stackrel{\hbar\cdot }{\to} F/\hbar^{2}\right] =
\op{im}\left[F/\hbar^{2} \stackrel{\hbar\cdot }{\to}  F/\hbar^{2}\right]$
which is immediate from the fact that $F$ is flat over $\bbD$.

\

\smallskip

\noindent
Part {\bf (b)} follows from \cite[Proposition~6.1]{Toda} applied to
$\mathfrak{N}/\hbar^{2}$. Indeed, according to
\cite[Proposition~6.1]{Toda}, for any sheaf $\Phi$ on $N$, the
existence of a flat extension of 
$\Phi$ to $\mathfrak{N}/\hbar^{2}$ is equivalent to the vanishing of
the product of the exponential $\exp(a(\Phi)) \in \oplus_{k}
\op{Ext}^{k}(G,G\otimes \Omega_{N}^{k})$ of the Atiyah class $a(\Phi)$ of
$\Phi$, and the $\bB \in H^{2}(N,\mathcal{O}_{N}) \subset
HT^{2}(N) = \oplus_{p+q = 2} H^{p}(N,\wedge^{q}T_{N})$. Here the
product is defined as the image under the
natural map 
\[
\xymatrix@1@C-1.2pc{
\left(\oplus_{k}
\op{Ext}^{k}(G,G\otimes \Omega_{N}^{k})\right)\otimes
\left(\oplus_{p,q} H^{p}(N,\wedge^{q}T_{N})\right) \ar[r]^-{\cup} & 
\oplus_{a,b}
\op{Ext}^{a}(G,G\otimes \Omega_{N}^{b}) \ar[r] & \oplus_{a}
\op{Ext}^{a}(G,G).
}
\]
In particular $\exp(a(\Phi))\cdot \bB$ can be zero only if $\bB$
restricts to zero on the support of $\Phi$. Applying this to $\Phi :=
i^{*}G$ and taking into account that $G$ is flat, we get the statement
{\bf (b)}. The claim is proven. \ \hfill $\Box$ 

\

\medskip

\noindent
The claim together with the classical Fourier-Mukai duality impose
non-trivial conditions on the support of a sheaf on $M$, that are
necessary for quantizing it. For instance, part {\bf (b)} of the claim
immediately implies that an ample line bundle on an abelian variety
can not be quantized since its Fourier-Mukai transform is a vector
bundle on the dual abelian variety.  The supports of modules over deformation 
quantizations was recently investigated in \cite{Tsygan2}.

\section{The gerbe $\mycal{P}ic^{0}(\bbX_{\bpi}/\bbD)$}  \label{s:Pic}
With the equivalence of categories in place we are now ready to
identify ${}_{\bB}\bbX^{\vee}$ geometrically as the relative Picard
stack $\mycal{P}ic^{0}(\bbX_{\bpi}/\bbD)$ of degree zero line bundles
on $\bbX_{\bpi} \to \bbD$.

\ 

\bigskip

\noindent
Before we define $\mycal{P}ic^{0}(\bbX_{\bpi}/\bbD)$ let us recall the
classical notion of a Picard variety.  In the classical situation, the
Picard variety $\op{Pic}^0(Z)$ of a smooth space $Z$ is defined as the
moduli space of isomorphism classes of degree zero line bundles. More
precisely, consider the site $\mycal{S}_{\op{an}}$ of analytic spaces
with the analytic topology and let $\mycal{P}ic^{0}(Z)$ denote the
moduli stack on $\mycal{S}_{\op{an}}$ associated to the prestack:
\begin{equation} \label{eq:Pic0_classical}
Y \to \left\{ \begin{minipage}[c]{2.5in} the groupoid
whose objects are
holomorphic line bundles on $Z
 \times Y$ of degree zero relative to $Y$, and morphisms are isomorphisms.
\end{minipage} \right\}.
\end{equation}
The associated sheaf of sets $\pi_{0}(\mycal{P}ic^{0}(Z))$ 
associated to $\mycal{P}ic^{0}(Z)$ is representable by a space
$\op{Pic}^{0}(Z)$.  The stack $\mycal{P}ic^{0}(Z)$ is an
$\mathcal{O}^{\times}$-gerbe over $\op{Pic}^{0}(Z)$ which encodes the
fact that the automorphisms of a line bundle are given by
multiplication by an invertible holomorphic function. Explicitly, for
any analytic space $Y$ one introduces an equivalence relation
$\sim_{Y}$ on the collection of all line bundles on $Z\times Y$. Two
line bundles $L, M \to Z\times Y$ are considered to be equivalent if
$L$ is isomorphic to $M\otimes p_{Y}^{*}A$ for some line bundle $A$ on
$Y$. The variety $\op{Pic}^0(Z)$ represents the functor 
\[
\pi_{0}(\mycal{P}ic^{0}(Z)) : \mycal{S}_{\op{an}} \to \sff{Sets}
\]
given by
\[
Y  \to \left\{ \begin{minipage}[c]{2.6in} the set of all
 $\sim_{Y}$-equivalence classes of holomorphic line bundles on $Z
 \times Y$ of degree zero relative to $Y$.
\end{minipage} \right\}.
\]
Moreover the moduli problem \eqref{eq:Pic0_classical} can be
rigidified in a simple way (by considering line bundles on $Z\times Y$
equipped with a trivialization $\{o\}\times Y$ for some fixed point $o
\in X$.) which trivializes this gerbe
$\mycal{P}ic^{0}(Z)$. Interestingly, in the non-commutative case we
can not resort to a rigidification trick since the non-commutative
space $\bbX_{\bpi}$ need not have any points. In fact, as we will see
below, $\mycal{P}ic^{0}(\bbX_{\bpi}/\bbD)$ is a gerbe on $\bbX^{\vee}$
which is no longer trivializable.

Do define line bundles of degree zero on $\bbX_{\bpi}$ we look at the
translation action of a torus $X$  on itself. Since the holomorphic
Poisson structures are constant, this action lifts to the sheaf
$\Ah{X,\bpi}$.  In particular, $X$ acts on the non-commutative space
$\bbX_{\bpi}$. In the classical situation, the degree zero line
bundles can be characterized as those that are translation
invariant. We use this as our definition of degree $0$ in the
non-commutative case.

\begin{defn}\label{defn:degzero}
A line bundle (= a locally free rank one left $\Ah{X,\bpi}$-module) on
$\bbX_{\bpi}$ is said to be {\em of degree zero} when it is
translation invariant.
\end{defn}

\

\noindent
It turns out that a line bundle $\mycal{L}$ on $\bbX_{\bpi}$ is of
degree zero if and only if its classical part
$\mycal{L}/\hbar\mycal{L}$ has zero first Chern class (see
Lemma~\ref{lem:quantizable}). It is a bit of an accident that such an
$\mycal{L}$ is also a bimodule.  We will not use these bimodule
structures.  As discussed before in section
\ref{prop:invertiblefamilies}, if one wants line bundles to vary in
families, the bimodule structures on the individual line bundles
cannot be chosen in a consistent way.

By definition $\mycal{P}ic^{0}(\bbX_{\bpi}/\bbD)$ is the moduli stack
of line bundles of relative degree zero on $\bbX_{\bpi} \to \bbD$. To
spell this out we will need the analytic site $(\mycal{FS}/\bbD)_{an}$
of formal analytic spaces over $\bbD$. Formal analytic spaces are the
analytic counterpart of Knutson's formal algebraic spaces. The theory
of formal analytic spaces is parallel to \cite[Chapter]{knutson} but
with commutative rings replaced with Stein algebras. A convenient way
to look at the formal analytic spaces over $\bbD$ is as commutative
deformation quantizations of analytic spaces. A morphism ${\mathbb f}
= (f,f^{\sharp}) : \bbX \to \bbY$ in the category $\mycal{FS}/\bbD$ is
an {\em analytic open immersion} if the map $f = {\mathbb f}/\hbar : X
\to Y$ is an open analytic map, and the map $f^{\sharp} :
f^{-1}{\mathcal O}_{Y} \to {\mathcal O}_{X}$ is an isomorphism. The
analytic topology on the category of formal analytic spaces is the
Grothendieck topology associated (in the sense of say \cite{knutson})
with the subcategory of {\em analytic open immersions}.

Given a formal analytic space $\bbY$ over $\bbD$, with closed fiber
$Y$ we form the ringed
space
\[\bbX_{\bpi} \times_{\bbD} \bbY = (X\times Y,
p_{X}^{-1}\Ah{X,\bpi}\widehat{\otimes}_{{\mathbb C}[[\hbar]]}
p_{Y}^{-1}\Ah{Y}) 
\]
\noindent
where $\widehat{\otimes}$ denotes the completed tensor product of
sheaves of nuclear Frechet algebras.

Given a formal analytic space $\bbY$ over $\bbD$, we say that a line
bundle $\mycal{L}$ on $\bbX_{\bpi} \times_{\bbD} \bbY$ is of degree
zero relative to $\bbY$ if for any section $\sigma: \bbD \to \bbY$,
the pullback $(1 \times \sigma)^{*}(\mycal{L})$ is of degree zero as
a line bundle on $\bbX_{\bpi}$.

\begin{defn} \label{defn:Pic0}
The moduli stack $\mycal{P}ic^{0}(\bbX_{\bpi}/\bbD)$ is the stack on
$(\mycal{FS}/\bbD)_{an}$ associated with the prestack:
\[
(\bbY \to \bbD) \to \left\{ \begin{minipage}[c]{2.5in} the groupoid
whose objects are line bundles over $\bbX_{\bpi} \times_{\bbD} \bbY$
of degree zero relative to $\bbY$, and morphisms are isomorphisms
\end{minipage} \right\}.
\]
\end{defn}

\

\bigskip

\noindent
By definition the Poincar\'{e} sheaf $\mycal{P}$ on
$\bbX_{\bpi}\times_{\bbD} ({}_{\bB}\bbX^{\vee})$ is 
line bundle of relative degree zero along $\bbX_{\bpi}$. Indeed if
we pull-back $\mycal{P}$ by a translation by a point in $X$, then we
get a line bundle isomorphic to $\mycal{P}$. This is easily seen in
terms of factors of automorphy. If $w
\in V$ is any point, then translation of the factor of automorphy by
$w$ results into a cohomologous factor of automorphy. The two factors
are related by the coboundary of the element
$g \in  C^{0}(\Lambda\times
\Gamma, \mycal{A}_{V\times \overline{V}^{\vee}}^{\times}(V\times
\overline{V}^{\vee}))  = \mycal{A}_{V\times
  \overline{V}^{\vee}}^{\times}(V\times 
\overline{V}^{\vee})$, where $g(l,v) := \exp(\pi\overline{\langle l,
  w\rangle})$.  

In particular the Poincar\'{e} sheaf $\mycal{P}$ gives rise to a
natural morphism of stacks 
\begin{equation}
\mathfrak{c} : {}_{\bB}\bbX^{\vee} \to \mycal{P}ic^{0}(\bbX_{\bpi}/\bbD), 
\end{equation}
which sends $\mathbb{f} : \mathbb{Y} \to {}_{\bB}\bbX^{\vee}$ to the
line bundle $\mathbb{f}^{*} \mycal{P}$ on $\bbX_{\bpi}\times_{\bbD}
\mathbb{Y}$. 

Our goal is to show that $\mathfrak{c}$ is an isomorphism. The first
step is to analyze the
relationship between the ${}_{\bB}\bbX^{\vee}$ and
$\mycal{P}ic^{0}(\bbX_{\bpi}/\bbD)$ on the level of $\bbD$-points.

\begin{prop}\label{prop:param} The map $\mathfrak{c}$ induces
an equivalence between the groupoid ${}_{\bB}\bbX^{\vee}(\bbD)$ of 
all sections of ${}_{\bB}\bbX^{\vee} \to \bbD$ and the groupoid 
$\mycal{P}ic^{0}(\bbX_{\bpi}/\bbD)(\bbD)$ of  global translation
invariant line  
bundles on $\bbX_{\bpi}$.
\end{prop}
{\bf Proof.} The set of isomorphism classes of the groupoid
${}_{\bB}\bbX^{\vee}(\bbD)$ is the set $\bbX^{\vee}(\bbD)$ of
$\bbD$-points of the formal space $\bbX^{\vee} \to \bbD$. The natural
map from ${}_{\bB}\bbX^{\vee}(\bbD)$ to the discrete groupoid
$\bbX^{\vee}(\bbD)$ is actually split. This is a special feature of
the groupoid of $\bbD$-points and no such splitting exist for general
test spaces $\bbY \in \mycal{F}\mycal{S}/\bbD$. For our purposes it
will be important to exhibit a distinguished splitting of
${}_{\bB}\bbX^{\vee}(\bbD) \to \bbX^{\vee}(\bbD)$. In fact we have
already done that in the proof of {\bf Step 1} of
Lemma~\ref{lem:phiPfaithful}. Specifically, to lift a section $\mathfrak{s}
: \bbD \to \bbX^{\vee}$ to a section $\mathbb{s} : \bbD \to
{}_{\bB}\bbX^{\vee}(\bbD)$, we have to produce a natural
trivialization of the pullback gerbe
$\mathfrak{s}^{*}({}_{\bB}\bbX^{\vee})$. For this we can use the
global quotient presentation of ${}_{\bB}\bbX^{\vee}(\bbD)$. Write $s
\in X^{\vee}$ for the image $\mathfrak{s}(0)$ of the closed point $0
\in \bbD$. Then in the notation of the proof of {\bf Step 1} of
Lemma~\ref{lem:phiPfaithful} we have a quotient presentation $[(F_{s}\times
\bbD)/\bGamma]$ of $\mathfrak{s}^{*}({}_{\bB}\bbX^{\vee}(\bbD))$. Now
a trivialization of the gerbe $[(F_{s}\times \bbD)/\bGamma]$ is simply
a $\Gamma$-equivariant structure on $\mathcal{O}_{F_{s}}[[\hbar]]$ in
which $z \in \mathbb{C}[[\hbar]] \subset \Gamma$ acts as
multiplication by $z$ (viewed as a section of
$\mathcal{O}_{F_{s}}[[\hbar]]$). The cocycle $c :
\wedge^{2}\Lambda^{\vee} \to \Ch^{\times}$ gives rise to such an
equivariant structure: a group element $(\xi,z) \in \Gamma$ acts on $f
= \{ f_{w} \}_{w \in F_{s}}$ via the formula $f \mapsto \{
z\tilde{c}(w,\xi)f_{w+\xi}\}_{w \in F_{s}}$.  Thus we get a natural
section $\mathfrak{a} : \bbX^{\vee}(\bbD) \to
{}_{\bB}\bbX^{\vee}(\bbD)$ and hence an isomorphism between the groupoid
${}_{\bB}\bbX^{\vee}(\bbD)$ and the groupoid $\bbX^{\vee}(\bbD)\times
B(\Ch^{\times})$. 

On the non-commutative side, the set
$\pi_{0}(\mycal{P}ic^{0}(\bbX_{\bpi}/\bbD)(\bbD))$ of isomorphism
classes of global degree zero quantum line bundles on $\bbX_{\bpi}\to
\bbD$, coincides with the set $\op{Pic}^{0}(\bbX/\bbD)(\bbD)$ of
isomorphism classes of global degree zero line bundles on $\bbX \to
\bbD$. This fact is not obvious but can be established as
follows. First of all, as explained in Lemma~\ref{lem:quantizable},
every element in the set $X^{\vee}\times
(\overline{V}^{\vee})^{\mathbb{Z}_{>0}}$ gives rise to a quantum line
bundle of degree zero, and this procedure induces a bijection between
$X^{\vee}\times (\overline{V}^{\vee})^{\mathbb{Z}_{>0}}$ and the set of
isomorphism classes of quantum line bundles of degree zero. On the
other hand, the set $\op{Pic}^{0}(\bbX/\bbD)(\bbD)$ is simply the set
of formal arcs in $X^{\vee} = \op{Pic}^{0}(X)$ and so can be described
\cite{vojta} explicitly in terms of Hasse-Schmidt higher derivations
\cite{hasse-schmidt,matsumura}.

Recall that a Hasse-Schmidt higher $\mathbb{C}$-derivation from 
$\mathcal{O}_{X^{\vee}}$ to $\mathbb{C}$ is a pair $(s,D)$, where $s
\in X^{\vee}$ is a point and $D$ is  
(bounded or unbounded) sequence 
$D = (D_{1}, D_{2}, D_{3}, \dots)$ of $\mathbb{C}$-linear maps  
$D_{i}: \mathcal{O}_{X^{\vee},s} \to \mathbb{C}$ satisfying
\[ 
D_{k}(fg) = \sum_{i = 0}^{k} D_{i}(f) D_{k-i}(g)
\]
for any two germs $f$ and $g$ of $\mathcal{O}_{X^{\vee},s}$. Here
$\mathcal{O}_{X^{\vee},s}$ denotes the local ring at $s \in X^{ \vee}$
and $D_{0}: \mathcal{O}_{X^{\vee},s} \to \mathbb{C}$ is the evaluation
map.

We denote the infinite order higher order derivations by
$\op{Der}^{\infty}_{\mathbb{C}}(\mathcal{O}_{X^{\vee}},
\mathbb{C})$. It is not hard to identify these derivations with the
set of formal arcs in $X^{\vee}$. Indeed, to specify a formal arc
$\mathfrak{s} : \bbD \to X^{\vee}$, we need to specify a point $s \in
X^{\vee}$ and a $\mathbb{C}$-algebra homomorphism
$\mathfrak{s}^{\sharp} : \mathcal{O}_{X^{\vee},s} \to \Ch$. In these
terms, the identification
\begin{equation} \label{eq:arcs-are-derivations} 
\op{Der}^{\infty}_{\mathbb{C}}
(\mathcal{O}_{X^{\vee}},\mathbb{C}) \cong 
\op{Hom}_{\bbD}(\bbD, \bbX^{\vee}) = \op{Hom}(\bbD,X^{\vee})
\end{equation}
is given by  
\[
(s,D) \mapsto \left(s, \mathfrak{s}^{\sharp} := \sum_{k = 0}^{\infty}
\hbar^{k} D_{k}\right).
\]
On the other hand we
have a bijection
\begin{equation} \label{eq:bundles-are-derivations}
X^{\vee}
\times (\overline{V}^{\vee})^{\mathbb{Z}_{>0}} \to
\op{Der}^{\infty}_{\mathbb{C}} 
(\mathcal{O}_{X^{\vee}},\mathbb{C})
\end{equation}
given by 
\[
(s, l_{1}, l_{2}, l_{3}, \dots) \mapsto 
(D_{0}=\op{ev}_{s}, D_{1}, D_{2}, D_{3}, \dots).
\]
Where 
\[
D_{0} \circ \exp\left(\pi\left( \hbar l_{1} + 
\hbar^{2} l_{2}+ \hbar^{3} l_{3} +  \cdots\right)\right)
= D_{0} + \hbar D_{1} + \hbar^{2} D_{2} + \cdots
\]
and the $l_{i} \in \overline{V}^{\vee} $ are thought of as 
translation invariant vector fields on $X^{\vee}$. The exponential in
the left hand side of this formula is 
defined by the usual power series using the composition 
of differential operators.  This is an extension of the standard fact 
that any tangent germ $l$, defined at a smooth point 
$m$ in some complex analytic space $M$, can be exponentiated to a
formal arc $\mathfrak{e} : \bbD \to M$ in $M$: 
\[
\mathfrak{e} := \left(D_{0}, D_{0} \circ l,  \frac{1}{2}D_{0} \circ l^{2}, 
\frac{1}{3!} D_{0} \circ l^{3}, \frac{1}{4!} D_{0} \circ l^{4},\dots
\right) \in  
\op{Der}^{\infty}_{\mathbb{C}}(\mathcal{O}_{M},\mathbb{C}),
\]
where $D_{0}$ denotes the evaluation map at the point $m$ and $l^{k}$
denotes the $k$-th iterated Lie derivative.

Combining the bijections \eqref{eq:arcs-are-derivations} and
\eqref{eq:bundles-are-derivations} with the fact that $X^{\vee}\times
(\overline{V}^{\vee})^{\mathbb{Z}_{>0}}$ parameterizes isomorphism
classes of degree zero quantum line bundles (see
Lemma~\ref{lem:ql}), we obtain a natural
identification
\[
\pi_{0}(\mycal{P}ic^{0}(\bbX_{\bpi}/\bbD)(\bbD))\cong
\bbX^{\vee}(\bbD). 
\]
Thus we get a map $\mycal{P}ic^{0}(\bbX_{\bpi}/\bbD)(\bbD) \to
\bbX^{\vee}(\bbD)$ which similarly to the map
${}_{\bB}\bbX^{\vee}(\bbD) \to \bbX^{\vee}(\bbD)$ admits a preferred
splitting $\mathfrak{b} : \bbX^{\vee}(\bbD) \to
\mycal{P}ic^{0}(\bbX_{\bpi}/\bbD)(\bbD)$, defined to be the
composition of the correspondence of Lemma~\ref{lem:quantizable} with
the identification \eqref{eq:bundles-are-derivations}. In particular we
get a natural equivalence of groupoids between
$\mycal{P}ic^{0}(\bbX_{\bpi}/\bbD)(\bbD)$ and the groupoid
$\bbX^{\vee}(\bbD)\times B(\Ch^{\times})$.

Therefore we get an explicit equivalence 
\begin{equation} \label{eq:spliteq}
{}_{\bB}\bbX^{\vee}(\bbD) \cong \bbX^{\vee}(\bbD)\times 
B(\Ch^{\times}) \cong \mycal{P}ic^{0}(\bbX_{\bpi}/\bbD)(\bbD).
\end{equation}

Finally, to show that the equivalence \eqref{eq:spliteq} is given by
the map $\mathfrak{c}$ we need to chase through the sequence of
bijections defining \eqref{eq:spliteq}. This is tedious but
straightforward. The key observation here is that $\mathfrak{c} :
{}_{\bB}\bbX^{\vee}(\bbD) \to \mycal{P}ic^{0}(\bbX_{\bpi}/\bbD)(\bbD)$
maps the discrete subgropoid $\bbX^{\vee}(\bbD) \subset
{}_{\bB}\bbX^{\vee}(\bbD)$ identically to the discrete subgroupoid
$\bbX^{\vee}(\bbD) \subset
\mycal{P}ic^{0}(\bbX_{\bpi}/\bbD)(\bbD)$. This is sufficient to
conclude that \eqref{eq:spliteq} is given by $\mathfrak{c}$ since
$\mathfrak{c}$ is a morphism of groupoids. The statement that
$\mathfrak{c}$ induces the identity on $\bbX^{\vee}(\bbD)$ is easy to
check directly. Specifically we need to check that the two sections
$\mathfrak{c}\circ \mathfrak{a}$ and $\mathfrak{b}$ of
$\mycal{P}ic^{0}(\bbX_{\bpi}/\bbD)(\bbD) \to \bbX^{\vee}(\bbD)$   are
isomorphic, i.e. we need to exhibit an isomorphism of functors
$\boldsymbol{\iota}$:
\[
\xymatrix@M+0.5pc{
\bbX^{\vee}(\bbD) \rrtwocell<10>^{\mathfrak{c}\circ
  \mathfrak{a}}_{\mathfrak{b}}{\boldsymbol{\iota}}  & &
\mycal{P}ic^{0}(\bbX_{\bpi}/\bbD)(\bbD).
}
\]
For this we will need the explicit form of the functors
$\mathfrak{c}\circ \mathfrak{a}$ and $\mathfrak{b}$.
By definition, the splitting $\mathfrak{b}$   the
composition of the correspondence of Lemma~\ref{lem:quantizable} with
the identification \eqref{eq:bundles-are-derivations}. Given an arc
$\mathfrak{s}$ in $X^{\vee}$, we can describe the degree zero line
bundle $\mathfrak{b}(\mathfrak{s})$ on $\bbX_{\bpi}$ explicitly by a
factor of automorphy.  Use, \eqref{eq:bundles-are-derivations} to
write $\mathfrak{s} = (s,\mathfrak{s}^{\sharp})$,
$\mathfrak{s}^{\sharp} = \op{ev}_{s}\circ \exp(\sum_{i=1}^{\infty}
\hbar^{i}l_{i}))$ for some collection of $l_{i} \in
\overline{V}^{\vee}$. Now according to Lemma~\ref{lem:quantizable} 
the factor of automorphy for
$\mathfrak{b}(\mathfrak{s})$ is the map 
\begin{equation} \label{eq:b(s)}
\xymatrix@R-2pc{
\Lambda \ar[r] & \Gamma(V,\Ah{V,\bpi}^{\times}) \\
\lambda \ar@{|->}[r] & \chi_{s}(\lambda)\exp\left( \pi \sum_{i=1}^{\infty}
\hbar^{i} \langle l_{i}, \lambda \rangle \right),
}
\end{equation} 
where $\chi_{s} : \Lambda \to U(1) \subset \mathbb{C}^{\times}$ is the 
unitary character corresponding to $s \in X^{\vee}$. To compare this
factor of automorphy with the one corresponding to $\mathfrak{c}\circ
\mathfrak{a}(\mathfrak{s})$, it is convenient to realize $\bbX_{\bpi}$
not as the quotient $\bbV_{\bpi}/\bLambda$ but rather as the quotient 
$(\bbV_{\bpi}\times F_{s})/(\bLambda\times \Lambda^{\vee})$. Using the
the identification $\bbX_{\bpi} = (\bbV_{\bpi}\times
F_{s})/(\bLambda\times \Lambda^{\vee})$ we can now rewrite the factor
of automorphy \eqref{eq:b(s)} as a factor of automorphy for the group
$\Lambda\times \Lambda^{\vee}$ with values in $\Gamma\left(V, \prod_{w \in
  F_{s}} \Ah{V,\bpi}^{\times}\right)$.  In these terms
\eqref{eq:b(s)} becomes 
\begin{equation} \label{eq:b(s)lifted}
\xymatrix@R-2pc{
\Lambda\times \Lambda^{\vee} \ar[r] & \prod_{w \in F_{s}} \Gamma(V, 
\Ah{V,\bpi}^{\times}) \\ 
(\lambda,\xi) \ar@{|->}[r] & \left\{
e^{\left(2\pi\sqrt{-1}\op{Im}\langle w, \lambda \rangle \right)}e^{\left( \pi
\sum_{i=1}^{\infty} 
\hbar^{i} \langle l_{i}, \lambda \rangle \right)} \right\}_{w \in F_{s}}.
}
\end{equation}

The factor of automorphy defining $\mathfrak{c}\circ
\mathfrak{a}(\mathfrak{s})$ is also easy to describe.  Since the map
$\mathfrak{c}$ is given by the Poincar\'{e} sheaf $\mycal{P}$, the
degree zero line bundle $\mathfrak{c}\circ \mathfrak{a}(\mathfrak{s})$
on $\bbX_{\bpi}$ can be described as first restricting the sheaf
$\mycal{P}$ to the product
$\bbX_{\bpi}\times_{\bbD}(\mathfrak{s}^{*}({}_{\bB}\bbX^{\vee}))$ and
then pulling back this restriction by our preferred trivialization of
the gerbe $\mathfrak{s}^{*}({}_{\bB}\bbX^{\vee})$. Equivalently, we
can tensor the restriction of $\mycal{P}$ to
$\bbX_{\bpi}\times_{\bbD} (\mathfrak{s}^{*}({}_{\bB}\bbX^{\vee}))$ by
the structure sheaf on $\bbD$ viewed as a weight $(-1)$ line bundle on
$\mathfrak{s}^{*}({}_{\bB}\bbX^{\vee})$ (via the preferred
trivialization). Writing
$\bbX_{\bpi}\times_{\bbD}\mathfrak{s}^{*}({}_{\bB}\bbX^{\vee})$  as
the quotient $[(\bbV\times F_{s})/(\bLambda\times_{\bbD} \bGamma)]$ 
we can now describe the restriction of $\mycal{P}$ by the
$\Lambda\times \Gamma$ factor of automorphy 
\[
\mathfrak{s}^{\sharp}\phi : \Lambda\times \Gamma \to \prod_{w \in
  F_{s}} \Gamma(V,\Ah{V,\bpi}^{\times}).
\]
Here as usual $\phi$ is the factor of automorphy for $\mycal{P}$ (see
\eqref{PSFA}) and the terms in $\mathfrak{s}^{\sharp}$ act on
the functions $\phi(\lambda,(\xi,z))$ as
  iterated Lie derivatives. Similarly, we can use this quotient
  presentation, to write the factor of automorphy for the structure
  sheaf on $\bbD$ viewed as a weight $(-1)$ line bundle on 
$\mathfrak{s}^{*}({}_{\bB}\bbX^{\vee})$. As explained in the proof of
  {\bf Step 1} of Lemma~\ref{lem:phiPfaithful} this twisted line
  bundle corresponds to the $\Gamma$ equivariant structure on
  $\prod_{w \in F_{s}}\mathcal{O}_{\bbD}$ given by the representation
  $\rho$. In particular, the pullback of the twisted line bundle to
  the product
  $\bbX_{\bpi}\times_{\bbD}\mathfrak{s}^{*}({}_{\bB}\bbX^{\vee})$ is
  given by the factor of automorphy obtained by applying $\rho$ to the
  constant section $1$, i.e. by 
\begin{equation} \label{eq:twisted}
\xymatrix@R-2pc{ 
\Lambda\times \Gamma \ar[r] &  \prod_{w \in
  F_{s}} \Gamma(V,\Ah{V,\bpi}^{\times}), \\
(\lambda,(\xi,z)) \ar[r] & \rho_{(\xi,z)}(1) = \{ z^{-1}\tilde{c}(w,\xi)
\}_{w \in F_{s}}
}
\end{equation}
Multiplying $\mathfrak{s}^{\sharp}\phi$ and
\eqref{eq:twisted} gives rise to a factor of automorphy on $\Lambda\times
\Gamma$, which is a pullback of a factor of automorphy on 
$\Lambda\times \Lambda^{\vee}$ which describes the non-commutative
line bundle $\mathfrak{c}\circ\mathfrak{a}(\mathfrak{s})$. The later
factor is easily 
computed. The formula \eqref{PSFA} describing $\phi$, together with
the definition of $\mathfrak{s}^{\sharp}$ give
\[
(\mathfrak{s}^{\sharp}\phi)(\lambda,(\xi,z))_{w}(v) =
\exp(\pi\langle\xi+w, \lambda 
\rangle) \exp(\pi\overline{\langle\xi, v
\rangle} )\exp\left( \pi
\sum_{i=1}^{\infty} 
\hbar^{i} \langle l_{i}, \lambda \rangle \right).
\]
Hence the factor of automorphy of
$\mathfrak{c}\circ\mathfrak{a}(\mathfrak{s})$is given by the formula:
\begin{equation} \label{eq:ca(s)}
\xymatrix@R-2pc{
\Lambda\times \Lambda^{\vee} \ar[r] & \prod_{w \in
  F_{s}} \Gamma(V,\Ah{V,\bpi}^{\times}), \\
(\lambda,\xi) \ar[r] & \left\{ e^{(\pi\langle\xi+w, \lambda
\rangle)}\tilde{c}(w,\xi) e^{(\pi\overline{\langle\xi, \bullet
\rangle})}e^{\left( \pi
\sum_{i=1}^{\infty} 
\hbar^{i} \langle l_{i}, \lambda \rangle \right)} \right\}_{w \in F_{s}}.
}
\end{equation} 
\

\noindent
We are now ready to describe the isomorphism of functors
$\boldsymbol{\iota}$. Specifying $\boldsymbol{\iota}$ is equivalent to
specifying isomorphisms $\boldsymbol{\iota}_{\mathfrak{s}} :
\mathfrak{c}\circ\mathfrak{a}(\mathfrak{s}) \to
\mathfrak{b}(\mathfrak{s})$ of non-commutative line bundles. In terms
of the factors of automorphy \eqref{eq:b(s)lifted} and
\eqref{eq:ca(s)} the isomorphism $\boldsymbol{\iota}_{\mathfrak{s}}$
can be viewed as a group cochain $\boldsymbol{\iota}_{\mathfrak{s}} \in
C^{0}(\Lambda\times \Lambda^{\vee}, \prod_{w \in
  F_{s}} \Gamma(V,\Ah{V,\bpi}^{\times})) = \prod_{w \in
  F_{s}} \Gamma(V,\Ah{V,\bpi}^{\times})$ which makes the factor of
automorphy \eqref{eq:b(s)lifted} cohomologous to
\eqref{eq:ca(s)}. A straightforward computations shows that 
\[
\boldsymbol{\iota}_{\mathfrak{s}} := \left\{ e^{-\pi\overline{\langle
  w,\bullet \rangle }}  \right\}_{w \in
  F_{s}}
\]
does the job. This completes the proof of the proposition.
\ \hfill $\Box$

\

\bigskip

\noindent
Proposition~\ref{prop:param} together with Theorem~\ref{thm:main}
readily imply the following

\begin{thm} \label{thm:picard} The morphism
 $\mathfrak{c} : {}_{\bB}\bbX^{\vee} \to
\mycal{P}ic^{0}(\bbX_{\bpi}/\bbD)$ is an isomorphism of analytic
stacks on $(\mycal{F}\mycal{S}/\bbD)_{\op{an}})$.  
\end{thm} 
{\bf Proof.} The proof follows the reasoning of
\cite[Theorem~11.2]{polishchuk} with some modifications necessary
since we work in the context of formal deformation quantizations. The
essential difficulties are already dealt with in
Proposition~\ref{prop:param}, Theorem~\ref{thm:main} and
\cite[Theorem~11.2]{polishchuk} but some additional work is required
to package the argument properly. Given a formal analytic space
$\mathbb{Y} \to \bbD$ and a line bundle $L \to \mathbb{Y}\times_{\bbD}
\bbX_{\bpi}$ we need to construct a morphism $\mathbb{f} : \mathbb{Y}
\to {}_{\bB}\bbX^{\vee}$ and an isomorphism $(\mathbb{f}\times
\op{id})^{*}\mycal{P} \cong L$. Consider the integral transform
\[
\phi_{\mycal{Q}}^{[\mathbb{Y}\times_{\bbD} \bbX_{\bpi}\to
\mathbb{Y}\times_{\bbD} ({}_{\bB}\bbX^{\vee}) ]} :
\sff{D}^{b}_{c}(\mathbb{Y}\times_{\bbD} \bbX_{\bpi}) \to
\sff{D}^{b}_{c}(\mathbb{Y}\times_{\bbD} ({}_{\bB}\bbX^{\vee}),-1). 
\]  
and the object
$\phi_{\mycal{Q}}^{[\mathbb{Y}\times_{\bbD} \bbX_{\bpi}\to
\mathbb{Y}\times_{\bbD} ({}_{\bB}\bbX^{\vee}) ]}(L) \in
\sff{D}^{*}_{c}( \mathbb{Y}\times_{\bbD} ({}_{\bB}\bbX^{\vee}),
-1)$. By Proposition~\ref{prop:param} it follows that for every $\bbD$
point $\mathfrak{y} : \bbD \to \mathbb{Y}$ of $\mathbb{Y}$, the
pullback $\mathfrak{y}^{*}L$ is isomorphic to $(\mathbb{s}\times
\op{id})^{*}\mycal{P}$ for some $\bbD$-point $\mathbb{s} : \bbD \to
{}_{\bB}\bbX^{\vee}$. Specifying the $\bbD$-point $\mathbb{s}$ amounts
to specifying a $\bbD$-point $\mathfrak{s} : \bbD \to \bbX^{\vee}$
together with 
a trivialization of the gerbe
$\mathfrak{s}^{*}({}_{\bB}\bbX^{\vee})$. Now $(\mathfrak{y}\times
\op{id})^{*}\phi_{\mycal{Q}}^{[\mathbb{Y}\times_{\bbD} \bbX_{\bpi}\to
\mathbb{Y}\times_{\bbD} ({}_{\bB}\bbX^{\vee}) ]}(L) =
\phi_{\mycal{Q}}^{[\bbX_{\bpi}\to
{}_{\bB}\bbX^{\vee}]}((\mathbb{s}\times \op{id})^{*}\mycal{P})$, which
in turn is isomorphic to
$(\mathfrak{s},-1)_{*}\mathcal{O}_{\bbD}[-g]$, as explained in the
proof of Theorem~\ref{thm:main}. Taking homs into the various elements
of our orthogonal spanning class for
$\sff{D}^{*}_{c}({}_{\bB}\bbX^{\vee},-1)$ we conclude that the object
$\phi_{\mycal{Q}}^{[\mathbb{Y}\times_{\bbD} \bbX_{\bpi}\to
\mathbb{Y}\times_{\bbD} ({}_{\bB}\bbX^{\vee}) ]}(L)$ is of the form
$\mycal{F}[-g]$ for some $(-1)$-twisted sheaf $\mycal{F}$ on
$\mathbb{Y}\times_{\bbD}({}_{\bB}\bbX^{\vee})$. The key point of the
argument is to show that the stack-theoretic support of $\mycal{F}$ is
the graph of a map $\mathbb{f} : \bbY \to {}_{\bB}\bbX^{\vee}$ and that
$\mycal{F}$ is a line bundle on its support. As in the proof of 
\cite[Theorem~11.2]{polishchuk}, this will follow from the property 
\begin{equation} \label{eq:pointslices}
(\mathfrak{y}\times
\op{id})^{*}\mycal{F} \cong (\mathfrak{s},-1)_{*}\mathcal{O}_{\bbD},
\quad \text{for all $\bbD$-points \ } \mathfrak{y} : \bbD \to \bbY
\end{equation}
provided that we can show that
$\mycal{F}$ is finite and flat over $\bbY$. 

The
property \eqref{eq:pointslices}
implies that $\mycal{F}$ is a $(-1)$-twisted line bundle on its
support.  In other words, we can find a closed analytic subspace
$\mathfrak{i} : \bbS \hookrightarrow \bbY\times_{\bbD} \bbX^{\vee}$,
so that the stack-theoretic support of
$\mycal{F}$ is the gerbe $\bbS\times_{\bbY\times_{\bbD} \bbX^{\vee}}
[\bbY\times_{\bbD} ({}_{\bB}\bbX^{\vee})]$  and $\mycal{F}$
trivializes this gerbe. In particular $\mycal{F}$ is (non-canonically)
isomorphic to a sheaf of the form $(\mathfrak{i},-1)_{*}\mycal{G}$ for
some line bundle $\mycal{G}$ on $\bbS$. Thus $\mycal{F}$ will be  finite and
flat over $\bbY$ if and only if the sheaf $\mathfrak{i}_{*}\mycal{G}$
on $\bbY\times_{\bbD} \bbX^{\vee}$ is finite and flat over $\bbY$. 
To check this we will use the following:

\begin{lem} \label{lem:commutativize} Suppose $L$ is a line bundle on
  $\bbX_{\bpi}\times_{\bbD} \bbY$ which is of degree zero relative to
  $\bbY$. Let $\bbS \subset \bbY\times_{\bbD} \bbX^{\vee}$ and
  $\mycal{G}$ be as above. Then there exists a line bundle $M$ on
  $\bbX\times_{\bbD} 
  \bbY$ which is of degree zero relative to $Y$ and for which 
\[
\phi_{\mathcal{P}^{\vee}[[\hbar]]}^{[\mathbb{Y}\times_{\bbD} \bbX\to
\mathbb{Y}\times_{\bbD} \bbX^{\vee} ]}(M) \cong
\mathfrak{i}_{*}\mycal{G}[-g]. 
\]
\end{lem}
{\bf Proof.} Consider the natural projection and addition maps
\[
\xymatrix@1@M+1pc{
\bbY\times_{\bbD} \bbX \times_{\bbD} \bbX_{\bpi}
\ar@<.5ex>[r]^-{\bbp_{1,2+3}} \ar@<-.5ex>[r]_-{\bbp_{1,3}} &
\bbY\times_{\bbD} \bbX_{\bpi}.
}
\]
Here $\bbp_{1,3}$ is the projection onto the first and third factors
and  $\bbp_{1,2+3} = \bbp_{1}\times \mathbb{m}_{(0.1)}$ is the product of the
projection onto $\bbY$ and the ringed space map $\mathbb{m}_{(0,1)} :
\bbX\times_{\bbD} \bbX_{\bpi} \to \bbX_{\bpi}$ described in
section~\ref{sssec:groupstr}. 

Since by hypothesis $L \to \bbY\times_{\bbD} \bbX_{\bpi}$ is of degree
zero relative to $\bbY$, i.e. $L$ is translation invariant in the
$\bbX_{\bpi}$ direction, it follows that the line bundles
$\bbp_{1,3}^{*}L$ and $\bbp_{1,2+3}^{*}L$ satisfy the assumptions of
the see-saw principle Proposition~\ref{prop:see-saw}. Therefore we can
find a line bundle $M$ on $\bbY\times_{\bbD} \bbX$ so that 
\[
\bbp_{1,2}^{*}M\otimes \bbp_{1,3}^{*}L \cong \bbp_{1,2+3}^{*}L.
\]
Now a straightforward diagram chase shows that the Fourier-Mukai
transform of $M$ with respect to the Poincar\'{e} sheaf $\bbp_{2,3}^{*}
\mathcal{P}^{\vee}[[\hbar]]$ on $\bbY\times_{\bbD} \bbX \times_{\bbD}
\bbX^{\vee}$ is the object $\mathfrak{i}_{*}\mycal{G}[-g]$. \ \hfill
$\Box$ 

\

\medskip

With $M$ in hand we can now proceed to reason as in the proof of 
\cite[Theorem~11.2]{polishchuk}. Since $\bbY\times_{\bbD} \bbX \to
\bbY$ and $\bbY\times_{\bbD} \bbX^{\vee} \to
\bbY$ are dual family of complex tori, the usual argument
\cite[Section~11.2]{polishchuk} shows that
$\phi_{\mathcal{P}^{\vee}[[\hbar]]}^{[\mathbb{Y}\times_{\bbD} \bbX\to 
\mathbb{Y}\times_{\bbD} \bbX^{\vee} ]}$ is left-adjoint to 
$\phi_{\mathcal{P}[[\hbar]]}^{[\mathbb{Y}\times_{\bbD} \bbX^{\vee} \to 
\mathbb{Y}\times_{\bbD} \bbX ]}$, and that the adjunction map
\[
\op{Id} \to \phi_{\mathcal{P}[[\hbar]]}^{[\mathbb{Y}\times_{\bbD}
    \bbX^{\vee} \to  
\mathbb{Y}\times_{\bbD} \bbX ]}\circ
    \phi_{\mathcal{P}^{\vee}[[\hbar]]}^{[\mathbb{Y}\times_{\bbD}
    \bbX\to  
\mathbb{Y}\times_{\bbD} \bbX^{\vee} ]}
\]
is an isomorphism.

Now to check that $\mathfrak{i}_{*}\mycal{G}$ is finite and flat over
$\bbY$ we can assume that the closed fiber $Y \subset \bbY$ of $\bbY
\to \bbD$ is a  Stein space and check that the global sections of
$\mathfrak{i}_{*}\mycal{G}$ on $\bbY\times_{\bbD} \bbX$ are a
flat module of finite rank over the Stein algebra
$\Gamma(\bbY,\mathcal{O}_{\bbY})$.

For the global sections $\Gamma(\bbY\times_{\bbD} \bbX,
\mathfrak{i}_{*}\mycal{G})$ we use the above adjunction to compute
\[
\begin{split}
\Gamma(\bbY\times_{\bbD} \bbX,
\mathfrak{i}_{*}\mycal{G}) & = \op{Hom}\left(\mathcal{O}_{\bbY\times_{\bbD}
  \bbX}, \mathfrak{i}_{*}\mycal{G}\right) \\
& = \op{Hom}\left(\phi_{\mathcal{P}^{\vee}[[\hbar]]}^{[\mathbb{Y}\times_{\bbD}
    \bbX\to  
\mathbb{Y}\times_{\bbD} \bbX^{\vee} ]}(\mathcal{O}_{\bbY\times_{\bbD}
  \bbX}), M\right) \\
& = \op{Hom}\left((\op{id}\times
\mathfrak{o})_{*}\mathcal{O}_{\bbY}[-g],M\right) 
\\
& = \op{Ext}^{g}\left((\op{id}\times
\mathfrak{o})_{*}\mathcal{O}_{\bbY},M\right), 
\end{split}
\]
where $\mathfrak{o} : \bbD \to \bbX$ is the constant section
corresponding to the origin $o \in X$. Now on $\bbY\times_{\bbD} \bbX$
we have $(\op{id}\times
\mathfrak{o})_{*}\mathcal{O}_{\bbY} =
\bbp_{2}^{*}\mathfrak{o}_{*}\mathcal{O}_{\bbD}$ and since
$\mathfrak{o}_{*}\mathcal{O}_{\bbD}$ is supported on a section of
$\bbX \to \bbD$ we get
\[
\begin{split}
\op{Ext}^{g}\left((\op{id}\times
\mathfrak{o})_{*}\mathcal{O}_{\bbY},M\right) & =
\op{Ext}^{g}\left(\bbp_{2}^{*}\mathfrak{o}_{*}\mathcal{O}_{\bbD},M\right)
\\
& =
\op{Ext}^{g}\left(\mathfrak{o}_{*}\mathcal{O}_{\bbD},\bbp_{2*}M\right)
\\
& =
H^{0}\left(\bbX,\mycal{E}xt^{g}\left(\mathfrak{o}_{*}\mathcal{O}_{\bbD},
\bbp_{2*}M\right) \right).
\end{split}
\]
Now, it only remains to observe that since $M$ is a line bundle
$\bbp_{2*}M$ is a flat $\mathcal{O}_{\bbX}$-module, and so $
\mycal{E}xt^{g}\left(\mathfrak{o}_{*}\mathcal{O}_{\bbD}, 
\bbp_{2*}M\right) =
\mycal{E}xt^{g}\left(\mathfrak{o}_{*}\mathcal{O}_{\bbD}, 
\mathcal{O}_{\bbX}\right)\otimes_{\mathcal{O}_{\bbX}} \bbp_{2*}M =
\mathfrak{o}_{*}\mathcal{O}_{\bbD}\otimes_{\mathcal{O}_{\bbX}}
\bbp_{2*}M$. In other words 
\[
\Gamma(\bbY\times_{\bbD} \bbX,
\mathfrak{i}_{*}\mycal{G}) =
H^{0}\left(\bbX,\mathfrak{o}_{*}\mathcal{O}_{\bbD}
\otimes_{\mathcal{O}_{\bbX}}\bbp_{2*}M \right)  
\]
which is clearly finite and flat as an
$\Gamma(\bbY,\mathcal{O}_{\bbY})$-module. 

This implies that the sheaf $\mycal{F}$ on $\bbY\times_{\bbD}
({}_{\bB}\bbX^{\vee})$ is finite and flat over $\bbY$, which combined
with the property \eqref{eq:pointslices} implies that the support
$\mathbb{S}$ is the graph of a morphism $\mathbb{f} : \bbY \to
{}_{\bB}\bbX^{\vee}$. Hence $\phi_{\mycal{P}}^{[\bbY\times_{\bbD}
({}_{\bB}\bbX^{\vee}) \to \bbY\times_{\bbD} \bbX_{\bpi}]}(\mycal{F})$
is a line bundle on $\bbY\times_{\bbD} \bbX_{\bpi}$ which is
isomorphic to $(\op{id}\times \mathbb{f})^{*}\mycal{P}$ along the fibers of
$\bbY\times_{\bbD} \bbX_{\bpi} \to \bbY$. However in the proof of
Theorem~\ref{thm:main} we checked that the composition functor
$\phi_{\mycal{P}}^{[{}_{\bB}\bbX^{\vee} \to \bbX_{\bpi}]}\circ
\phi_{\mycal{Q}}^{[\bbX_{\bpi} \to {}_{\bB}\bbX^{\vee}]}$ is isomorphic to
$\op{id}[-g]$. This implies that the canonical adjunction morphism 
\[
\phi_{\mycal{P}}^{[\bbY\times_{\bbD}
({}_{\bB}\bbX^{\vee}) \to \bbY\times_{\bbD} \bbX_{\bpi}]}(\mycal{F})
\to L 
\]
is an isomorphism since it is an isomorphism over every $\bbD$ point
of $\bbY$. Applying again the see-saw principle
Proposition~\ref{prop:see-saw} we conclude that $L$ is isomorphic to 
$\bbp_{1}^{*}A\otimes (\op{id}\times \mathbb{f})^{*}\mycal{P}$ for a line
bundle $A$ on $\bbY$ which is unique up to a unique isomorphism. This
shows that
$\mathfrak{c}$ is an isomorphism of stacks and concludes the proof of
the theorem. \ \hfill $\Box$

\

\bigskip
\bigskip

\appendix

\

\bigskip

\noindent
{\Large\bfseries Appendices}

\

\noindent
In the following appendicies, we prove some general facts concerning
the module theory for deformation quantizations of Hausdorff analytic
spaces.  In some cases we specialize to complex tori, or the Moyal
quantization in particular.  It should be noted that many of the
proofs in these appendices carry over immediately to other contexts.
In particular we expect that our results will be applicable to
deformation quantizations in characteristic $p$ \cite{bk-fedosov}
and the case of $C^{\infty}$ real manifolds. In particular our
deformation theory analysis recovers some results from
\cite{bw-hermitian,bw-bimodule,bw-induction}. 

\

\bigskip

\

\Appendix{Quantum invertible sheaves} \label{app:quantum}

\

\noindent
Let $S$ be a Hausdorff analytic space. Let $\mathbb{S} := (S,\A)$
be a deformation quantization of $(S,\mathcal{O})$. We will assume
that $\A$ and $\mathcal{O}[[\hbar]]$ are isomorphic as sheaves of
$\Ch$-modules and that the product on $\A$ is given by a
$\star$-product on $\mathcal{O}[[\hbar]]$.  Suppose further that for
some subsheaf of $\mathbb{C}$-algebras $\mathcal{C} \subset
\mathcal{O}$, the left or right $\star$-multiplication action of
$\mathcal{C} [[\hbar]] \subset \A$ on $\A$ agree with the commutative
action of multiplication of $\mathcal{C} [[\hbar]] \subset
\mathcal{O}[[\hbar]]$ on $\mathcal{O}[[\hbar]]$.  In particular
$\mathcal{C} [[\hbar]] \subset \A$ is a subsheaf of central
$\Ch$-subalgebras in $A$.

To any locally free left rank one $\mathcal{C}[[\hbar]]$-module $M$ we
can associate a natural line bundle $L :=
\A\otimes_{\mathcal{C}[[\hbar]]} M$ on $\mathbb{S}$. In this way we
get a functor from the groupoid of locally free left rank one
$\mathcal{C}[[\hbar]]$-modules to the groupoid of line bundles on
$\mathbb{S}$. We would like to understand the map this functor induces
on isomorphism classes.

The isomorphism classes of locally free left rank one
$\mathcal{C}[[\hbar]]$-modules are in a natural bijection with the
cohomology group $H^1(S, \mathcal{C} [[\hbar]]^\times)$, whereas the
isomorphism classes of line bundles on $\mathbb{S}$ are in bijection
with the cohomology set $H^1(S, \A^\times)$.  Denote by $G$ the image
of $H^1(S, \mathcal{C} [[\hbar]]^\times)$ in $H^1(S, \A^\times)$, by
$I_0$ the image of $H^1(S,\mathcal{C}^{\times})$ in $H^1(S,
\mathcal{O}^\times)$, and by $I$ the image of $H^1(S,\mathcal{C})$ in
$H^1(S, \mathcal{O})$.  Then we seek to understand the structure of
$G$, in terms of $I_{0}$ and $I$.

\begin{lem}\label{lem:ql}
There is a natural bijection
$I_0 \times \prod_{k=1}^{\infty} I \to G$.
In particular $G$ is a commutative group.
\end{lem}
{\bf Proof.} 
First consider the isomorphism of sheaves of groups
\[
\xymatrix@R-1.5pc@C-1pc{
\op{Exp} : & \mathcal{C}^\times \times \prod_{k=1}^{\infty} \mathcal{C} 
\ar[r] &  \mathcal{C}[[\hbar]]^\times \\
& (a_0, a_1, a_2, a_3, \ldots)
\ar@{|->}[r] &  a_0 \exp(a_1 \hbar + a_2 \hbar^2 + a_3 \hbar^3 +
\cdots).
}
\]
\noindent
For any open covering  $\mathfrak{U}$ of $S$ we have an induced
isomorphism of groups
\[
\xymatrix@R-1pc@C-1pc@M+0.5pc{
\op{Exp}: & \check{Z}^{1}\left(S,\mathfrak{U}, \mathcal{C}^\times \right) 
\times \prod_{i=1}^{\infty}
\check{Z}^{1}(S,\mathfrak{U}, \mathcal{C}) \ar@{=}[d] 
\ar[r] & 
\check{Z}^{1}\left(S,\mathfrak{U}, \mathcal{C}[[\hbar]]^\times  \right)
\\
& \check{Z}^{1}\left(S, \mathfrak{U}, \mathcal{C}^\times \times
\prod_{k=1}^{\infty}   
\mathcal{C}\right) & 
}
\]
We also have a bijection of sets
\[
\xymatrix@R-1pc@C-1pc@M+0.5pc{
\op{Exp}: & \check{C}^{0}\left(S,\mathfrak{U}, \mathcal{O}^{\times}\right) 
\times \prod_{k=1}^{\infty}
\check{C}^{0}(S,\mathfrak{U}, \mathcal{O}) \ar@{=}[d] 
\ar[r] & 
\check{C}^{0}(S,\mathfrak{U}, \mathcal{O}[[\hbar]]^{\times} ) 
\ar[d]^-{\cong} \\
& \check{C}^{0}\left(S,\mathfrak{U}, \mathcal{O}^{\times} \times
\prod_{k=1}^{\infty} \mathcal{O} \right) & \check{C}^{0}(S, \mathfrak{U},
\A^{\times}) 
}
\]
given by the same formula (for the regular exponential, not the
$\star$-exponential).  Suppose now that we have two cocycles $a, a'
\in \check{Z}^{1}\left(S,\mathfrak{U}, \mathcal{C}^\times \right) 
\times \prod_{k=1}^{\infty}
\check{Z}^{1}(S,\mathfrak{U},  \mathcal{C})$ which satisfy
\[
a ({a'})^{-1} = \delta(b)
\]
in $\check{Z}^1\left(S, \mathfrak{U}, \mathcal{O}^\times\right) \times
\prod_{k=1}^{\infty} \check{Z}^{1}(S, \mathfrak{U}, \mathcal{O} )$,
with some 
$b \in \check{C}^0(S,\mathfrak{U}, \mathcal{O}^\times ) \times
\prod_{k=1}^{\infty} \check{C}^{0}(S, \mathfrak{U}, \mathcal{O} )$.  Then
we have $a_{ij} b_j = b_i a'_{ij}$ and so 
\[
\op{Exp}(a_{ij})\cdot \op{Exp}(b_j) = 
\op{Exp} {(b_i)}\cdot \op{Exp} {(a'_{ij})}
\]
Since $\op{Exp}(a_{ij})$ and 
$\op{Exp}(a'_{ij})$ are in $\mathcal{C}[[\hbar]]$, this means that
\[
(\op{Exp} (b_i))^{-1} \star 
\op{Exp} (a_{ij})\star \op{Exp} ({b_j}) = \op{Exp}(a'_{ij})
\]
in $\check{Z}^{1}(S, \mathfrak{U}, \A^{\times})$,
where the inverse here is taken with respect to the $\star$-product.
Therefore we 
have produced a well defined and surjective map
\[
\op{Exp}: I_{0} \times \prod_{k=1}^{\infty} I
\rightarrow G.
\]
Now suppose $a \in
\check{Z}^{1}(S, \mathfrak{U},  \mathcal{C}^\times ) 
\times \prod_{k=1}^{\infty}
\check{Z}^{1}(S, \mathfrak{U}, \mathcal{C})$ is a cocycle for which
$\op{Exp}(a) \in \check{Z}^{1}(S,\mathfrak{U},\A^{\times})$ is a
coboundary. Then 
$\op{Exp}(a_{ij}) = (\op{Exp}{(a)})_{ij} = c_i \star c_j^{-1}$ for some
$c \in \check{C}^{0}(S, \mathfrak{U} , \A^{\times} )$.  Consequently
$(\op{Exp}{(a_{ij})}) c_{j} = c_{i}$ and so if we choose
$b \in \check{C}^0(S, \mathfrak{U}, \mathcal{O}^\times )
\times \prod_{k=1}^{\infty}
\check{C}^{0}(S,\mathfrak{U}, \mathcal{O})$ satisfying 
$\op{Exp}{(b)} = c$, then
we have $a = \delta(b)$.  Therefore $\op{Exp}$ is in fact an isomorphism.
By taking the limit over all coverings $\mathfrak{U}$,
we arrive at the desired statement. \ \hfill $\Box$

\

\begin{rem}
One important example of the above is the case when $S$ is a complex
manifold, $\A$ is a deformation quantization of the structure sheaf
$\mathcal{O}_{S}$ which is globally isomorphic to
$\mathcal{O}_{S}[[\hbar]]$ and $\bpi$ is the corresponding global
holomorphic Poisson structure on $S$.  In this case we can take
$\mathcal{C} \subset \mathcal{O}$ to be the subsheaf of functions
constant along the leaves of the foliation of $S$ by symplectic
leaves. Explicitly $\mathcal{C}$ is the kernel of the map $\mathcal{O}
\to T_{S}$ of $\mathbb{C}$-sheaves given by the composition of the de
Rham differential $d: \mathcal{O} \to \Omega_{S, \op{cl}}^{1} \subset
\Omega_{S}^{1}$ with the contraction $\bpi \cntrct : \Omega_{S}^{1}
\to T_{S}$.  If we further define $\Omega_{S, \bpi, \op{cl}}^{1}
\subset \Omega_{S, \op{cl}}^{1}$ as the kernel of the map $\bpi
\cntrct : \Omega_{S, \op{cl}}^{1} \to T_{S}$, then we have a short
exact sequence of sheaves of groups
\[
{1} \to \mathcal{C}[[\hbar]]^{\times} \to \mathcal{O}[[\hbar]]^{\times} 
\to (\Omega_{S, \op{cl}}^{1}/\Omega_{S, \bpi, \op{cl}}^{1})[[\hbar]]
\to {0}. 
\]
Combining the relevant part of the associated 
long exact sequence of cohomology groups with Lemma~\ref{lem:ql}
we obtain a five term exact sequence of groups 
\begin{equation} \label{eq:fiveterm}
\les
{H^{0}\left(S, \mathcal{C}[[\hbar]]^{\times}\right)}
{H^{0}\left(S,\A^{\times}\right)}
{H^{0}\left(S, (\Omega_{S, \op{cl}}^{1}/\Omega_{S, \bpi,
    \op{cl}}^{1}\right)[[\hbar]]}
{H^{1}\left(S, \mathcal{C}[[\hbar]]^{\times}\right)} 
{I_0 \times \prod_{k=1}^{\infty} I}
\end{equation}
where the identification $I_0 \times \prod_{k=1}^{\infty} I \cong G
\subset H^{1}(S, \A^{\times})$ is described in Lemma~\ref{lem:ql}, and
we have used the natural identification
\[\op{image}[H^{1}(S, \mathcal{C}[[\hbar]]^{\times}) 
\to H^{1}(S, \A^{\times})] 
\cong \op{image}[H^{1}(S, \mathcal{C}[[\hbar]]^{\times}) 
\to H^{1}(S, \mathcal{O}[[\hbar]]^{\times})].
\]
Note that often the first two groups in \eqref{eq:fiveterm} are
isomorphic, e.g. for a compact $S$. In such cases, the sequence
\eqref{eq:fiveterm} reduces to a short exact sequence and we get a
concrete description of $H^{1}\left(S,
\mathcal{C}[[\hbar]]^{\times}\right)$. 

For line bundles on a quantization of a complex torus we will obtain
more specific results in Lemma \ref{lem:quantizable} using factors of
automorphy on the universal cover, instead of \u{C}ech cohomology.
  
\end{rem}

\

\bigskip
\bigskip

\Appendix{Quantizing vector bundles}
\label{app:trivial}

\

\noindent
In this appendix, we provide the promised proof of the triviality of
locally free left $\Ah{V}$-modules where $\Ah{V}$ is a Moyal
quantization of a vector space $V$.  In fact we prove a more general
statement, which is applicable in many contexts that arise in
deformation quantization.

Let $S$ be a Hausdorff analytic space and let $\A$ be a sheaf of
$\mathbb{C}[[\hbar]]$-algebras such that $\A
\otimes_{\mathbb{C}[[\hbar]]} \mathbb{C} \cong \mathcal{O}_{S}$.
Consider the multiplicative sheaf of groups
$\cAut_{\A-\op{mod}}(\A^{\oplus m})$ of left $\A$-module automorphisms
of $\A^{\oplus m}$.  Isomorphism classes of left $\A$-modules which
are locally free of rank $m$ are in a natural bijective correspondence
with $H^1(S, \cAut_{\A-\op{mod}}(\A^{\oplus m}))$.  It turns out that
the vanishing of the first cohomology of the multiplicative sheaf of
groups $\cAut_{\mathcal{O}}(\mathcal{O}^{\oplus m})$ and the first
cohomology of the additive sheaf of groups
$\cEnd_{\mathcal{O}}(\mathcal{O}^{\oplus m})$ guarantee that all such
left $\A$-modules are isomorphic to the trivial one.  Below we will be
treating only left modules, and so we will simply write
$\cAut_{\A}(\A^{\oplus m})$ in place of
$\cAut_{\A-\op{mod}}(\A^{\oplus m})$ and also $\cEnd_{\A}(\A^{\oplus
m})$ in place of $\cEnd_{\A-\op{mod}}(\A^{\oplus m})$.

\begin{lem} \label{lem:acyclic}
If $H^1(S, \cEnd_{\mathcal{O}}(\mathcal{O}^{\oplus m})) = \{0\}$ and
$H^1(S, \cAut_{\mathcal{O}}(\mathcal{O}^{\oplus m}))=\{1\}$ then the
set \linebreak
$H^1(S, \cAut_{\A-\op{mod}}(\A^{\oplus m}))$ has only one element.
\end{lem}
{\bf Proof.}  Fix a cofinal system of coverings
$\boldsymbol{\mathfrak{S}}$, such that if $\mathfrak{U}$ is any
covering in $\boldsymbol{\mathfrak{S}}$ we have
$\check{H}^1(S,\mathfrak{U}, \cEnd_{\mathcal{O}}(\mathcal{O}^{\oplus
m})) = \{0\}$, $\check{H}^1(S,\mathfrak{U},
\cAut_{\mathcal{O}}(\mathcal{O}^{\oplus m}))=\{1\}$, and also such
that for any covering $\mathfrak{U}$ in $\boldsymbol{\mathfrak{S}}$
and for any open set $U$ in the covering $\mathfrak{U}$, there is a
splitting (a morphism of sheaves of $\mathbb{C}$-vector spaces)
$\sigma_{U}$ over $U$ of the projection $\rho: \A \to \mathcal{O}$.
Then if we have a covering $\mathfrak{U} = \{U_{i} | i \in I \}$ in
the system $\boldsymbol{\mathfrak{S}}$ it is enough to show that any
\u{C}ech 1-cochain $G \in \check{Z}^1(S,\mathfrak{U},
\cAut_{\A}(\A^{\oplus m}))$ is a coboundary.  We will do this by
constructing a sequence $\{ G =: G^{0}, G^{1}, G^{2}, \dots \}$  of
cocycles in $\check{Z}^1(S,\mathfrak{U},
\cAut_{\A}(\A^{\oplus m}))$ with the following 
properties: 
\begin{itemize}
\item $G^{j+1}$ is cohomologous to $G^{j}$ for $j \geq 0$;
\item $\rho(G^{j}) = \op{id}$ for $j \geq 1$;
\item $\rho((G^{j} - \op{id})/\hbar^{k}) = 0$ for every $j \geq 1$ and
  every $k = 1, \dots , j-1$.   
\end{itemize}
To carry out the construction chose splittings $\sigma_i :=
\sigma_{U_{i}} : \mathcal{O}_{|U_{i}} \to \A_{|U_{i}}$ over each
$U_{i}$ as above.  The cocycle relationship says that $G_{ij} \in
\cAut_{\A}(\A^{\oplus m})(U_{ij})$ satisfy
\[
G_{ik} = G_{ij} \circ G_{jk}
\]
Reducing modulo $\hbar$ we see that
\[
\rho(G) \in
\check{Z}^{1}(S, \mathfrak{U},\cAut_{\mathcal{O}}(\mathcal{O}^{\oplus n})).
\]
\noindent
Since $\check{H}^{1}(S,\cAut_{\mathcal{O}}(\mathcal{O}^{\oplus m})) =
0$ we can chose
\[
\phi^{0} = \{ \phi^{0}_{i} \}_{i\in I} \in \check{C}^0(S, \mathfrak{U}, 
\cAut_{\mathcal{O}}(\mathcal{O}^{\oplus m}))
\]
satisfying
\[
{\left({\phi^{0}_i}\right)}^{-1} \circ  \rho(G_{ij}) \circ {\phi^{0}_j}
= 1_{\mathcal{O}|_{U_{ij}}}.
\]  
Consider the cochain
$G^{1} \in \check{C}^1(S,\mathfrak{U},
\cAut_{\A}(\A^{\oplus m}))$ defined by
\[
G^{1}_{ij}= {\left(\sigma_{i}\left(\phi^{0}_i\right)\right)}^{-1} 
\circ  G_{ij} \circ {\left(\sigma_{j}\left(\phi^{0}_j\right)\right)}.
\]  
Clearly $G^{1}$ is a cocycles which is cohomologous to $G$. By
construction $G^{1}$ 
satisfies $\rho(G^{1}) = \op{id}$ and therefore $\rho(G^{1} - \op{id})
= 0$.  

Assume that by induction we have produced 
$G^{1}, \dots, G^{n}$ which are
cohomologous to each other and 
such that for each $1 \leq j \leq n$ we have
$\rho(G^{j}) = \op{id}$ and for $k < j$ we have
$\rho((G^{j} - \op{id})/\hbar^{k}) = 0$.  Now we have
\begin{equation} \label{eq:nowwehave}
(G^{n}_{ij} - \op{id}) \circ (G^{n}_{jk} - \op{id}) = 
G^{n}_{ik} - G^{n}_{ij} - G^{n}_{jk} + \op{id} = 
(G^{n}_{ik} - \op{id})-(G^{n}_{ij} - \op{id})-(G^{n}_{jk} - \op{id}).
\end{equation}
Notice that 
\[\rho\left(\left((G^{n}_{ij} - \op{id}) \circ (G^{n}_{jk} - 
\op{id})\right)/\hbar^{n}\right) = 
\rho\left(\left(G^{n}_{ij} - \op{id} \right)/\hbar^{n-1} \right) 
\circ \rho\left(\left(G^{n}_{jk} - \op{id} \right)/\hbar \right) = 0.  
\]
\noindent
Therefore dividing by $\hbar^{n}$ and applying $\rho$ to 
equation \eqref{eq:nowwehave} we conclude that
\[\rho((G^{n}_{ik} - \op{id})/\hbar^{n})
= \rho((G^{n}_{ij} - \op{id})/\hbar^{n}) + \rho((G^{n}_{jk} -
\op{id})/\hbar^{n}) 
\]
In other words $\rho((G^{n} - \op{id})/\hbar^{n})
\in \check{Z}^1(S,\mathfrak{U},  
\cEnd_{\mathcal{O}}(\mathcal{O}^{\oplus m}))$.  Therefore, taking into
account the fact that $\check{H}^{1}(S,
\mycal{E}nd_{\mathcal{O}}(\mathcal{O}^{\oplus m})) = 0$, we can
chose a cochain
\[
\phi^{n} = \{ \phi^{n}_{i} \}_{i \in I} \in \check{C}^{0}(S,\mathfrak{U}, 
\cEnd_{\mathcal{O}}(\mathcal{O}^{\oplus m}))
\]
satisfying
\[\phi^{n}_{i} - \phi^{n}_{j} = \rho((G^{n}_{ij} - \op{id})/\hbar^{n})
\]
Consider now the cocycle
$G^{n+1} \in \check{Z}^1(S,\mathfrak{U},
\cAut_{\A}(\A^{\oplus m}))$ defined by
\[
G^{n+1}_{ij} = \left(\op{id} + \hbar^{n}
\sigma_{i}\left(\phi^{n}_{i}\right)\right)^{-1} 
\circ G^{n}_{ij} \circ \left(\op{id} + \hbar^{n}
\sigma_{j}\left(\phi^{n}_{j}\right)\right) 
\]
Then modulo $\hbar^{n+1}$ we have that
$\left(\op{id} + \hbar^{n} \sigma_{i}\left(\phi^{n}_{i}\right)\right)^{-1}$
is equivalent to $\left(\op{id} - \hbar^{n}
\sigma_{i}\left(\phi^{n}_{i}\right)\right)$ 
and so  $G^{n+1}_{ij}$  is cohomologous to
$\left(\op{id} - \hbar^{n} \sigma_{i}\left(\phi^{n}_{i}\right) + (G^{n}_{ij}
- \op{id}) 
+ \hbar^{n} \sigma_{j}\left(\phi^{n}_{j}\right)\right)$ 
where $(G^{n}_{ij} -
\op{id}) \in \linebreak
\check{C}^{1}(\mathfrak{U},\hbar^{n}\cEnd_{\A}(\A^{\oplus m}))$.
\noindent
Clearly $G^{n+1}$ is cohomologous to $G^{n}$, 
$\rho(G^{n+1}) = \op{id}$ and $\rho((G^{n+1} -
\op{id})/\hbar^{k}) = 0$ 
for $0 \leq k < n$.  Finally,
\[
\rho((G^{n+1} - \op{id})/\hbar^{n}) = \rho((G^{n}_{ij} -
\op{id})/\hbar^{n}) - 
\phi^{(n)}_{i} + \phi^{(n)}_{j} = 0.
\] 
Therefore, we have produced the required element $G^{n+1}$ completing
the induction step and the proof of the 
lemma.
\ \hfill $\Box$

\

\medskip

\noindent
The previous lemma shows that if there are no non-trivial rank $m$
vector bundles on $S$ and if the trivial rank $m$ vector bundle has no
infinitesimal deformations, then there are no non-trivial quantum
vector bundles of rank $m$ on the deformation quantization $\mathbb{S}
:= (S,\A)$.

Next we will discuss the obstructions for quantizing a rank $m$
holomorphic vector bundle $W$ on $S$ along a given deformation
quantization $\mathbb{S} = (S,\A)$. We will do this order-by-order in
the formal parameter $\hbar$.  One may consider these investigations as 
first steps in understanding $\A-\op{mod}$ as the formal deformation of 
the abelian category $\mathcal{O}_{S}-\op{mod}$.  Hopefully, this will 
eventually lead to an interpretation in terms of a more systematic study of 
the deformation theory of abelian categories.  This study has appeared
for instance in \cite{VandenBergh/Lowen}.  

\

As before, we will write
$\cAut_{\A/\hbar^{q}\A}({(\A/\hbar^{q}\A)}^{\oplus m})$ instead of
$\cAut_{{\A/\hbar^{q}\A}-\op{mod}}({(\A/\hbar^{q}\A)}^{\oplus m})$. 
Consider the projection map 
\[
\op{proj_{q}}: H^{1}(S, 
\cAut_{\A/\hbar^{q}\A}({(\A/\hbar^{q}\A)}^{\oplus m})) 
\to H^{1}(S, \cAut_{\mathcal{O}}(\mathcal{O}^{\oplus m})).
\]
We define the set of length $q$ quantizations of $W$ along
$\mathbb{S}$ as 
\[
\op{Quant}_{q}(W) := \op{proj}_{q}^{-1}([W]) 
\]
\

\begin{lem} \label{lem:extend}
Let $W$ be a holomorphic vector bundle on a Hausdorff analytic space
$S$.  Then  
there is a map 
\[
\op{ob}_{n+1} : \op{Quant}_{n+1}(W) \to H^{2}(S, \cEnd(W))
\]  
which measures the obstruction for a length $n+1$ quantization of $W$
to prolong to a length $n+2$ quantization. The ambiguity in choosing 
such a prolongation is given by 
$H^{1}(S, \cEnd(W))$.  In other words we have an exact sequence of sets
\[
\op{Quant}_{n+2}(W)
\to 
\op{Quant}_{n+1}(W)
\stackrel{\op{ob}_{n+1}}{\to} H^{2}(S, \cEnd(W))
\] 
and a free action of the additive group $H^{1}(S, \cEnd(W))$ on
$\op{Quant}_{n+2}(W)$, so that 
\[
\op{Quant}_{n+2}(W)/H^{1}(S,
\cEnd(W)) = \op{im}[\op{Quant}_{n+2}(W) \to \op{Quant}_{n+1}(W)].
\] 
\end{lem}
{\bf Proof.}  We first define the map $\op{ob}_{n+1}:
\op{Quant}_{n+1}(W) \to H^{2}(S, \cEnd(W))$ and prove that
$\op{ob}_{n+1}^{-1}(0)$ is the image of $\op{Quant}_{n+2}(W)$.  Fix a
fine enough open cover $\mathfrak{U} = \{U_{i} | i \in I \}$ of
$S$. Represent $[W] \in H^{1}(S,
\cAut_{\mathcal{O}}(\mathcal{O}^{\oplus{m}}))$ by a cocycle $g \in
\check{Z}^{1}(S, \mathfrak{U},
\cAut_{\mathcal{O}}(\mathcal{O}^{\oplus{m}}))$ where $g_{ij} = \mu_{i}
\circ \mu_{j}^{-1}$ for some trivializations $\mu_{i}:W_{{|}_{U_{i}}}
\to \mathcal{O}^{\oplus m}_{{|}_{U_{i}}}$.  Similarly, we represent $G
\in \op{Quant}_{n+1}(W)$ by an element $\{G_{ij} \} \in
\check{Z}^{1}(S, \mathfrak{U},
\cAut_{\A/\hbar^{n+1}\A}({(\A/\hbar^{n+1}\A)}^{\oplus m}))$

Define $\op{ob}_{n+1}(G)$ to be the element of $H^{2}(S, \cEnd(W))$
induced by taking the limit over all open covers of the elements $\{
\op{ob}_{n+1}(G)_{ijk} \} \in \check{Z}^{1}(S, \mathfrak{U},
\cEnd_{\mathcal{O}}(W))$ where
\[
\op{ob}_{n+1}(G)_{ijk} = \mu_{k}^{-1} \circ \rho\left(
\frac{\widetilde{G}_{ki} \circ  
\widetilde{G}_{ij} \circ \widetilde{G}_{jk} - \op{id}}{ \hbar^{n+1}} \right)
\circ  \mu_{k} 
= \mu_{i}^{-1} \circ \rho\left(  \frac{\widetilde{G}_{ij}\circ 
\widetilde{G}_{jk} - \widetilde{G}_{ik}}{\hbar^{n+1}}  \right) \circ \mu_{k}.
\]
Here $\rho : \cAut_{\A/\hbar^{n+1}}\left(\A^{\oplus
m}/\hbar^{n+1}\right) \to \cAut_{\mathcal{O}}(\mathcal{O}^{\oplus m})$
is the reduction modulo $\hbar$, and $\{ \widetilde{G}_{ij} \}$ is a lift
of $\{ G_{ij} \}$ to $\check{C}^{1}(S,\mathfrak{U},
\cAut_{\A/\hbar^{n+2}\A}({(\A/\hbar^{n+2}\A)}^{\oplus m}))$.  We will
check that $\{ \op{ob}_{n+1}(G)_{ijk} \}$ is closed. The check that
the above definition of $\op{ob}_{n+1}$ is independent of all choices
made in the construction is easy but tedious and is left to the
reader.
\[
\begin{split}
(\delta & \{ \op{ob}_{n+1}(G)_{abc} \})_{ijkl} = \op{ob}_{n+1}(G)_{jkl} - 
\op{ob}_{n+1}(G)_{ikl} + \op{ob}_{n+1}(G)_{ijl} -
\op{ob}_{n+1}(G)_{ijk} \\
& = \mu_{j}^{-1} \circ \rho\left(\frac{\widetilde{G}_{jk} \circ 
\widetilde{G}_{kl} - \widetilde{G}_{jl}}{\hbar^{n+1}}\right) \circ \mu_{l}
- \mu_{i}^{-1} \circ \rho\left(\frac{\widetilde{G}_{ik} \circ 
\widetilde{G}_{kl} - \widetilde{G}_{il}}{\hbar^{n+1}}\right) \circ \mu_{l} \\
& +\mu_{i}^{-1} \circ \rho\left(\frac{\widetilde{G}_{ij} \circ 
\widetilde{G}_{jl} - \widetilde{G}_{il}}{\hbar^{n+1}}\right) \circ \mu_{l}
-\mu_{i}^{-1} \circ \rho\left(\frac{\widetilde{G}_{ij} \circ 
\widetilde{G}_{jk} - \widetilde{G}_{ik}}{\hbar^{n+1}}\right) \circ \mu_{k} 
\end{split}
\]
By making the substitution 
$\widetilde{G}_{ik} = \widetilde{G}_{ij} \circ \widetilde{G}_{jk} 
+ (\widetilde{G}_{ik} -  \widetilde{G}_{ij} \circ \widetilde{G}_{jk})$ 
in the second term, we can rewrite this term as 
\[
-\mu_{i}^{-1} \circ \rho 
\left(\frac{\widetilde{G}_{ij} \circ \widetilde{G}_{jk}\circ\widetilde{G}_{kl} - 
\widetilde{G}_{il}}{\hbar^{n+1}}\right) \circ \mu_{l} 
-\mu_{i}^{-1} \circ \rho\left(\frac{\widetilde{G}_{ik} - 
\widetilde{G}_{ij} \circ \widetilde{G}_{jk}}{\hbar^{n+1}}\right) \circ \mu_{k}.
\]
Similarly, by making the substitution 
$\widetilde{G}_{jl} = \widetilde{G}_{jk} \circ \widetilde{G}_{kl} 
+ (\widetilde{G}_{jl} -  \widetilde{G}_{jk} \circ \widetilde{G}_{kl})$ 
into the third term, we can rewrite this term as 
\[
\mu_{i}^{-1} \circ \rho 
\left (\frac{\widetilde{G}_{ij} \circ \widetilde{G}_{jk}\circ\widetilde{G}_{kl} 
- \widetilde{G}_{il}}{\hbar^{n+1}} \right) \circ \mu_{l}
+ \mu_{j}^{-1} \circ \rho\left(\frac{\widetilde{G}_{jl} - 
\widetilde{G}_{jk} \circ \widetilde{G}_{kl}}{\hbar^{n+1}}\right) \circ \mu_{l}.
\]
It is now clear that all the terms cancel, and we have shown that
$\delta \{\op{ob}_{n+1}(G)_{abc}\} = 0$.  Clearly if $G$ comes from
$\op{Quant}_{n+1}(W)$ then we can choose the cochain $\{
\widetilde{G}_{ij} \}$ to be a cocycle, and so the cohomology class
$\op{ob}_{n+1}(G)$ is $0$ in $H^{2}(S,\cEnd_{\mathcal{O}}(W))$.
Conversely, if $\op{ob}_{n+1}(G) = 0$ then for a fine enough
$\mathfrak{U}$ we can find an element $C \in \check{C}^{1}(S, \mathfrak{U},
\cEnd_{\mathcal{O}}(W))$, with $\delta(C)_{ijk} =
\op{ob}_{n+1}(G)_{ijk}$.  Then we can define
a new element $\widetilde{G}' \in \check{C}^{1}(S, \mathfrak{U},
\cAut_{\A/\hbar^{n+2}\A}({(\A/\hbar^{n+2}\A)}^{\oplus m}))$ by the rule
\[
\widetilde{G}'_{ij} = \widetilde{G}_{ij} - \hbar^{n+1}
\sigma_{i}(\mu_{i} \circ C_{ij} \circ  
\mu_{j}^{-1}).
\]
Here, $\sigma_{i}$ is our chosen $\mathbb{C}-$module splitting of
$\rho:\A \to \mathcal{O}$ over $U_{i}$.  Notice that $\widetilde{G}'$
maps to $G$ modulo $\hbar^{n+1}$ and so in order to see that it is
closed, we merely observe the following vanishing:
\[
\begin{split}
\rho\left(\left((\widetilde{G}_{ij}\right.\right. & -
\hbar^{n+1} \sigma_{i}(\mu_{i}  
\circ C_{ij} \circ \mu_{j}^{-1})) 
\circ (\widetilde{G}_{jk} - \hbar^{n+1}\sigma_{j}(\mu_{j} \circ C_{jk}
\circ \mu_{k}^{-1}))  \\
& \left.\left. - (\widetilde{G}_{ik} - \hbar^{n+1} \sigma_{i}(\mu_{i} 
\circ C_{ik} \circ \mu_{k}^{-1}))\right)/\hbar^{n+1}\right) \\
& = \mu_{i} \circ \op{ob}_{n+1}(G)_{ijk} \circ \mu_{k}^{-1} - 
g_{ij} \circ \mu_{j} \circ C_{jk} \circ \mu_{k}^{-1} - 
\mu_{i} \circ C_{ij} \circ \mu_{j}^{-1} \circ g_{jk} + 
\mu_{i} \circ C_{ik} \circ \mu_{k}^{-1} \\
& = \mu_{i} \circ \left( \op{ob}_{n+1}(G)_{ijk} - C_{jk} - C_{ij} 
+ C_{ik} \right) \circ \mu_{k}^{-1} = \mu_{i} \circ 
\left( (\op{ob}_{n+1}(G) - \delta(C))_{ijk} \right) \circ \mu_{k}^{-1} \\ 
& = 0.
\end{split}
\]
The group $H^{1}(S, \cEnd_{\mathcal{O}}(W))$ 
acts on $\op{Quant}_{n+2}(W)$ by 
\[
\{K_{ij} \} \mapsto \{ K_{ij} \circ 
(\op{id} - \hbar^{n+1} \sigma_{i}(\mu_{j} \circ h_{ij} \circ \mu_{j}^{-1})) \}
=\{K_{ij} - \hbar^{n+1} \sigma_{i}(\mu_{i} \circ h_{ij} \circ \mu_{j}^{-1})) \}
\]
for $h \in H^{1}(S, \cEnd_{\mathcal{O}}(W))$ and $K \in \op{Quant}_{n+2}(W)$.  
This action is clearly free and preserves the fibers of the map 
$\op{Quant}_{n+2}(W) \to \op{Quant}_{n+1}(W)$.  In order to see that
 it is transitive, consider 
two elements $K, K' \in \op{Quant}_{n+2}(W)$ in the same 
fiber.  They define an unique element of 
$h \in H^{1}(S, \cEnd_{\mathcal{O}}(W))$ by the formula
\[h_{ij} = \mu_{j}^{-1} \circ \rho \left( \frac{\op{id} -K_{ji} 
\circ (K'_{ji})^{-1}}{\hbar^{n+1}} \right) \circ \mu_{j}.
\]
and it is easily seen that $h$ maps $K$ to $K'$.  
In order to see that $h$ is closed, we simply calculate
\[
\begin{split}
h_{ij} & = \mu_{j}^{-1} \circ \rho \left(\frac{K_{jk} \circ 
(K_{kj} - K_{kj} \circ K_{ji} \circ (K'_{ji})^{-1})}{\hbar^{n+1}} \right) 
\circ \mu_{j} \\
&  =  \mu_{k}^{-1} \circ \rho \left(\frac{ 
K_{kj} - K_{ki} 
\circ (K'_{ji})^{-1}}{\hbar^{n+1}} \right) \circ \mu_{j} \\
& = \mu_{k}^{-1} \circ \rho \left(\frac{ 
(K_{kj} - K_{ki} 
\circ (K'_{ji})^{-1}) \circ K'_{jk}}{\hbar^{n+1}} \right) \circ
\mu_{k} \\
& = \mu_{k}^{-1} \circ \rho \left(\frac{ 
K_{kj} \circ (K'_{kj})^{-1} - K_{ki}\circ (K'_{ki})^{-1}}{\hbar^{n+1}}
\right)  \circ \mu_{k} \\
& = h_{ik} - h_{jk}.
\end{split}
\]
This concludes the proof of the lemma.
\ \hfill $\Box$

\

\medskip

\begin{rem} \label{rem:principal} The set
  $\op{Quant}_{q}(W)$ can be naturally identified with the 
set of isomorphism classes of objects in a complex analytic stack.
The map $\op{Quant}_{n+2}(W) \to \op{Quant}_{n+1}(W)$ is induced from
a morphism of stacks and the action of $H^{1}(S,
\cEnd_{\mathcal{O}}(W))$  on $\op{Quant}_{n+2}(W)$ can be refined to
an analytic action on the stack corresponding to
$\op{Quant}_{n+2}(W)$. 
\end{rem}

\

\begin{rem}
The above lemma implies that if $W$ is a classical vector bundle for
which $H^{2}(S, \cEnd(W)) \cong \{0\}$ then $W$ extends to a locally
free left $\A$ module of the same rank.  At each stage of extension,
the ambiguity is precisely $H^{1}(S, \cEnd(W))$.  However, the
vanishing of the obstruction space is not necessary for
quantizability. There are many bundles with non-trivial obstruction
spaces which quantize to all orders.  For instance any flat bundle
does.
\end{rem}

\

\bigskip
\bigskip

\Appendix{The quantum Appell-Humbert theorem}

\

\noindent
We now focus on the case of the Moyal deformation quantization of the
sheaf of holomorphic functions on a complex torus.  Here we can use
factors of automorphy to obtain more precise formulae for the various
obstruction maps.  We begin by computing explicitly the 
obstruction $\op{ob}_{0}(W)$ for a line bundle $W$ on $X$ to quantize
to first order 
in $\hbar$. The assignment $W \mapsto \op{ob}_{0}(W)$ gives rise to a
short exact sequence of pointed cohomology sets
\begin{equation} \label{eq:obs0-sequence}
\xymatrix@1{
H^{1}(X, (\A/\hbar^{2}\A)^{\times}) \ar[r] &  H^{1}(X, \mathcal{O}^{\times})
\ar[r]^-{\op{ob}_{0}} & H^{2}(X,\mathcal{O}).
}
\end{equation}
If $W$ is a line bundle on $X$, then the `relative to $W$' part of
this sequence is precisely the sequence appearing in the statement of
Lemma~\ref{lem:extend}. Indeed $\op{Quant}_{0}(W) = \{ [W] \} \subset
H^{1}(X, \mathcal{O}^{\times})$, and $\op{Quant}_{1}(W)$ is just the
fiber of the map $H^{1}(X, (\A/\hbar^{2}\A)^{\times}) \to H^{1}(X,
\mathcal{O}^{\times})$ over the point $[W]$.

We can understand the sequence \eqref{eq:obs0-sequence} in terms of
the group cohomology of $\Lambda$ acting on functions on the universal
cover $V$. Taking into account the fact that
$H^{1}(V,(\A/\hbar^{2}\A)^{\times})) = 0$ (see the proof of
Lemma~\ref{lem:acyclic}) we can rewrite \eqref{eq:obs0-sequence} as
the exact sequence of pointed group cohomology sets:
\[
H^{1}(\Lambda, H^{0}(V,(\A/\hbar^{2}\A)^{\times})) \to 
H^{1}(\Lambda, H^{0}(V, \mathcal{O}^{\times})) \to 
H^{2}(\Lambda, H^{0}(V, \mathcal{O})).
\]
Now recall from  \eqref{cocyclecondition} that 
$Z^{1}(\Lambda, H^{0}(V,(\A/\hbar^{2}\A)^{\times}))$ consists 
of maps 
\[
\phi = \phi_{0} + \hbar \phi_{1} : 
\Lambda \to H^{0}(V, (\A/\hbar^{2}\A)^{\times})
\] 
satisfying 
\begin{equation} \label{eq-closedcond}
(\delta \phi)(\lambda_{1}, \lambda_{2})= 
\phi(\lambda_{1} + \lambda_{2}) - \phi(\lambda_{2}) \star 
(\phi(\lambda_{1}) \cdot \lambda_{2}). 
\end{equation}
Two cocycles $\phi$ and $\psi$ are cohomologous if there 
exists $f \in  H^{0}(V, (\A/\hbar^{2}\A)^{\times})$ which 
satisfies, for all $\lambda \in \Lambda$, the relationship 
\begin{equation} \label{eq:coboundary}
\psi(\lambda) = f^{-1} \star \phi(\lambda) \star (f \cdot \lambda).
\end{equation}
In the following, we will use the notation
\[
f \star g = fg + \sum_{j=1}^{\infty} \hbar^{j} (f \star g)_{j},
\]
for the components of a  star product.

Suppose now we are given a holomorphic line bundle $W$ on $X$
represented by a particular cocycle $\phi_{0} \in
Z^{1}(\Lambda,H^{0}(V,\mathcal{O}^{\times}))$. By the classical
Appell-Humbert theorem we can always replace $\phi_{0}$ by a
cohomologous cocycle which is given by the Appell-Humbert
formula:
\begin{equation} \label{eq-whoisreadingthis?}
\AH_{(H,\chi)}(\lambda)(v) =  \chi(\lambda)
\op{exp}\left(\pi H(v, \lambda) +  
\frac{\pi}{2}H(\lambda, \lambda)\right).
\end{equation}
Here $H$ is an element in the Neron-Severi group of $X$ thought of as
a Hermitian form on $V$ which satisfies $\op{Im} H(\Lambda, \Lambda)
\subset \mathbb{Z}$, and $\chi$ is an $H$-semi-character of
$\Lambda$.

Denote by  $\mathcal{P}(\Lambda)$ the group of all pairs $(H,\chi)$
where $H \in NS(X)$ and $\chi$ is a semicharacter for $H$
\cite[Section~1.2]{lange}.  By the Appell-Humbert theorem, the assignment
\[
\xymatrix@R-2pc{
\mathcal{P}(\Lambda) \ar[r] & Z^{1}(\Lambda,
H^{0}(V,\mathcal{O}^{\times}))\\
(H,\chi) \ar[r] & \AH_{(H,\chi)},
}
\]
is an injective group homomorphism, which after a composition with the
projection 
\[
Z^{1}(\Lambda, H^{0}(V,\mathcal{O}^{\times}) \twoheadrightarrow
H^{1}(\Lambda,H^{0}(V,\mathcal{O}^{\times}) \cong \op{Pic}(X)
\] 
becomes
an isomorphism.

We now have the following lemma which computes the
obstruction and ambiguity to extending $\phi_{0}$ to a non-commutative
cocycle $\phi = \phi_{0} + \hbar \phi_{1}$:

\begin{lem} \label{firststep}
\begin{itemize} 
\item[{\bf (a)}]
The obstruction map
\[
\op{ob}_{0}: H^{1}(\Lambda, H^{0}(V, \mathcal{O}^{\times})) \to 
H^{2}(\Lambda, H^{0}(V, \mathcal{O}))
\]
Can be lifted to a map on Appell-Humbert data:
\[
\xymatrix{
\mathcal{P}(\Lambda) \ar[r]^-{\mathfrak{ob}_{0}} \ar@{_{(}->}[d] &
Z^{2}(\Lambda,H^{0}(V, \mathcal{O})) \ar[d]\\
H^{1}(\Lambda, H^{0}(V, \mathcal{O}^{\times})) \ar[r]_-{\op{ob}_{0}} & 
H^{2}(\Lambda, H^{0}(V, \mathcal{O}))
}
\] 
where
\[
\mathfrak{ob}_{0}(\AH_{(H,\chi)})(\lambda_{1}, \lambda_{2}) = 
\{h_{\lambda_{2}}, h_{\lambda_{1}}\}.
\]
\item[{\bf (b)}]  Suppose that $W \in \op{Pic}(X)$ is such that
  $\op{ob}_{0}([W]) = 0$ in $H^{2}(X,\mathcal{O})$. Let $(H,\chi)$ be
  the Appell-Humbert data corresponding to $[W]$. Then
  $\mathfrak{ob}_{0}\left(\AH_{(H,\chi)}\right) = 0$ in
  $Z^{2}(\Lambda,H^{0}(V,\mathcal{O}))$.
\end{itemize}
\end{lem}
{\bf Proof.} Suppose that we can find a $\phi \in Z^{1}(\Lambda,
H^{0}(V,(\A/\hbar^{2}\A)^{\times})$, so that $\phi = \phi_{0}$ modulo
$\hbar$. If $\psi_{0}$ is a cocycle, cohomologous to $\phi_{0}$, then
we can find a cocycle $\psi \in Z^{1}(\Lambda,
H^{0}(V,(\A/\hbar^{2}\A)^{\times})$, so that $\psi = \psi_{0}$ modulo
$\hbar$, and $\psi$ is cohomologous to $\phi$. Indeed, if $f \in
H^{0}(V,\mathcal{O}^{\times})$ is a global holomorphic function for
which  $\psi_{0}(\lambda)/\phi_{0}(\lambda) = (f\cdot \lambda)/f$,
then by viewing $f$ as an element in $H^{0}(V,\A)$ we can define
a new cocycle $\psi$ according to the rule
\eqref{eq:coboundary} using $\phi$ and $f$. This new $\psi$ clearly
has the required properties. 

Hence, without a loss of generality, we may assume that $\phi =
\AH_{(H,\chi)} + \hbar \phi_{1}$ for some appropriately
chosen Appell-Humbert data $(H,\chi)$.

\medskip

\noindent
Let now  
$\delta: C^{1}(\Lambda, H^{0}(V, \mathcal{O})) \to 
C^{2}(\Lambda, H^{0}(V, \mathcal{O}))$ denote the group cohomology
differential given by 
\[
\theta \mapsto [(\lambda_{1}, \lambda_{2}) \mapsto 
\theta(\lambda_{1} + \lambda_{2}) - \theta(\lambda_{2}) - 
\theta(\lambda_{1}) \cdot \lambda_{2}].
\]
A non-commutative 
cochain $\phi = \phi_{0} + \hbar
\phi_{1} \in C^{1}(\Lambda,H^{0}(V,(\A/\hbar^{2}\A)^{\times})$ is a
cocycle 
if and only if   
$\phi_{1}' = \frac{\phi_{1}}{\phi_{0}} \in
C^{1}(\Lambda,H^{0}(V,\mathcal{O}))$  satisfies the condition
\[
(\delta \phi_{1}')(\lambda_{1}, \lambda_{2}) = 
\frac{(\phi_{0}(\lambda_{2})\star ((\phi_{0}(\lambda_{1}))\cdot 
\lambda_{2} ))_{1}}{(\phi_{0}(\lambda_{2})) 
(\phi_{0}(\lambda_{1}) \cdot \lambda_{2})}
=
\frac{(\phi_{0}(\lambda_{2})\star ((\phi_{0}(\lambda_{1}))\cdot 
\lambda_{2}))_{1} }{\phi_{0}(\lambda_{1} + \lambda_{2})}
\]
After substituting $\phi_{0} = \AH_{(H,\chi)}$ into this
formula, several terms cancel and we get:
\[
\begin{split}
(\delta \phi_{1}')(\lambda_{1}, \lambda_{2})  & = 
\frac{(\op{exp}(\pi H(v,\lambda_{2})) \star 
\op{exp}(\pi H(v+ \lambda_{2}, \lambda_{1})) )_{1}}
{\op{exp}(\pi H(v, \lambda_{2}) + \pi H(v+ \lambda_{2}, \lambda_{1}))}
\\
& \\
& =
\frac{(\op{exp}(\pi H(v,\lambda_{2})) \star 
\op{exp}(\pi H(v, \lambda_{1})) )_{1}}
{\op{exp}(\pi H(v, \lambda_{1}+ \lambda_{2}))} \\
& \\
& = \frac{(\exp(h_{\lambda_{2}})\star
  \exp(h_{\lambda_{1}}))_{1}}{\exp(h_{\lambda_{1} + \lambda_{2}})} \\
& \\
& = \left( \op{exp}\left(\hbar  \{ h_{\lambda_{2}}, h_{\lambda_{1}} \}  
\right) \right)_{1} \\
& \\
& = \{ h_{\lambda_{2}}, h_{\lambda_{1}} \}.
\end{split}
\]
Here $h_{\lambda}\in V^{\vee}$ denotes the
 $\mathbb{C}$-linear function 
$v \mapsto \pi H(v, \lambda)$ and in the last equality we used the identity
\eqref{eq-star-linear}.  Due to the equality $\delta \phi_{1}' =
\{ h_{\lambda_{2}}, h_{\lambda_{1}} \}$ we conclude that 
we will be able to extend $\AH_{(H,\chi)}$ to a
non-commutative cocycle
$\phi = \AH_{(H,\chi)} + \hbar \phi_{1}$ if and only if the
cocycle 
$[(\lambda_{1},\lambda_{2}) \mapsto \{ h_{\lambda_{2}},
  h_{\lambda_{1}} \}] \in 
Z^{2}(\Lambda,H^{0}(V,\mathcal{O})$ is a coboundary. This shows that 
$[(\lambda_{1},\lambda_{2}) \mapsto \{ h_{\lambda_{2}},
  h_{\lambda_{1}} \}]$ represents the 
obstruction class $\op{ob}_{0}([\AH_{(H,\chi)}])$ and proves
part {\bf (a)} of the lemma.

\

\noindent
For the proof of part {\bf (b)} notice that by construction the
cocycle $\mathfrak{ob}_{0}\left(\AH_{(H,\chi)}\right)$ is
actually in $Z^{2}(\Lambda,\mathbb{C}) \subset
Z^{2}(\Lambda,H^{0}(V,\mathcal{O}))$. Furthermore if we consider the
canonical Hodge decomposition 
\[
\begin{split}
H^{2}(X,\mathbb{C}) & =
H^{2}(X,\mathcal{O})\oplus H^{1}(X,\Omega^{1})\oplus H^{0}(X,
\Omega^{2}) \\
& = \wedge^{2}\overline{V}^{\vee} \oplus (\overline{V}^{\vee} \otimes
V^{\vee}) \oplus \wedge^{2} V^{\vee},
\end{split}
\]
then the image of
$\mathfrak{ob}_{0}\left(\AH_{(H,\chi)}\right)$ in
$H^{2}(X,\mathbb{C})$ lands entirely in the piece
$H^{2}(X,\mathcal{O}) = \wedge^{2}\overline{V}^{\vee}$. Indeed,
thinking of $H$ as an element in $\overline{V}^{\vee} \otimes
V^{\vee}$ we can rewrite the image
$\mathfrak{ob}_{0}\left(\AH_{(H,\chi)}\right)$ in purely
linear algebraic terms as the contraction $H 
\cntrct \bpi \cntrctother H$. Indeed, the additive map
$\lambda_{1}\wedge \lambda_{2} \mapsto \{h_{\lambda_{2}},
h_{\lambda_{1}} \}$ extends by linearity to a unique conjugate linear
homomorphism $\wedge^{2}V \to \mathbb{C}$ which equals
$H \cntrct \bpi \cntrctother H$ as an element in
$\wedge^{2}\overline{V}^{\vee}$. 

Therefore the obstruction $\op{ob}_{0}([\AH_{(H,\chi)}])$
vanishes if and only if $H \cntrct \bpi \cntrctother H = 0 $ in
$H^{2}(X, \mathbb{C})$. Since $H \cntrct \bpi \cntrctother H$ was
proportional to the anti-linear extension of the map \linebreak
$\mathfrak{ob}(\AH_{(H,\chi)}) : \Lambda \times \Lambda \to
\mathbb{C}$ this concludes the proof of part {\bf (b)}. \ \hfill $\Box$

\

\medskip

\noindent
The formula for $\op{ob}_{0}$ given in the lemma can also be deduced
from the first order analysis carried out in Toda's paper
\cite{Toda}. However, in our case, the specific geometry of the Moyal
quantization of a complex torus allows us to push the analysis
further. In fact, it turns out that for a line bundle $W$ on $X$,
the vanishing of $\op{ob}_{0}([W])$ is both necessary and sufficient for 
$W$ to quantize to all orders:

\begin{lem} \label{lem:quantizable}
A line bundle $W$ on $X$ can be extended to a line
  bundle on $\mathbb{X}_{\bpi}$ if and only if 
\[
\op{ob}_{0}([W]) = c_{1}(W)\cntrct \bpi \cntrctother c_{1}(W) = 0
\quad \text{in} \quad H^{2}(X,\mathcal{O}).
\]
\end{lem}
{\bf Proof.} Let $(H,\chi)$ be the Appell-Humbert data for the
isomorphism class of line bundles $[W]$. Consider the map
\[
\xymatrix@R-2pc{ \Lambda \ar[r]^{\phi} & H^{0}(V,\Ah{V,\bpi}^{\times}) \\
\lambda \ar[r] & \left[ v \mapsto \chi(\lambda)\op{exp}\left(\pi H(v,
\lambda) + \frac{\pi}{2}H(\lambda, \lambda) \right)\right] }
\]
where the exponential now is the $\star$-exponential. By definition 
$\phi \in C^{1}(\Lambda,H^{0}(V,\Ah{V,\bpi}^{\times}))$ is a
non-commutative cochain. In order for $\phi$ to be a cocycle, we must
have that 
\begin{equation} \label{eq:defect}
\phi(\lambda_{1}+\lambda_{2})^{-1}\star \phi(\lambda_{2}) \star
(\phi(\lambda_{1})\cdot \lambda_{2}) = 1.
\end{equation}
However, in  the proof of Lemma~\ref{firststep} we evaluated the left
hand side of \eqref{eq:defect} and showed that it is equal to 
$\exp(\hbar\{h_{\lambda_{2}}, h_{\lambda_{1}} \})$. This
proves our assertion since the constant $\{h_{\lambda_{2}},
h_{\lambda_{1}} \}$ is equal to the value of $H \cntrct
\bpi \cntrctother H \in \wedge^{2}\overline{V}^{\vee}$ on the element
$\lambda_{1}\wedge \lambda_{2}$. \ \hfill $\Box$

\

\medskip

\noindent
Suppose $(H,\chi) \in \mathcal{P}(\Lambda)$ is some
  Appell-Humbert data and let $l(\hbar) = \sum_{i=1}^{\infty}
  \hbar^{i}l_{i} \in \hbar
  \overline{V}^{\vee}[[\hbar]]$. Consider 
  the map $\qAH_{((H,\chi),l(\hbar))} : \Lambda \to
  H^{0}(V,\Ah{V,\bpi}^{\times})$ given by 
\begin{equation} \label{eq:qAH}
\qAH_{((H,\chi),\, l(\hbar))}(\lambda)(v) = \chi(\lambda)
\op{exp}\left(\pi H(v, \lambda) +  
\frac{\pi}{2}H(\lambda, \lambda) + 
\sum_{j=1}^{\infty} \hbar^{j}\pi \langle l_{j},\lambda \rangle \right).
\end{equation}
A straightforward check shows that if $(H,\chi)$ satisfies $H \cntrct
\bpi \cntrctother H = 0$, then the map \linebreak
 $\qAH_{((H,\chi),\, l(\hbar))}$ is
a non-commutative 1-cocycle, i.e. $\qAH_{((H,\chi),\, l(\hbar))} \in
Z^{1}(\Lambda, H^{0}(V,\A^{\times}))$.

\
\medskip

\noindent
Consider now the subset $\mathcal{P}(\Lambda,\bpi) \subset
\mathcal{P}(\Lambda)$ defined by 
\[
\mathcal{P}(\Lambda,\bpi) = \left\{ (H,\chi) \in  \mathcal{P}(\Lambda) 
\left| H \cntrct
\bpi \cntrctother H = 0 \right. \right\}.
\]
With this notation we have the following quantum version of the
Appell-Humbert theorem:
\begin{prop} The map 
\[
\xymatrix@R-2pc{
\mathcal{P}(\Lambda,\bpi)\times \hbar \overline{V}^{\vee}[[\hbar]]
\ar[r] & H^{1}(X,\Ah{X,\bpi}^{\times}) \\
\left( (H,\chi), \sum_{i = 1}^{\infty} \hbar^{i}l_{i} \right)
\ar@{|->}[r] & \qAH_{((H,\chi),\, l(\hbar))}(\lambda)
}
\]
is a bijection of pointed sets. 
\end{prop}
{\bf Proof.}  It suffices to show that for each $j$ the map $\qAH$
induces a bijection
\[
\mathcal{P}(\Lambda,\bpi)\times \hbar
\overline{V}^{\vee}[[\hbar]]/\hbar^{j} \to  
H^{1}(X,(\Ah{X,\bpi}/\hbar^{j})^{\times}).
\]
Note that the case $j=1$ is the usual Appell-Humbert theorem.  If we
assume that this has been shown for $j \leq n$ then to show that it
holds for $j=n+1$ we can use the argument from the proof of 
Lemma~\ref{firststep} where it was shown that the case $j=1$
implies the case $j=2$.  The key point is that for $\psi \in
C^{1}(\Lambda, H^{0}(V, \mathcal{O}))$ the cochain in $C^{1}(\Lambda,
H^{0}(V,(\Ah{\bpi}/\hbar^{n+1})^{\times})$ given by
\[\lambda \mapsto \chi(\lambda)
\op{exp}\left(\pi H(v, \lambda) +  
\frac{\pi}{2}H(\lambda, \lambda) + 
\sum_{j=1}^{n} \hbar^{j}\pi \langle l_{j},\lambda \rangle +
\hbar^{n+1} \psi \right) 
\]
is a cocycle if and only if $\psi$ is in $Z^{1}(\Lambda, H^{0}(V,
\mathcal{O}))$. 
\hfill $\Box$

\

\medskip
 
\begin{rem} \label{rem:qAHd0} 
The quantum Appell-Humbert theorem gives preferred group cocycle
representatives for the isomorphism classes of quantum line bundles.
For degree zero line bundles the representatives are given by
\[\lambda \mapsto \chi(\lambda) \op{exp}\left(
\sum_{j=1}^{\infty} \hbar^{j} \pi \langle l_{j},\lambda \rangle \right).
\]
Here $\chi \in \op{Hom}(\Lambda, U(1))$ and $l_{j} \in
\overline{V}^{\times}$.   Thus  the connected component 
$H^{1}(X, \Ah{X,\bpi}^{\times})_{o}$ 
of the quantum Picard $H^{1}(X, \Ah{X,\bpi}^{\times})$ 
is in bijection with 
$X^{\vee} \times (\overline{V}^{\vee})^{\mathbb{Z}_{>0}} 
\cong H^{1}(X, \mathcal{O}[[\hbar]]^{\times} )_{o}$.  
\end{rem}

\begin{rem}
As an example consider the product of elliptic curves $E_{1}\times
E_{2}$ with coordinates $(z_{1},z_{2})$ and Poisson structure
$\frac{\partial}{\partial z_{1}} \wedge \frac{\partial}{\partial
z_{2}}$.  Let $L$ be the line bundle corresponding to the divisor
$E_{1} \times \{ 0 \}$ and $M$ the line bundle corresponding to the
divisor $\{ 0 \} \times E_{2}$.  Then $L$ and $M$ are quantizable but
$L \otimes M$ is not.  Notice also that $L \oplus M$ is a quantizable
vector bundle yet has a non-zero second Chern class, given by the
Poincar\'{e} dual to the intersection of the two divisors.
\end{rem}

\

\bigskip
\bigskip

\Appendix{On the cohomology of left $\A$-modules}
\label{ss:cohomologyleft} 

\

\noindent
Let $S$ be a Hausdorff analytic space and $\A$ be a sheaf of 
$\mathbb{C}[[\hbar]]-$algebras such that 
$\A\otimes_{\mathbb{C}[[\hbar]]} \mathbb{C} \cong \mathcal{O}_{S}$.  
The usual arguments in \cite{Hartshorne} go through to show that if
$\mycal{L}$ is a left $\A$-module, then there are well-defined
cohomology groups of $\mycal{L}$ computed from the derived functors
of the global sections functor in the categories of sheaves of abelian
groups, sheaves of $\mathbb{C}$-vector spaces, sheaves of
$\Ch$-modules, and sheaves of $\A$-modules.  Furthermore, all these
cohomologies are naturally isomorphic to each other.  If $\mycal{L}$
is a locally free left $\A$-module of finite rank, then these
cohomologies also agree with the \u{C}ech cohomology of $\mycal{L}$.
Indeed, chose a cover 
$\mathfrak{U} = \{ U_{i} | i \in \mathcal{I} \}$ by contractible open 
sets such that for $U_{I} = \cap_{j \in I} U_{j}$ we have  
$H^{j}(S, U_{I},\mathcal{O}) = \{ 0 \}$ and
$H^{1}(S, U_{I},\mathcal{O}^{\times}) = \{ 1 \}$ for all finite subsets
$I \subset \mathcal{I}$ and all $j \geq 1$.  Denote the inclusion maps by 
$\kappa_{I}: U_{I} \to S$.  Then, by
Lemma~\ref{lem:acyclic} we have $H^{j}(S,U_{I}, \mycal{L}) \cong
H^{j}(S,U_{I}, \A) \cong H^{j}(S,U_{I}, \mathcal{O})[[\hbar]] = \{ 0 \}$
for all $I$ and $j \geq 1$.  Therefore we can compute cohomology from 
the following acyclic resolution
\[
0 \to \mycal{L} \to \bigoplus_{i \in \mathcal{I}}
\kappa_{i*}(\mycal{L}|_{U_i})
\to \bigoplus_{\{i,j\} \in \mathcal{I}^{\times 2}}
\kappa_{ij*}(\mycal{L}|_{U_{i,j}})
\to \cdots
\]
This cohomology is precisely the \u{C}ech cohomology.

\begin{lem} \label{lem:cohomology} Let $\mycal{L}$ be a degree zero
  line bundle on the non-commutative torus $\bbX_{\bpi}$. View
  $\mycal{L}$ as a sheaf of left $\Ah{X,\bpi}$ modules on the
  underlying torus $X$. Then
\begin{itemize}
\item[{\bf (a)}] $\mycal{L}$ is non-trivial if and only if
  $H^{0}(X,\mycal{L}) = 0$.
\item[{\bf (b)}] If $\mycal{L}/\hbar\mycal{L} \not\cong \mathcal{O}$, then
  $H^{i}(X,\mycal{L}) = 0$ for all $i \geq 0$.
\end{itemize}
\end{lem}
{\bf Proof.} First we prove {\bf (a)}. Since $\mycal{L}$ is a
translation invariant line bundle on $\bbX_{\bpi}$, it is given by a
constant factor of automorphy (see Remark~\ref{rem:qAHd0}) and so the
sheaf $\mycal{L}$ has a preferred flat connection. Denote the
corresponding local system of $\Ch$-modules by
$\underline{\mycal{L}}$. Furthermore the sheaf $\mycal{L}$ has a
natural structure of an $\mathcal{O}[[\hbar]]$-module. This follows from
the identification $\mycal{L} \cong \underline{\mycal{L}}\otimes_{\mathbb{C}}
\mathcal{O} = \underline{\mycal{L}}\otimes_{\Ch} \mathcal{O}[[\hbar]]
$. Now observe that  
\[
H^{0}(X,\mycal{L})  =
\op{Hom}_{\mathcal{O}-\op{mod}}(\mathcal{O},\mycal{L}) =
  \lim_{\longleftarrow}\displaylimits
  \op{Hom}_{\mathcal{O}-\op{mod}}(\mathcal{O},\mycal{L}/\hbar^{k}) =
  \lim_{\longleftarrow}\displaylimits H^{0}(X,\mycal{L}/\hbar^{k}). 
\]
Note that $\mycal{L}/\hbar^{k}$ is a free
$\mathcal{O}[[\hbar]]/\hbar^{k}$ module and so is a holomorphic vector
bundle of rank $k$ on $X$. The vector bundle $\mycal{L}/\hbar^{k}$ has
a preferred flat connection coming from the natural identification
$\mycal{L}/\hbar^{k} \cong
(\underline{\mycal{L}}/\hbar^{k})\otimes_{\mathbb{C}}
\mathcal{O}$. From the formula for the factor of automorphy for
$\mycal{L}$ (see Remark~\ref{rem:qAHd0}) it is clear that the flat
connection on $\mycal{L}/\hbar^{k}$ is unitary. Since $X$ is
K\"{a}hler we can now use that the Hodge decomposition to compare the
cohomology of the local system $\underline{\mycal{L}}/\hbar^{k}$ and
the holomorphic bundle $\mycal{L}/\hbar^{k}$. In particular we have
that the natural map $\underline{\mycal{L}}/\hbar^{k} \to
\mycal{L}/\hbar^{k}$ induces an isomorphism on $H^{0}$. Therefore 
\[
H^{0}(X,\mycal{L}) = \lim_{\longleftarrow}\displaylimits
H^{0}(X,\mycal{L}/\hbar^{k}) = \lim_{\longleftarrow}\displaylimits
H^{0}(X,\underline{\mycal{L}}/\hbar^{k}) = H^{0}(X,\underline{\mycal{L}}),
\]
and so a non-zero global section of $\mycal{L}$ is nowhere
vanishing. Since $\mycal{L}$ is a locally free rank one
$\Ah{X,\bpi}$-module this implies that $\mycal{L}$ is trivial. The
opposite implication is obvious. This completes the proof of part {\bf (a)}.

\

\medskip

\noindent
For part {\bf (b)} recall the classical result that if $L$ is a
non-trivial degree zero line bundle on $X$, then $H^{j}(L) = 0$ for
all $j\geq 0$. Suppose now 
$\mycal{L}$ is a quantum line bundle of degree zero for which
$\mycal{L}/\hbar \mycal{L}$ is non-trivial.  Consider the short exact
sequence 
\[
0 \rar{} \mycal{L} \rar{\hbar} \mycal{L}
\rar{} \mycal{L}/\hbar \mycal{L} \rar{} 0.
\]
Since all cohomology groups  of $\mycal{L}/\hbar \mycal{L}$ vanish,
the long exact
cohomology sequence implies 
that $\hbar$ induces an isomorphism of
$\mathbb{C}[[\hbar]]$-modules $H^j(\mycal{L}) \to H^j(\mycal{L})$
for all $j$.
This implies that $H^j(\mycal{L}) = 0$ for all $j$.  Indeed
if not, chose the largest possible $p \in \{0,1,2, \dots\}$ such that
$H^j(\mycal{L})/\hbar^{p} H^j(\mycal{L}) = (0)$, then applying the
isomorphism gives that
$H^j(\mycal{L})/\hbar^{p + 1} H^j(\mycal{L}) = 0$,
a contradiction. \ \hfill $\Box$

\

\medskip

\noindent Notice that if $\mycal{L}$ is a degree zero quantum line
bundle with $\mycal{L}/\hbar \mycal{L} \cong \mathcal{O}$ $\mycal{L}$
is non-trivial.  In this case the higher cohomology groups of
$\mycal{L}$ need not vanish.  The easiest way to see this is to note that
\begin{equation} \label{eq:hL}
0 \rar{} \mycal{L} \rar{\hbar} \mycal{L}
\rar{} \mathcal{O} \rar{} 0.
\end{equation}
is a short exact sequence of $\mathcal{O}$-modules and so the
extension class of this sequence
lies in $H^1(X,\mycal{L})$.  Modulo $\hbar^{2}$ the sequence
\eqref{eq:hL} induces a
short exact sequence of $\mathcal{O}$-modules
\begin{equation} \label{eq:hhL}
\xymatrix@R-1.3pc{
0 \ar[r] &  \mycal{L}/\hbar \ar@{=}[d] \ar[r] & \mycal{L}/\hbar^{2}
\ar[r] & \mycal{L}/\hbar \ar@{=}[d] \ar[r] & 0 \\
& \mathcal{O} & & \mathcal{O}
}
\end{equation}
whose extension class is in $H^{1}(X,\mathcal{O})$. 
If $H^1(X,\mycal{L})$ was zero, then this sequence \eqref{eq:hL}, and
hence the sequence \eqref{eq:hhL},
will split. However it is immediate to check that if $\mycal{L}$ is
represented by quantum Appell-Humbert data $\left((0,1), \sum_{i =
  1}^{\infty} \hbar^{i} l_{i} \right)$, then the extension class of
\eqref{eq:hhL} is given by the group cohomology class $[l_{1}] \in
H^{1}(\Lambda, H^{0}(V,\mathcal{O}))$. Since we are completely free to
chose the $l_{i}$'s this shows that \eqref{eq:hhL}, and hence
\eqref{eq:hL}, are non-split in general. This implies that a
general $\mycal{L}$ with $\mycal{L}/\hbar \cong \mathcal{O}$ has
non-trivial first cohomology.


\newcommand{\etalchar}[1]{$^{#1}$}

\

\bigskip

\noindent
Department of Mathematics\\
University of Pennsylvania \\
209 South 33rd Street \\
Philadelphia, PA 19104-6395

\

\bigskip

\noindent
email: orenb@math.upenn.edu, blockj@math.upenn.edu, tpantev@math.upenn.edu

\end{document}